\input amstex
\documentstyle{amsppt}
%
\catcode`@=11
\redefine\output@{%
  \def\break{\penalty-\@M}\let\par\endgraf
  \ifodd\pageno\global\hoffset=105pt\else\global\hoffset=8pt\fi  
  \shipout\vbox{%
    \ifplain@
      \let\makeheadline\relax \let\makefootline\relax
    \else
      \iffirstpage@ \global\firstpage@false
        \let\rightheadline\frheadline
        \let\leftheadline\flheadline
      \else
        \ifrunheads@ 
        \else \let\makeheadline\relax
        \fi
      \fi
    \fi
    \makeheadline \pagebody \makefootline}%
  \advancepageno \ifnum\outputpenalty>-\@MM\else\dosupereject\fi
}
\catcode`\@=\active
\nopagenumbers
\def\negskp{\hskip -2pt}
\def\MatGrSL{\operatorname{SL}}
\def\MatGrSO{\operatorname{SO}}
\def\Alpha{\operatorname{A}}
\def\idop{\operatorname{\bold{id}}}
\def\tr{\operatorname{tr}}
\def\compos{\,\raise 1pt\hbox{$\sssize\circ$} \,}
\def\vtttrule{\vrule height 12pt depth 19pt}
\def\boxit#1#2{\vcenter{\hsize=122pt\offinterlineskip\hrule
  \line{\vtttrule\hss\vtop{\hsize=120pt\centerline{#1}\vskip 5pt
  \centerline{#2}}\hss\vtttrule}\hrule}}
\def\msum#1{\operatornamewithlimits{\sum^#1\!{\ssize\ldots}\!\sum^#1}}
\accentedsymbol\tx{\tilde x}
\accentedsymbol\bd{\kern 2pt\bar{\kern -2pt d}}
\accentedsymbol\bbd{\kern 2pt\bar{\kern -2pt\bold d}}
\accentedsymbol\bPsi{\kern 1pt\overline{\kern -1pt\boldsymbol\Psi
\kern -1pt}\kern 1pt}
\accentedsymbol\bulletH{\overset{\kern 2pt\sssize\bullet}\to H}
\accentedsymbol\circH{\overset{\kern 2pt\sssize\circ}\to H}
\accentedsymbol\bulletBH{\overset{\sssize\bullet}\to{\bold H}}
\accentedsymbol\circBH{\overset{\sssize\circ}\to{\bold H}}
\accentedsymbol\bulletd{\overset{\kern 3pt\sssize\bullet\kern -2pt}\to d}
\accentedsymbol\circd{\overset{\kern 3pt\sssize\circ\kern -2pt}\to d}
\accentedsymbol\bulbulgamma{\overset{\kern 1pt\sssize\bullet%
\bullet}\to\gamma}
\accentedsymbol\bulcircgamma{\overset{\kern 1pt\sssize\bullet%
\circ}\to\gamma}
\accentedsymbol\circbulgamma{\overset{\kern 1pt\sssize\circ%
\bullet}\to\gamma}
\accentedsymbol\circcircgamma{\overset{\kern 1pt\sssize\circ%
\circ}\to\gamma}
\accentedsymbol\bulbulA{\overset{\sssize\bullet%
\bullet}\to{\operatorname{A}}}
\accentedsymbol\bulcircA{\overset{\sssize\bullet%
\circ}\to{\operatorname{A}}}
\accentedsymbol\circbulA{\overset{\sssize\circ%
\bullet}\to{\operatorname{A}}}
\accentedsymbol\circcircA{\overset{\sssize\circ%
\circ}\to{\operatorname{A}}}
\def\blue#1{#1}
\catcode`#=11\def\diez{#}\catcode`#=6
\catcode`_=11\def\podcherkivanie{_}\catcode`_=8
\def\mycite#1{\cite{\blue{#1}}\immediate\special{ps:
     ShrHPSdict begin /ShrBORDERthickness 0 def}}

\def\mytag#1{%
    \tag#1}
\def\mythetag#1{\thetag{\blue{#1}}\immediate\special{ps:
     ShrHPSdict begin /ShrBORDERthickness 0 def}}
\def\myrefno#1{\no#1}
\def\myhref#1#2{\blue{#2}\immediate\special{ps:
     ShrHPSdict begin /ShrBORDERthickness 0 def}}
\def\myEarXivlink{\myhref{http://arXiv.org}{http:/\negskp/arXiv.org}}
\def\myGeoCities{\myhref{http://www.geocities.com}{GeoCities}}
\def\mytheorem#1{\csname proclaim\endcsname{Theorem #1}}
\def\mythetheorem#1{\blue{#1}\immediate\special{ps:
     ShrHPSdict begin /ShrBORDERthickness 0 def}}
\def\mylemma#1{\csname proclaim\endcsname{Lemma #1}}

\def\mycorollary#1{\csname proclaim\endcsname{Corollary #1}}
\def\mythecorollary#1{\blue{#1}\immediate\special{ps:
     ShrHPSdict begin /ShrBORDERthickness 0 def}}
\def\mydefinition#1{\definition{Definition #1}}
\def\mythedefinition#1{\blue{#1}\immediate\special{ps:
     ShrHPSdict begin /ShrBORDERthickness 0 def}}

\pagewidth{360pt}
\pageheight{606pt}
\topmatter
\title
A note on metric connections for chiral and Dirac spinors.
\endtitle
\author
R.~A.~Sharipov
\endauthor
\address 5 Rabochaya street, 450003 Ufa, Russia\newline
\vphantom{a}\kern 12pt Cell Phone: +7-(917)-476-93-48
\endaddress
\email \vtop to 30pt{\hsize=280pt\noindent
\myhref{mailto:R\podcherkivanie Sharipov\@ic.bashedu.ru}
{R\_\hskip 1pt Sharipov\@ic.bashedu.ru}\newline
\myhref{mailto:r-sharipov\@mail.ru}
{r-sharipov\@mail.ru}\newline
\myhref{mailto:ra\podcherkivanie sharipov\@lycos.com}{ra\_\hskip 1pt
sharipov\@lycos.com}\vss}
\endemail
\urladdr
\vtop to 20pt{\hsize=280pt\noindent
\myhref{http://www.geocities.com/r-sharipov}
{http:/\negskp/www.geocities.com/r-sharipov}\newline
\myhref{http://www.freetextbooks.boom.ru/index.html}
{http:/\negskp/www.freetextbooks.boom.ru/index.html}\vss}
\endurladdr
\abstract
    It is known that the bundle of Dirac spinors is produced 
as a direct sum of two bundles --- the bundle of chiral spinors 
and its Hermitian conjugate bundle. In this paper some aspects 
of metric connections for chiral and Dirac spinors are resumed 
and their relation is studied.
\endabstract
\endtopmatter
\loadbold
\loadeufb
\TagsOnRight
\document
\accentedsymbol\tbvartheta{\tilde{\overline{\boldsymbol\vartheta}}
\vphantom{\boldsymbol\vartheta}}

\rightheadtext{A note on metric connections \dots}
\head
1. Chiral spinors. 
\endhead
    The construction of two-component Weyl spinors, they are also called 
chiral spinors, is based on the following  well-known group homomorphism
$$
\hskip -2em
\varphi\!:\MatGrSL(2,\Bbb C)\to\MatGrSO^+(1,3,\Bbb R).
\mytag{1.1}
$$
The homomorphism \mythetag{1.1} is defined through the formula
$$
\hskip -2em
\goth S\cdot\boldsymbol\sigma_m\cdot\goth S^{\sssize\dagger}=\sum^3_{k=0}
S^k_m\,\boldsymbol\sigma_k,
\mytag{1.2}
$$
where $\boldsymbol\sigma_1,\,\boldsymbol\sigma_2,\,\boldsymbol\sigma_3$ 
are Pauli matrices complemented with the unit matrix $\boldsymbol\sigma_0$:
$$
\xalignat 2
&\hskip -2em
\boldsymbol\sigma_0=\Vmatrix 1 & 0\\0 & 1\endVmatrix,
&&\boldsymbol\sigma_2=\Vmatrix 0 & -i\\i & 0\endVmatrix,\\
\vspace{-1.4ex}
&&&\mytag{1.3}\\
\vspace{-1.4ex}
&\hskip -2em
\boldsymbol\sigma_1=\Vmatrix 0 & 1\\1 & 0\endVmatrix,
&&\boldsymbol\sigma_3=\Vmatrix 1 & 0\\0 & -1\endVmatrix.
\endxalignat
$$
By means of \mythetag{1.2} and \mythetag{1.3} each matrix $\goth S\in
\MatGrSL(2,\Bbb C)$ is associated with some matrix $S\in\MatGrSO^+(1,3,
\Bbb R)$ so that we can write $S=\varphi(\goth S)$, see \mycite{1},
\mycite{2}, and \mycite{3} for detailed description of this
construction.\par
    Let $M$ be a {\it space-time} manifold, i\.\,e\. a four-dimensional
orientable manifold equipped with a pseudo-Euclidean Minkowski-type 
metric $\bold g$ and carrying a special smooth geometric structure 
which is called a {\it polarization}. Once some polarization is fixed, 
one can distinguish the {\it Future light cone\/} from the {\it Past 
light cone\/} at each point $p\in M$ (see \mycite{4} for more details). 
A moving frame $(U,\,\boldsymbol\Upsilon_0,\,\boldsymbol\Upsilon_1,\,
\boldsymbol\Upsilon_2,\,\boldsymbol\Upsilon_3)$ of the tangent bundle 
$TM$ is an ordered set of four smooth vector fields $\boldsymbol
\Upsilon_0$, $\boldsymbol\Upsilon_1$, $\boldsymbol\Upsilon_2$, 
$\boldsymbol\Upsilon_3$ which are defined and $\Bbb R$-linearly 
independent at each point $p$ of some open subset $U\subset M$. This 
moving frame is called a {\it positively polarized right orthonormal 
frame} if the following conditions are fulfilled:
\roster
\rosteritemwd=1pt
\item the value of the first vector filed $\boldsymbol\Upsilon_0$ at each
point $p\in U$ belongs to the interior of the Future light cone determined 
by the polarization of $M$;
\item it is a right frame in the sense of the orientation of $M$;
\item the metric tensor $\bold g$ is given by the standard
Minkowski matrix in this frame:
$$
g_{ij}=g(\boldsymbol\Upsilon_i,\boldsymbol\Upsilon_j)=\Vmatrix
\format \l&\quad\r&\quad\r&\quad\r\\
1 &0 &0 &0\\ 0 &-1 &0 &0\\ 0 &0 &-1 &0\\ 0 &0 &0 &-1\endVmatrix.
$$
\endroster
Apart from positively polarized right orthonormal frames, below we shall
consider the following three special types of frames in $TM$:
\rosteritemwd=5pt
\roster
\item"---" positively polarized left orthonormal frames;
\item"---" negatively polarized right orthonormal frames;
\item"---" negatively polarized left orthonormal frames.
\endroster
The definitions of these types of frames are easily obtained by 
alternating the above condition \therosteritem{1} and \therosteritem{2}
with the opposite ones.\par
     Let $(U,\,\boldsymbol\Upsilon_0,\,\boldsymbol\Upsilon_1,\,\boldsymbol
\Upsilon_2,\,\boldsymbol\Upsilon_3)$ and $(\tilde U,\,\tilde{\boldsymbol
\Upsilon}_0,\,\tilde{\boldsymbol\Upsilon}_1,\,\tilde{\boldsymbol
\Upsilon}_2,\,\tilde{\boldsymbol\Upsilon}_3)$ be arbitrary two frames of
the tangent bundle $TM$ such that $U\cap\tilde U\neq\varnothing$. Then at
each point $p\in U\cap\tilde U$ we can write the following relationships 
for their frame vectors:
$$
\xalignat 2
&\hskip -2em
\tilde{\boldsymbol\Upsilon}_i=\sum^3_{j=0}S^j_i\,\boldsymbol\Upsilon_j,
&&\boldsymbol\Upsilon_i=\sum^3_{j=0}T^j_i\,\tilde{\boldsymbol\Upsilon}_j.
\mytag{1.4}
\endxalignat
$$
The relationships \mythetag{1.4} are called {\it transition formulas},
while the coefficients $S^j_i$ and $T^j_i$ in them are the components
of two mutually inverse transition matrices $S$ and $T=S^{-1}$. If both
frames $(U,\,\boldsymbol\Upsilon_0,\,\boldsymbol\Upsilon_1,\,\boldsymbol
\Upsilon_2,\,\boldsymbol\Upsilon_3)$ and $(\tilde U,\,\tilde{\boldsymbol
\Upsilon}_0,\,\tilde{\boldsymbol\Upsilon}_1,\,\tilde{\boldsymbol
\Upsilon}_2,\,\tilde{\boldsymbol\Upsilon}_3)$ are positively polarized
right orthonormal frames, then the transition matrices $S$ and $T$ both 
are orthochronous Lorentzian matrices with $\det S=1$ and $\det T=1$. 
Such matrices form the {\it special orthochronous matrix Lorentz group}
$\MatGrSO^+(1,3,\Bbb R)$.\par
    Assume that $SM$ is a two-dimensional smooth complex vector bundle 
over the space-time $M$ equipped with a non-vanishing skew-symmetric
bilinear form $\bold d$ at each point $p\in M$. This bilinear form 
$\bold d$ is called the {\it spin-metric tensor}. A moving frame $(U,
\,\boldsymbol\Psi_1,\,\boldsymbol\Psi_2)$ of $SM$ is an ordered set 
of two smooth sections $\boldsymbol\Psi_1$ and $\boldsymbol\Psi_2$ over
some open subset $U\subset M$ which are $\Bbb C$-linearly independent 
at each point $p\in U$. A moving frame $(U,\,\boldsymbol\Psi_1,\,
\boldsymbol\Psi_2)$ is called an {\it orthonormal frame\/} if
$$
\hskip -2em
d_{ij}=d(\boldsymbol\Psi_i,\boldsymbol\Psi_j)
=\Vmatrix 0 & 1\\ 
\vspace{1ex} -1 & 0\endVmatrix,
\mytag{1.5}
$$
i\.\,e\. if the spin-metric tensor $\bold d$ is given by the skew-symmetric
matrix \mythetag{1.5} in this frame. For two arbitrary frames $(U,\,
\boldsymbol\Psi_1,\,\boldsymbol\Psi_2)$ and $(\tilde U,\,\tilde{\boldsymbol
\Psi}_1,\,\tilde{\boldsymbol\Psi}_2)$ of the bundle $SM$ with overlapping 
domains $U\cap\tilde U\neq\varnothing$ we can write the following
transition formulas:
$$
\xalignat 2
&\hskip -2em
\tilde{\boldsymbol\Psi}_i=\sum^2_{j=1}\goth S^j_i\,
\boldsymbol\Psi_j,
&&\boldsymbol\Psi_i=\sum^2_{j=1}\goth T^j_i\,
\tilde{\boldsymbol\Psi}_j.
\mytag{1.6}
\endxalignat
$$
Like in \mythetag{1.4}, the coefficients $\goth S^j_i$ and $\goth T^j_i$
in \mythetag{1.6} are the components of two mutually inverse transition
matrices $\goth S$ and $\goth T=\goth S^{-1}$. If $(U,\,\boldsymbol
\Psi_1,\,\boldsymbol\Psi_2)$ and $(\tilde U,\,\tilde{\boldsymbol\Psi}_1,
\,\tilde{\boldsymbol\Psi}_2)$ are two orthonormal frames, then both
matrices $\goth S$ and $\goth T$ in \mythetag{1.6} belong to the special
linear matrix group $\MatGrSL(2,\Bbb C)$. 
\mydefinition{1.1} A two-dimensional complex vector 
bundle $SM$ over the space-time manifold $M$ equipped with a nonzero
spin-metric $\bold d$ is called a {\it spinor bundle} if each
orthonormal frame $(U,\,\boldsymbol\Psi_1,\,\boldsymbol\Psi_2)$ of
$SM$ is associated with some positively polarized right orthonormal 
frame $(U,\,\boldsymbol\Upsilon_0,\,\boldsymbol\Upsilon_1,\,
\boldsymbol\Upsilon_2,\,\boldsymbol\Upsilon_3)$ of the tangent bundle
$TM$ such that for any two orthonormal frames $(U,\,\boldsymbol\Psi_1,
\,\boldsymbol\Psi_2)$ and $(\tilde U,\,\tilde{\boldsymbol\Psi}_1,\,
\tilde{\boldsymbol\Psi}_2)$ with overlapping domains $U\cap\tilde U\neq
\varnothing$ the associated tangent frames $(U,\,\boldsymbol\Upsilon_0,
\,\boldsymbol\Upsilon_1,\,\boldsymbol\Upsilon_2,\,\boldsymbol\Upsilon_3)$
and $(\tilde U,\,\tilde{\boldsymbol\Upsilon}_0,\,
\tilde{\boldsymbol\Upsilon}_1,\,\tilde{\boldsymbol\Upsilon}_2,\,
\tilde{\boldsymbol\Upsilon}_3)$ are related to each other by means of
the formulas \mythetag{1.4}, where the transition matrices $S$ and $T$
are obtained from the transition matrices $\goth S$ and $\goth T$ in
\mythetag{1.6} by applying the homomorphism \mythetag{1.1}, i\.\,e\. 
$S=\varphi(\goth S)$ and $T=\varphi(\goth T)$.
\enddefinition
    The definition~\mythedefinition{1.1} reflects the basic feature of
all spinor bundles. They are closely related to tangent bundle and
this relation is implemented through associated frame pairs of some
definite types.
\head
2. Tensorial and spin-tensorial fields.
\endhead
    Tensorial an spin-tensorial fields are introduced in a standard way
as described in \mycite{3}, \mycite{5}, \mycite{6}, and many other papers.
First of all we introduce the complexified tangent and cotangent bundles
$\Bbb CTM$ and $\Bbb CT^*\!M$:
$$
\xalignat 2
&\hskip -2em
\Bbb CT_p(M)=\Bbb C\otimes T_p(M),
&&\Bbb CT^*_p(M)=\Bbb C\otimes T^*_p(M).
\mytag{2.1}
\endxalignat
$$
The complex bundles \mythetag{2.1} are obviously dual to each other.
Then we introduce the conjugate and Hermitian conjugate bundles for the
spinor bundle $SM$:
$$
\xalignat 3
&\hskip -2em
S^*\!M,
&&S^{\sssize\dagger}\!M,
&&S^{*\sssize\dagger}\!M=S^{{\sssize\dagger}*}\!M.
\mytag{2.2}
\endxalignat
$$
Using \mythetag{2.1} and \mythetag{2.2}, we define the following
tensor products:
$$
\gather
\hskip -2em
\Bbb CT^m_n\!M=\overbrace{\Bbb CTM\otimes\ldots\otimes
\Bbb CTM}^{\text{$m$ times}}\otimes\underbrace{\Bbb CT^*\!M
\otimes\ldots\otimes \Bbb CT^*\!M}_{\text{$n$ times}}\,,
\mytag{2.3}\\
\hskip -2em
S^\alpha_\beta M=\overbrace{SM\otimes\ldots\otimes 
SM}^{\text{$\alpha$ times}}\otimes
\underbrace{S^*\!M\otimes\ldots\otimes 
S^*\!M}_{\text{$\beta$ times}}\,,
\mytag{2.4}\\
\hskip -2em
\bar S^\nu_\gamma M=\overbrace{S^{{\sssize\dagger}*}\!M\otimes
\ldots\otimes S^{{\sssize\dagger}*}\!M}^{\text{$\nu$ times}}
\otimes\underbrace{S^{\sssize\dagger}\!M)\otimes\ldots\otimes 
S^{\sssize\dagger}\!M}_{\text{$\gamma$ times}}.
\mytag{2.5}
\endgather
$$     
Note that \mythetag{2.3} is the {\it complexified tensor bundle 
of the type $(m,n)$}, \mythetag{2.4} is the {\it spin-tensorial 
bundle of the type $(\alpha,\beta)$}, and \mythetag{2.5} is
the {\it barred spin-tensorial bundle of the type $(\nu,\gamma)$}.
Combining these three bundles, we define the following spin-tensorial 
bundle of the mixed type $(\alpha,\beta|\nu,\gamma|m,n)$:
$$
\hskip -2em
S^\alpha_\beta\bar S^\nu_\gamma T^m_n\!M=S^\alpha_\beta M
\otimes\bar S^\nu_\gamma M\otimes\Bbb CT^m_n\!M.
\mytag{2.6}
$$
Smooth sections of the bundle \mythetag{2.6} are called 
{\it spin-tensorial} fields of the type $(\alpha,\beta|\nu,\gamma|m,n)$.
The metric tensor $\bold g$ of the base space-time manifold $M$ now
is interpreted as a spin-tensorial field of the type $(0,0|0,0|0,2)$,
while the spin-metric tensor $\bold d$ is a spin-tensorial field of the
type $(0,2|0,0|0,0)$.\par
     Note that an arbitrary spin-tensorial field of the type $(0,0|0,0|0,2)$
is a complex field, while $\bold g$ is a real field. Therefore, we have
$$
\hskip -2em
\bold g=\tau(\bold g),
\mytag{2.7}
$$
where $\tau$ is the semilinear involution of complex conjugation:
$$
\hskip -2em
\CD
@>\tau>>\\
\vspace{-4ex}
S^\alpha_\beta\bar S^\nu_\gamma T^m_n M@.
S^\nu_\gamma\bar S^\alpha_\beta T^m_n M.\\
\vspace{-4.2ex}
@<<\tau< 
\endCD
\mytag{2.8}
$$
Both mappings \mythetag{2.8} are denoted by the same symbol, hence, 
formally we have the involution identity $\tau^2=\tau\compos\tau=\idop$.
More detailed description of the involution $\tau$ can be found in
\mycite{3} and \mycite{6}.\par
     A coordinate description of spin-tensorial fields is obtained in
terms of frame pairs. Let $(U,\,\boldsymbol\Upsilon_0,\,\boldsymbol
\Upsilon_1,\,\boldsymbol\Upsilon_2,\,\boldsymbol\Upsilon_3)$ be an
arbitrary frame of the tangent bundle $TM$ and let $(U,\,\boldsymbol
\Psi_1,\,\boldsymbol\Psi_2)$ be an arbitrary frame of the spinor bundle
$SM$. Denote by $(U,\,\boldsymbol\eta^0,\,\boldsymbol\eta^1,\,\boldsymbol
\eta^2,\,\boldsymbol\eta^3)$ and $(U,\,\boldsymbol\vartheta^{\,1},\,
\boldsymbol\vartheta^{\,2})$ the dual frames for $(U,\,\boldsymbol
\Upsilon_0,\,\boldsymbol\Upsilon_1,\,\boldsymbol\Upsilon_2,\,\boldsymbol
\Upsilon_3)$ and $(U,\,\boldsymbol\Psi_1,\,\boldsymbol\Psi_2)$ 
respectively. Then let's denote
$$
\xalignat 2
&\hskip -2em
\bPsi_i=\tau(\boldsymbol\Psi_i),
&&\overline{\boldsymbol\vartheta}\vphantom{\boldsymbol\vartheta}^{\,i}
=\tau(\boldsymbol\vartheta^{\,i}).
\mytag{2.9}
\endxalignat
$$
The barred spinor fields in \mythetag{2.9} compose two frames $(U,\,
\bPsi_1,\,\bPsi_2)$ and $(U,\,\overline{\boldsymbol
\vartheta}\vphantom{\boldsymbol\vartheta}^{\,1},\,\overline{\boldsymbol
\vartheta}\vphantom{\boldsymbol\vartheta}^{\,2})$ in $S^{\sssize
\dagger}\!M$ and $S^{*\sssize\dagger}\!M$ respectively. Now, according 
to the formulas \mythetag{2.3}, \mythetag{2.4}, and \mythetag{2.5}, we
define the following tensor products:
$$
\align
&\hskip -2em
\boldsymbol\Upsilon^{k_1\ldots\,k_n}_{h_1\ldots\,h_m}
=\boldsymbol\Upsilon_{h_1}\otimes\ldots\otimes\boldsymbol\Upsilon_{h_m}
\otimes\boldsymbol\eta^{k_1}\otimes\ldots\otimes\boldsymbol\eta^{k_n},
\mytag{2.10}\\
\vspace{2ex}
&\hskip -2em
\boldsymbol\Psi^{j_1\ldots\,j_\beta}_{i_1\ldots\,i_\alpha}
=\boldsymbol\Psi_{i_1}\otimes\ldots\otimes\boldsymbol\Psi_{i_\alpha}
\otimes\boldsymbol\vartheta^{\,j_1}\otimes\ldots\otimes
\boldsymbol\vartheta^{\,j_\beta},
\mytag{2.11}\\
\vspace{2ex}
&\hskip -2em
\bPsi^{\bar j_1\ldots\,\bar j_\gamma}_{\bar i_1\ldots\,\bar i_\nu}
=\bPsi_{\bar i_1}\otimes\ldots\otimes\bPsi_{\bar i_\nu}\otimes
\overline{\boldsymbol\vartheta}^{\,j_1}\otimes\ldots\otimes
\overline{\boldsymbol\vartheta}^{\,j_\gamma}.
\mytag{2.12}
\endalign
$$
And finally, according to \mythetag{2.6}, from \mythetag{2.10},
\mythetag{2.11}, and \mythetag{2.12} we produce
$$
\hskip -2em
\boldsymbol\Psi^{j_1\ldots\,j_\beta\,\bar j_1\ldots\,\bar j_\gamma
\,k_1\ldots\, k_n}_{i_1\ldots\,i_\alpha\,\bar i_1\ldots\,\bar i_\nu
\,h_1\ldots\,h_m}
=\boldsymbol\Psi^{j_1\ldots\,j_\beta}_{i_1\ldots\,i_\alpha}
\otimes\bPsi^{\bar j_1\ldots\,\bar j_\gamma}_{\bar i_1\ldots\,
\bar i_\nu}\otimes\boldsymbol\Upsilon^{k_1\ldots\,k_n}_{h_1\ldots\,h_m}.
\mytag{2.13}
$$
Using \mythetag{2.13}, for any spin-tensorial field of the type
$(\alpha,\beta|\nu,\gamma|m,n)$ we write
$$
\bold X=\msum{2}\Sb i_1,\,\ldots,\,i_\alpha\\ j_1,\,\ldots,\,j_\beta
\endSb\msum{2}\Sb\bar i_1,\,\ldots,\,\bar i_\nu\\ \bar j_1,\,\ldots,\,
\bar j_\gamma\endSb\msum{3}\Sb h_1,\,\ldots,\,h_m\\ k_1,\,\ldots,\,k_n
\endSb
X^{i_1\ldots\,i_\alpha\,\bar i_1\ldots\,\bar i_\nu\,h_1\ldots\,h_m}_{j_1
\ldots\,j_\beta\,\bar j_1\ldots\,\bar j_\gamma\,k_1\ldots\, k_n}\ 
\boldsymbol\Psi^{j_1\ldots\,j_\beta\,\bar j_1\ldots\,\bar j_\gamma\,
k_1\ldots\, k_n}_{i_1\ldots\,i_\alpha\,\bar i_1\ldots\,\bar i_\nu\,
h_1\ldots\,h_m}.\quad
\mytag{2.14}
$$
Since all of the above frames $(U,\,\boldsymbol\eta^0,\,\boldsymbol
\eta^1,\,\boldsymbol\eta^2,\,\boldsymbol\eta^3)$, $(U,\,\boldsymbol
\vartheta^{\,1},\,\boldsymbol\vartheta^{\,2})$, $(U,\,\bPsi_1,\,
\bPsi_2)$, and $(U,\,\overline{\boldsymbol\vartheta}\vphantom{\boldsymbol
\vartheta}^{\,1},\,\overline{\boldsymbol\vartheta}\vphantom{\boldsymbol
\vartheta}^{\,2})$ and their tensor products \mythetag{2.10}, 
\mythetag{2.11}, \mythetag{2.12}, and \mythetag{2.13} are produced from 
two initial frames $(U,\,\boldsymbol\Upsilon_0,\,\boldsymbol\Upsilon_1,\,
\boldsymbol\Upsilon_2,\,\boldsymbol\Upsilon_3)$ and $(U,\,\boldsymbol
\Psi_1,\,\boldsymbol\Psi_2)$, the coefficients  $X^{i_1\ldots\,i_\alpha\,
\bar i_1\ldots\,\bar i_\nu\,h_1\ldots\,h_m}_{j_1\ldots\,j_\beta\,\bar j_1
\ldots\,\bar j_\gamma\,k_1\ldots\, k_n}$ in \mythetag{2.14} are called
the {\it coordinate representation\/} of the field $\bold X$ in the frame
pair $(U,\,\boldsymbol\Upsilon_0,\,\boldsymbol\Upsilon_1,\,\boldsymbol
\Upsilon_2,\,\boldsymbol\Upsilon_3)$ and $(U,\,\boldsymbol\Psi_1,\,
\boldsymbol\Psi_2)$. When passing from this frame pair to another frame
pair $(\tilde U,\,\tilde{\boldsymbol\Upsilon}_0,\,\tilde{\boldsymbol
\Upsilon}_1,\,\tilde{\boldsymbol\Upsilon}_2,\,\tilde{\boldsymbol
\Upsilon}_3)$ and $(\tilde U,\,\tilde{\boldsymbol\Psi}_1,\,
\tilde{\boldsymbol\Psi}_2)$ these coefficients are transformed as follows:
$$
\align
&\hskip -4em
\aligned
&\tilde X^{i_1\ldots\,i_\alpha\,\bar i_1\ldots\,\bar i_\nu\,
h_1\ldots\,h_m}_{j_1\ldots\,j_\beta\,\bar j_1\ldots\,\bar j_\gamma\,
k_1\ldots\,k_n}
=\dsize\msum{4}\Sb a_1,\,\ldots,\,a_\alpha\\ b_1,\,\ldots,\,b_\beta\endSb
\dsize\msum{4}\Sb \bar a_1,\,\ldots,\,\bar a_\nu\\ 
\bar b_1,\,\ldots,\,\bar b_\gamma\endSb
\dsize\msum{3}
\Sb c_1,\,\ldots,\,c_m\\ d_1,\,\ldots,\,d_n\endSb
\hat{\goth T}^{\,i_1}_{a_1}\ldots\,\hat{\goth T}^{\,i_\alpha}_{a_\alpha}
\,\times\\
&\kern 40pt
\times\,\hat{\goth S}^{b_1}_{j_1}\ldots\,\hat{\goth S}^{b_\beta}_{j_\beta}
\ \overline{\hat{\goth T}^{\,\bar i_1}_{\bar a_1}}
\ldots\,\overline{\hat{\goth T}^{\,\bar i_\nu}_{\bar a_\nu}}\ 
\ \overline{\hat{\goth S}^{\,\bar b_1}_{\bar j_1}}
\ldots\,\overline{\hat{\goth S}^{\,\bar b_\gamma}_{\bar j_\gamma}}\ 
T^{h_1}_{c_1}\ldots\,T^{h_m}_{c_m}\,\times\\
\vspace{1.5ex}
&\kern 60pt
\times\,S^{\,d_1}_{k_1}\ldots\,S^{\,d_n}_{k_n}\ 
X^{\,a_1\ldots\,a_\alpha\,\bar a_1\ldots\,\bar a_\nu\,
c_1\ldots\,c_m}_{\,b_1\ldots\,b_\beta\,\bar b_1\ldots\,\bar b_\gamma
\,d_1\ldots\,d_n},
\endaligned
\mytag{2.15}\\
\vspace{2ex}
&\hskip -4em
\aligned
&X^{i_1\ldots\,i_\alpha\,\bar i_1\ldots\,\bar i_\nu\,
h_1\ldots\,h_m}_{j_1\ldots\,j_\beta\,\bar j_1\ldots\,\bar j_\gamma\,
k_1\ldots\,k_n}
=\dsize\msum{4}\Sb a_1,\,\ldots,\,a_\alpha\\ b_1,\,\ldots,\,b_\beta\endSb
\dsize\msum{4}\Sb \bar a_1,\,\ldots,\,\bar a_\nu\\ 
\bar b_1,\,\ldots,\,\bar b_\gamma\endSb
\dsize\msum{3}
\Sb c_1,\,\ldots,\,c_m\\ d_1,\,\ldots,\,d_n\endSb
\hat{\goth S}^{\,i_1}_{a_1}\ldots\,\hat{\goth S}^{\,i_\alpha}_{a_\alpha}
\,\times\\
&\kern 40pt
\times\,\hat{\goth T}^{b_1}_{j_1}\ldots\,\hat{\goth T}^{b_\beta}_{j_\beta}
\ \overline{\hat{\goth S}^{\,\bar i_1}_{\bar a_1}}
\ldots\,\overline{\hat{\goth S}^{\,\bar i_\nu}_{\bar a_\nu}}\ 
\ \overline{\hat{\goth T}^{\,\bar b_1}_{\bar j_1}}
\ldots\,\overline{\hat{\goth T}^{\,\bar b_\gamma}_{\bar j_\gamma}}\ 
S^{h_1}_{c_1}\ldots\,S^{h_m}_{c_m}\,\times\\
\vspace{1.5ex}
&\kern 60pt
\times\,T^{\,d_1}_{k_1}\ldots\,T^{\,d_n}_{k_n}\ 
\tilde X^{\,a_1\ldots\,a_\alpha\,\bar a_1\ldots\,\bar a_\nu\,
c_1\ldots\,c_m}_{\,b_1\ldots\,b_\beta\,\bar b_1\ldots\,\bar b_\gamma
\,d_1\ldots\,d_n}.
\endaligned
\mytag{2.16}
\endalign
$$
The formulas \mythetag{2.15} and \mythetag{2.16} express the general
transformation rules for the components of a chiral spin-tensorial field
under a change of frame pairs. The matrices $\goth S$, $\goth T$, $S$, 
and $T$ in them are taken from \mythetag{1.4} and \mythetag{1.6}.\par
    Note that $(U,\,\boldsymbol\Upsilon_0,\,\boldsymbol\Upsilon_1,\,
\boldsymbol\Upsilon_2,\,\boldsymbol\Upsilon_3)$ and $(U,\,\boldsymbol
\eta^0,\,\boldsymbol\eta^1,\,\boldsymbol\eta^2,\,\boldsymbol\eta^3)$
are real frames. They are invariant under the action of the semilinear
involution $\tau$:
$$
\xalignat 2
&\hskip -2em
\tau(\boldsymbol\Upsilon_i)=\boldsymbol\Upsilon_i,
&&\tau(\boldsymbol\eta_i)=\boldsymbol\eta_i.
\mytag{2.17}
\endxalignat
$$
Applying $\tau$ to \mythetag{2.14} and taking into account \mythetag{2.9} 
and \mythetag{2.17}, we obtain
$$
\gather
\hskip -1em
\tau(\bold X)=\msum{2}\Sb i_1,\,\ldots,\,i_\nu\\ j_1,\,\ldots,\,j_\gamma\\
\bar i_1,\,\ldots,\,\bar i_\alpha\\ \bar j_1,\,\ldots,\,
\bar j_\beta\endSb\msum{3}\Sb h_1,\,\ldots,\,h_m\\ k_1,\,\ldots,\,k_n
\endSb
\tau X^{i_1\ldots\,i_\nu\,\bar i_1\ldots\,\bar i_\alpha\,h_1\ldots\,
h_m}_{j_1\ldots\,j_\gamma\,\bar j_1\ldots\,\bar j_\beta\,k_1\ldots\, 
k_n}\ \boldsymbol\Psi^{j_1\ldots\,j_\gamma\,\bar j_1\ldots\,\bar j_\beta
\,k_1\ldots\, k_n}_{i_1\ldots\,i_\nu\,\bar i_1\ldots\,\bar i_\alpha\,
h_1\ldots\,h_m},\qquad
\mytag{2.18}\\
\vspace{-1ex}
\intertext{where}
\vspace{-2ex}
\tau X^{i_1\ldots\,i_\nu\,\bar i_1\ldots\,\bar i_\alpha\,h_1\ldots\,
h_m}_{j_1\ldots\,j_\gamma\,\bar j_1\ldots\,\bar j_\beta\,k_1\ldots\, 
k_n}
=\overline{X^{\bar i_1\ldots\,\bar i_\alpha\,i_1\ldots\,i_\nu\,h_1
\ldots\,h_m}_{\bar j_1\ldots\,\bar j_\beta\,j_1\ldots\,j_\gamma\,k_1
\ldots\,k_n}}.\qquad
\mytag{2.19}
\endgather
$$
Note that the formulas \mythetag{2.18} and \mythetag{2.19} are in
agreement with \mythetag{2.8}. They mean that the involution $\tau$
acts upon the components of spin-tensorial fields as the complex
conjugation exchanging barred and non-barred spinor indices.

\head
3. Basic spin-tensorial fields of chiral spinors.
\endhead
    The metric tensor $\bold g$ is the basic field for both chiral and
Dirac spinors. As it was already mentioned above, $\bold g$ is interpreted 
as a spin-tensorial field of the type $(0,0|0,0|0,2)$ satisfying the reality
condition \mythetag{2.7}. Apart from $\bold g$, in the theory of chiral
spinors there are two other basic spin-tensorial fields. The first of
them is the spin-metric tensor $\bold d$. It is a field of the type 
$(0,2|0,0|0,0)$. In canonically associated frame pairs (see 
definition~\mythedefinition{1.1}) its components are given by the matrix
\mythetag{1.5}. The second basic spin tensorial field in the theory of
chiral spinors is denoted by $\bold G$. It is called the {\it Infeld-van
der Waerden field}, its components $G^{i\kern 0.5pt\bar i}_q$ are called
the {\it Infeld-van der Waerden symbols}. In canonically associated frame
pairs the Infeld-van der Waerden symbols are given explicitly by the
formulas 
$$
\xalignat 4
&\hskip -2em
G^{11}_0=1, &&G^{11}_1=0, &&G^{11}_2=0,  &&G^{11}_3=1,\\
&\hskip -2em
G^{12}_0=0, &&G^{12}_1=1, &&G^{12}_2=-i, &&G^{12}_3=0,\\
\vspace{-1.4ex}
\mytag{3.1}\\
\vspace{-1.4ex}
&\hskip -2em
G^{21}_0=0, &&G^{21}_1=1, &&G^{21}_2=i,  &&G^{21}_3=0,\\
&\hskip -2em
G^{22}_0=1, &&G^{22}_1=0, &&G^{22}_2=0,  &&G^{22}_3=-1.\\
\endxalignat
$$
These formulas \mythetag{3.1} are derived from \mythetag{1.3} due 
to the formula \mythetag{1.2}. The Infeld-van der Waerden field is 
a spin-tensorial field of the type $(1,0|1,0|0,1)$.\par
     Applying the index lowering and index raising procedures to 
the Infeld-van der Waerden symbols \mythetag{3.1} we get the 
inverse\footnotemark\ Infeld-van der Waerden symbols:
$$
\hskip -2em
G^{\,q}_{i\kern 0.5pt\bar i}=\sum^2_{j=1}\sum^2_{\bar j=1}
\sum^3_{k=0}G^{j\kern 0.5pt\bar j}_k\ g^{kq}\,d_{ji}
\ \bd_{\bar j\bar i}.
\mytag{3.2}
$$
The inverse Infeld-van der Waerden given by the formula \mythetag{3.2}
are the components of a spin-tensorial field of the type $(0,1|0,1|1,0)$,
it is denoted by the same symbol $\bold G$ as the initial Infeld-van der
Waerden field. Here are the numeric values of the inverse Infeld-van der 
Waerden symbols:
$$
\xalignat 4
&\hskip -2em
G^{\,0}_{11}=1, &&G^{\,0}_{12}=0, &&G^{\,0}_{21}=0,  
&&G^{\,0}_{22}=1,\\
&\hskip -2em
G^{\,1}_{11}=0, &&G^{\,1}_{12}=1, &&G^{\,1}_{21}=1,  
&&G^{\,1}_{22}=0,\\
\vspace{-1.4ex}
\mytag{3.3}\\
\vspace{-1.4ex}
&\hskip -2em
G^{\,2}_{11}=0, &&G^{\,2}_{12}=i, &&G^{\,2}_{21}=-i,  
&&G^{\,2}_{22}=0,\\
&\hskip -2em
G^{\,3}_{11}=1, &&G^{\,3}_{12}=0, &&G^{\,3}_{21}=0,  
&&G^{\,3}_{22}=-1.
\endxalignat
$$
\footnotetext{\  Note that the definition \mythetag{3.2} of the inverse 
Infeld-van der Waerden symbols is different from that of \mycite{3}. Their
numeric values \mythetag{3.3} are also different from \thetag{7.9} in 
\mycite{3}. However, they differ only by the numeric factor $1/2$.}
The barred spin-metric tensor $\bbd$ is derived from $\bold d$ by applying
$\tau$:
$$
\hskip -2em
\bbd=\tau(\bold d).
\mytag{3.4}
$$
It is a spin-tensorial field of the type $(0,0|0,2|0,0)$. The components
of the field \mythetag{3.4} in \mythetag{3.2} are derived by means of the 
formula
$$
\hskip -2em
\bd_{\bar j\bar i}=\overline{d_{\bar j\bar i}\mathstrut}.
\mytag{3.5}
$$
This formula \mythetag{3.5} is a special case of the general formula
\mythetag{2.19}. \pagebreak In the case of the Infeld-van der Waerden 
field the reality condition for it is given by the formula:
$$
\hskip -2em
\tau(\bold G)=\bold G.
\mytag{3.6}
$$
Applying \mythetag{2.19} to \mythetag{3.6}, we get the following equalities:
$$
\xalignat 2
&\hskip -2em
G^{i\kern 0.5pt\bar i}_q=\overline{G^{\bar i\kern 0.5pt i}_q},
&&G^{\,q}_{i\kern 0.5pt\bar i}=\overline{G^{\,q}_{\bar i\kern 0.5pt i}}.
\mytag{3.7}
\endxalignat
$$
These equalities \mythetag{3.7} express the reality condition
\mythetag{3.6} in a coordinate form.\par
     The Infeld-van der Waerden symbols satisfy various identities relating
these symbols with the metric and spin-metric tensors:
$$
\gather
\hskip -2em
\sum^3_{p=0}\sum^3_{q=0}g_{p\,q}\ G^{\,p}_{i\kern 0.5pt\bar i}
\ G^{\,q}_{j\kern 0.5pt\bar j}=2\,d_{ij}\,\bd_{\,\bar i\kern 0.5pt\bar j},
\mytag{3.8}\\
\hskip -2em
\sum^2_{i=1}\sum^2_{j=1}\sum^2_{\bar i=1}\sum^2_{\bar j=1}
d_{ij}\ \bd_{\,\bar i\kern 0.5pt\bar j}
\ G^{i\kern 0.5pt\bar i}_p\ G^{j\kern 0.5pt\bar j}_q=2\,g_{p\,q},
\mytag{3.9}\\
\hskip -2em
\sum^3_{p=0}\sum^3_{q=0}g^{p\,q}\ G^{i\kern 0.5pt\bar i}_p
\ G^{j\kern 0.5pt\bar j}_q
=2\,d^{\,ij}\,\bd^{\,\bar i\kern 0.5pt \bar j},
\mytag{3.10}\\
\hskip -2em
\sum^2_{i=1}\sum^2_{j=1}\sum^2_{\bar i=1}\sum^2_{\bar j=1}
d^{\,ij}\,\bd^{\,\bar i\kern 0.5pt \bar j}
\,G^{\,p}_{i\kern 0.5pt\bar i}
\ G^{\,q}_{j\kern 0.5pt\bar j}=2\,g^{p\,q}.
\mytag{3.11}
\endgather
$$
The following equalities are easily derived from \mythetag{3.8}, 
\mythetag{3.9}, \mythetag{3.10}, and \mythetag{3.11}:
$$
\xalignat 2
&\hskip -2em
\sum^2_{i=1}\sum^2_{\bar i=1}G^{i\kern 0.5pt\bar i}_p\ 
G^{\,q}_{i\kern 0.5pt\bar i}=2\,\delta^q_p,
&&\sum^3_{q=0}G^{i\kern 0.5pt\bar i}_q\ 
G^{\,q}_{j\kern 0.5pt\bar j}=2\,\delta^i_j
\,\delta^{\bar i}_{\bar j}.\qquad
\mytag{3.12}
\endxalignat
$$
The equalities \mythetag{3.8}, \mythetag{3.9}, \mythetag{3.10}, 
\mythetag{3.11}, and \mythetag{3.12} here are slightly different 
from \thetag{7.14}, \thetag{7.15}, \thetag{7.16}, \thetag{7.17},
and \thetag{7.13} in \mycite{3} since the quantities $G^{\,q}_{i
\kern 0.5pt\bar i}$ here differ from those of \mycite{3} by the
numeric factor $1/2$. The components of the dual spin-metric tensors
in \mythetag{3.10} and \mythetag{3.11} are given by the matrices
inverse to $d_{ij}$ and $\bd_{\bar i\bar j}$:
$$
\xalignat 2
&\sum^2_{q=1}d_{i\kern 0.5pt q}\,d^{\,qj}=\delta^{\,j}_i,
&&\sum^2_{\bar q=1}\bd_{\bar i\kern 0.5pt\bar q}\ 
\bd^{\,\bar q\bar j}=\delta^{\,\bar j}_{\,\bar i}.
\endxalignat
$$
The identities \mythetag{3.8}, \mythetag{3.9}, \mythetag{3.10}, 
\mythetag{3.11} are easily derived in canonically associated
frame pairs (see definition~\mythedefinition{1.1}). However, due 
to the spin-tensorial nature of the quantities in them they remain 
valid for arbitrary two frames of $TM$ and $SM$.\par
\head
4. Metric connections for chiral spinors.
\endhead
     Let $(U,\,\boldsymbol\Upsilon_0,\,\boldsymbol\Upsilon_1,\,
\boldsymbol\Upsilon_2,\,\boldsymbol\Upsilon_3)$ and $(U,\,\boldsymbol
\Psi_1,\,\boldsymbol\Psi_2)$ be two frames with a common domain $U$ of 
the bundles $TM$ and $SM$ respectively. Let $(\tilde U,\,
\tilde{\boldsymbol\Upsilon}_0,\,\tilde{\boldsymbol\Upsilon}_1,\,
\tilde{\boldsymbol\Upsilon}_2,\,\tilde{\boldsymbol\Upsilon}_3)$ 
and $(\tilde U,\,\tilde{\boldsymbol\Psi}_1,\,\tilde{\boldsymbol\Psi}_2)$
be other two such frames. Assume that $U\cap\tilde U\neq\varnothing$.
Then at each point $p$ of the intersection $U\cap\tilde U$ one can write
the transition formulas \mythetag{1.4} and \mythetag{1.6}. Assume that
the domains $U$ and $\tilde U$ are small enough so that one can introduce
local coordinates $x^0,\,x^1,\,x^2,\,x^3$ and $\tx^0,\,\tx^1,\,\tx^2,\,
\tx^3$ in them. Then, apart from the frames $(U,\,\boldsymbol\Upsilon_0,\,
\boldsymbol\Upsilon_1,\,\boldsymbol\Upsilon_2,\,\boldsymbol\Upsilon_3)$
and $(\tilde U,\,\tilde{\boldsymbol\Upsilon}_0,\,\tilde{\boldsymbol
\Upsilon}_1,\,\tilde{\boldsymbol\Upsilon}_2,\,\tilde{\boldsymbol
\Upsilon}_3)$, which are non-holonomic in general case, we have two
holonomic coordinate frames $(U,\,\bold E_0,\,\bold E_1,\,\bold E_2,\,
\bold E_3)$ and $(U,\,\tilde{\bold E}_0,\,\tilde{\bold E}_1,\,
\tilde{\bold E}_2,\,\tilde{\bold E}_3)$ composed by the vector fields
$$
\xalignat 2
&\hskip -2em
\bold E_i=\frac{\partial}{\partial x^i},
&&\tilde{\bold E}_i=\frac{\partial}{\partial\tx^i}.
\mytag{4.1}
\endxalignat
$$  
Taking the expansions of $\boldsymbol\Upsilon_i$ and $\tilde{\boldsymbol
\Upsilon}_i$ in these holonomic frames
$$
\xalignat 2
&\hskip -2em
\boldsymbol\Upsilon_i=\sum^3_{j=0}\Upsilon^j_i\ \bold E_j,
&&\tilde{\boldsymbol\Upsilon}_i=\sum^3_{j=0}\tilde\Upsilon^j_i
\ \tilde{\bold E}_j,
\mytag{4.2}
\endxalignat
$$  
due to \mythetag{4.1} and \mythetag{4.2} we can represent
$\boldsymbol\Upsilon_i$ and $\tilde{\boldsymbol\Upsilon}_i$ as linear 
differential operators
$$
\align
\hskip -2em
\boldsymbol\Upsilon_i&=\Upsilon^0_i\,\frac{\partial}{\partial x^0}
+\Upsilon^1_i\,\frac{\partial}{\partial x^1}+\Upsilon^2_i\,
\frac{\partial}{\partial x^2}+\Upsilon^3_i\,\frac{\partial}
{\partial x^3},
\mytag{4.3}\\
\vspace{2ex}
\hskip -2em
\tilde{\boldsymbol\Upsilon}_i&=\tilde\Upsilon^0_i\,\frac{\partial}
{\partial\tx^0}+\tilde\Upsilon^1_i\,\frac{\partial}{\partial\tx^1}
+\tilde\Upsilon^2_i\,\frac{\partial}{\partial\tx^2}
+\tilde\Upsilon^3_i\,\frac{\partial}{\partial\tx^3}.
\mytag{4.4}
\endalign
$$
Applying the differential operators \mythetag{4.3} and \mythetag{4.4}
to a smooth scalar function $f$, we denote the resulting functions
as the Lie derivatives:
$$
\xalignat 2
&\hskip -2em
L_{\boldsymbol\Upsilon_i}(f)=\sum^3_{j=0}\Upsilon^j_i\ 
\frac{\partial f}{\partial x^i},
&&L_{\tilde{\boldsymbol\Upsilon}_i}(f)=\sum^3_{j=0}\tilde
\Upsilon^j_i\  \frac{\partial f}{\partial\tx^i}.
\mytag{4.5}
\endxalignat
$$
Note that the components of the transition matrices $S$, $T$,
$\goth S$, $\goth T$ from \mythetag{1.4} and \mythetag{1.6} are 
smooth functions within the intersection domain $U\cap\tilde U$.
Therefore, one can substitute them for $f$ into \mythetag{4.5}.
As a result we can define the following functions:
$$
\align
&\hskip -2em
\tilde\theta^k_{ij}=\sum^3_{a=0}T^k_a\,L_{\tilde{\boldsymbol\Upsilon}_i}
\!(S^a_j)=-\sum^3_{a=0}L_{\tilde{\boldsymbol\Upsilon}_i}
\!(T^k_a)\,S^a_j,
\mytag{4.6}\\
&\hskip -2em
\tilde\vartheta^k_{ij}=\sum^2_{a=1}\goth T^k_a\,
L_{\tilde{\boldsymbol\Upsilon}_i}(\goth S^a_j)
=-\sum^2_{a=1}L_{\tilde{\boldsymbol\Upsilon}_i}
\!(\goth T^k_a)\,\goth S^a_j,
\mytag{4.7}\\
&\hskip -2em
\theta^k_{ij}=\sum^3_{a=0}S^k_a\,L_{\boldsymbol\Upsilon_i}
\!(T^a_j)=-\sum^3_{a=0}L_{\boldsymbol\Upsilon_i}
\!(S^k_a)\,T^a_j,
\mytag{4.8}\\
&\hskip -2em
\vartheta^k_{ij}=\sum^2_{a=1}\goth S^k_a\,
L_{\boldsymbol\Upsilon_i}(\goth T^a_j)
=-\sum^2_{a=1}L_{\boldsymbol\Upsilon_i}
\!(\goth S^k_a)\,\goth T^a_j.
\mytag{4.9}
\endalign
$$
The $\theta$-parameters with and without tilde introduced by the
above formulas \mythetag{4.6}, \mythetag{4.7}, \mythetag{4.8},
\mythetag{4.9} are related to each other through the following 
formulas:
$$
\allowdisplaybreaks
\align
&\hskip -2em
\theta^k_{ij}=-\sum^3_{a=0}\sum^3_{b=0}\sum^3_{c=0}
T^a_i\,\tilde\theta^{\,c}_{\!ab}\ S^k_c\,T^b_j,
\mytag{4.10}\\
&\hskip -2em
\tilde\theta^k_{ij}=-\sum^3_{a=0}\sum^3_{b=0}\sum^3_{c=0}
S^a_i\,\theta^{\,c}_{\!ab}\ T^k_c\,S^b_j,
\mytag{4.11}\\
&\hskip -2em
\vartheta^k_{ij}=-\sum^3_{a=0}\sum^2_{b=1}\sum^2_{c=1}
T^a_i\,\tilde\vartheta^{\,c}_{\!ab}\ \goth S^k_c\,\goth T^b_j,
\mytag{4.12}\\
&\hskip -2em
\tilde\vartheta^k_{ij}=-\sum^3_{a=0}\sum^2_{b=1}\sum^2_{c=1}
S^a_i\,\vartheta^{\,c}_{\!ab}\ \goth T^k_c\,\goth S^b_j.
\mytag{4.13}
\endalign
$$
In general case $\theta^k_{ij}$ and $\tilde\theta^k_{ij}$ are 
asymmetric in their lower indices and the extent of this asymmetry 
is characterized by the formulas
$$
\xalignat 2
&\hskip -2em
\theta^k_{ij}-\theta^k_{j\,i}=c^k_{ij},
&&\tilde\theta^k_{ij}-\tilde\theta^k_{j\,i}=\tilde c^k_{ij},
\mytag{4.14}
\endxalignat
$$
where $c^k_{ij}$ and $\tilde c^k_{ij}$ are defined by the following
commutator relationships:
$$
\xalignat 2
&\hskip -2em
[\boldsymbol\Upsilon_i,\,\boldsymbol\Upsilon_j]=
\sum^3_{k=0}c^{\,k}_{ij}\,\boldsymbol\Upsilon_k,
&&[\tilde{\boldsymbol\Upsilon}_i,\,\tilde{\boldsymbol\Upsilon}_j]=
\sum^3_{k=0}\tilde c^{\,k}_{ij}\,\tilde{\boldsymbol\Upsilon}_k.
\mytag{4.15}
\endxalignat
$$
The quantities $c^{\,k}_{ij}$ and $\tilde c^{\,k}_{ij}$ in \mythetag{4.14}
and \mythetag{4.15} are similar to structural constants of Lie algebras.
For this reason they are called the {\it structural constants} of the
frames $(U,\,\boldsymbol\Upsilon_0,\,\boldsymbol\Upsilon_1,\,\boldsymbol
\Upsilon_2,\,\boldsymbol\Upsilon_3)$ and $(\tilde U,\,\tilde{\boldsymbol
\Upsilon}_0,\,\tilde{\boldsymbol\Upsilon}_1,\,\tilde{\boldsymbol
\Upsilon}_2,\,\tilde{\boldsymbol\Upsilon}_3)$, though actually they are 
not constants, but smooth real-valued functions within the domains $U$ 
and $\tilde U$ respectively. As for the identities \mythetag{4.10},
\mythetag{4.11}, \mythetag{4.12}, and \mythetag{4.13}, they are easily
derived from \mythetag{4.6}, \mythetag{4.7}, \mythetag{4.8}, and 
\mythetag{4.9} due to \mythetag{4.5}.
\mydefinition{4.1} A {\it spinor connection} of the bundle of chiral 
spinors $SM$ is a geometric object such that in each frame pair $(U,\,
\boldsymbol\Upsilon_0,\,\boldsymbol\Upsilon_1,\,\boldsymbol\Upsilon_2,
\,\boldsymbol\Upsilon_3)$ and $(U,\,\boldsymbol\Psi_1,\,\boldsymbol
\Psi_2)$ of $TM$ and $SM$ it is given by three arrays of smooth
complex-valued functions
$$
\align
\Gamma^k_{ij}&=\Gamma^k_{ij}(p),\quad i,j,k=0,\,\ldots,\,3,\\
\Alpha^k_{ij}&=\Alpha^k_{ij}(p),\quad i=0,\,\ldots,\,3,\quad j,k=1,2,\\
\bar{\Alpha}\vphantom{\Alpha}^k_{ij}
&=\bar{\Alpha}\vphantom{\Alpha}^k_{ij}(p),\quad i=0,\,\ldots,\,3,\quad 
j,k=1,2,
\endalign
$$
where $p\in U$, such that when passing from $(U,\,\boldsymbol\Upsilon_0,
\,\boldsymbol\Upsilon_1,\,\boldsymbol\Upsilon_2,\,\boldsymbol\Upsilon_3)$
and $(U,\,\boldsymbol\Psi_1,\,\boldsymbol\Psi_2)$ to some other frame 
pair $(\tilde U,\,\tilde{\boldsymbol\Upsilon}_0,\,\tilde{\boldsymbol
\Upsilon}_1,\,\tilde{\boldsymbol\Upsilon}_2,\,\tilde{\boldsymbol
\Upsilon}_3)$ and $(\tilde U,\,\tilde{\boldsymbol\Psi}_1,\,
\tilde{\boldsymbol\Psi}_2)$ with $U\cap\tilde U\neq\varnothing$ 
these functions are transformed as follows:
$$
\align
&\hskip -2em
\Gamma^k_{ij}=\dsize\sum^3_{b=0}\sum^3_{a=0}\sum^3_{c=0}
S^k_a\,T^b_j\,T^c_i\ \tilde\Gamma^a_{c\,b}+\theta^k_{ij},
\mytag{4.16}\\
&\hskip -2em
\Alpha^k_{ij}=\dsize\sum^2_{b=1}\sum^2_{a=1}\sum^3_{c=0}
\goth S^k_a\,\goth T^b_j\,T^c_i\ \tilde{\Alpha}\vphantom{\Alpha}^a_{c\,b}
+\vartheta^k_{ij},
\mytag{4.17}\\
&\hskip -2em
\bar{\Alpha}\vphantom{\Alpha}^k_{ij}=\sum^2_{b=1}\sum^2_{a=1}
\sum^3_{c=0}\overline{\goth S^k_a}\ \overline{\goth T^b_j}\,T^c_i\  
\tilde{\bar{\Alpha}}\vphantom{\Alpha}^a_{c\,b}
+\overline{\vartheta^k_{ij}}.
\mytag{4.18}
\endalign
$$
\enddefinition
The components of transition matrices $S$, $T$, $\goth S$, and $\goth T$
in \mythetag{4.16}, \mythetag{4.17}, and \mythetag{4.18} are taken from
\mythetag{1.4} and \mythetag{1.6}, while the quantities $\theta^k_{ij}$
and $\vartheta^k_{ij}$ are defined in \mythetag{4.8} and \mythetag{4.9}.
Spinor connections introduced by the definition~\mythedefinition{4.1} 
are used in order to define covariant differentiations acting upon
spin-tensorial fields and producing other spin-tensorial fields from 
them. The covariant differential $\nabla$ associated with the spinor
connection $(\Gamma,\Alpha,\bar{\Alpha})$ is a differential operator
$$
\hskip -2em
\nabla\!:\,S^\alpha_\beta\bar S^\nu_\gamma T^m_n M
\to S^\alpha_\beta\bar S^\nu_\gamma T^m_{n+1} M.
\mytag{4.19}
$$
In a frame pair $(U,\,\boldsymbol\Upsilon_0,\,\boldsymbol\Upsilon_1,
\,\boldsymbol\Upsilon_2,\,\boldsymbol\Upsilon_3)$ and $(U,\,\boldsymbol
\Psi_1,\,\boldsymbol\Psi_2)$, i\.\,e\. in a coordinate form, the operator
\mythetag{4.19} is represented by the corresponding covariant derivative 
$$
\hskip -5em
\gathered
\nabla_{\!k_{n+1}}X^{i_1\ldots\,i_\alpha\,\bar i_1\ldots\,\bar i_\nu
\,h_1\ldots\,h_m}_{j_1\ldots\,j_\beta\,\bar j_1\ldots\,\bar j_\gamma
\,k_1\ldots\, k_n}
=L_{\boldsymbol\Upsilon_{k_{n+1}}}\!\bigl(X^{i_1\ldots\,i_\alpha\,\bar i_1
\ldots\,\bar i_\nu\,h_1\ldots\,h_m}_{j_1\ldots\,j_\beta\,\bar j_1
\ldots\,\bar j_\gamma\,k_1\ldots\, k_n}\bigr)\,-\\
\vspace{2ex}
\gathered
\kern -9em
+\sum^\alpha_{\mu=1}\sum^2_{v_\mu=1}\Alpha^{i_\mu}_{k_{n+1}\,v_\mu}\ 
X^{i_1\ldots\,v_\mu\,\ldots\,i_\varepsilon\,\bar i_1\ldots\,\bar i_\nu
\,h_1\ldots\,h_m}_{j_1\ldots\,\ldots\,\ldots\,j_\beta\,\bar j_1\ldots\,
\bar j_\gamma\,k_1\ldots\,k_n}\,-\\
\kern 9em-\sum^\beta_{\mu=1}\sum^2_{w_\mu=1}\Alpha^{w_\mu}_{k_{n+1}\,j_\mu}
\ X^{i_1\ldots\,\ldots\,\ldots\,i_\alpha\,\bar i_1\ldots\,\bar i_\nu
h_1\ldots\,h_m}_{j_1\ldots\,w_\mu\,\ldots\,j_\beta\,\bar j_1\ldots\,
\bar j_\gamma k_1\ldots\,k_n}\,+\\
\kern -9em
+\sum^\nu_{\mu=1}\sum^2_{v_\mu=1}
\bar{\Alpha}\vphantom{\Alpha}^{\bar i_\mu}_{k_{n+1}\,v_\mu}\ 
X^{i_1\ldots\,i_\alpha\,\bar i_1\ldots\,v_\mu\,\ldots\,\bar i_\nu
\,h_1\ldots\,h_m}_{j_1\ldots\,j_\beta\,\bar j_1\ldots\,\ldots\,\ldots\,
\bar j_\gamma\,k_1\ldots\,k_n}\,-\\
\kern 9em-\sum^\gamma_{\mu=1}\sum^2_{w_\mu=1}
\bar{\Alpha}\vphantom{\Alpha}^{w_\mu}_{k_{n+1}\,\bar j_\mu}\
X^{i_1\ldots\,i_\alpha\,\bar i_1\ldots\,\ldots\,\ldots\,\bar i_\nu
\,h_1\ldots\,h_m}_{j_1\ldots\,j_\beta\,\bar j_1\ldots\,w_\mu\,\ldots\,
\bar j_\gamma\,k_1\ldots\,k_n}\,+\\
\kern -9em+\sum^m_{\mu=1}\sum^3_{v_\mu=0}\Gamma^{h_\mu}_{k_{n+1}\,v_\mu}\ 
X^{i_1\ldots\,i_\alpha\,\bar i_1\ldots\,\bar i_\nu\,
h_1\ldots\,v_\mu\,\ldots\,h_m}_{j_1\ldots\,j_\beta\,\bar j_1\ldots\,
\bar j_\gamma\,k_1\ldots\,\ldots\,\ldots\,k_n}\,-\\
\kern 9em-\sum^n_{\mu=1}\sum^3_{w_\mu=0}\Gamma^{w_\mu}_{k_{n+1}\,k_\mu}\
X^{i_1\ldots\,i_\alpha\,\bar i_1\ldots\,\bar i_\nu\,
h_1\ldots\,\ldots\,\ldots\,h_m}_{j_1\ldots\,j_\beta\,\bar j_1\ldots\,
\bar j_\gamma\,k_1\ldots\,w_\mu\,\ldots\,k_n}.
\endgathered\kern 4em
\endgathered\hskip -4em
\mytag{4.20}
$$
The formula \mythetag{4.20} should be understood in the following way.
If $\bold X$ is a spin-tensorial field of the type $(\alpha,\beta|\nu,
\gamma|m,n)$ and $X^{i_1\ldots\,i_\alpha\,\bar i_1\ldots\,\bar i_\nu
\,h_1\ldots\,h_m}_{j_1\ldots\,j_\beta\,\bar j_1\ldots\,\bar j_\gamma
\,k_1\ldots\, k_n}$ is its coordinate representation in the expansion
\mythetag{2.14}, then for the spin-tensorial field $\bold Y=
\nabla\bold X$ its coordinate representation is given by the formula
$$
Y^{i_1\ldots\,i_\alpha\,\bar i_1\ldots\,\bar i_\nu
\,h_1\ldots\,h_m}_{j_1\ldots\,j_\beta\,\bar j_1\ldots\,\bar j_\gamma
\,k_1\ldots\,k_{n+1}}=\nabla_{\!k_{n+1}}X^{i_1\ldots\,i_\alpha\,
\bar i_1\ldots\,\bar i_\nu\,h_1\ldots\,h_m}_{j_1\ldots\,j_\beta\,
\bar j_1\ldots\,\bar j_\gamma\,k_1\ldots\, k_n}.
$$
\mydefinition{4.2} A spinor connection $(\Gamma,\Alpha,\bar{\Alpha})$
of the bundle of chiral spinors $SM$ is called {\it concordant with 
the complex conjugation\/} if the corresponding covariant differential 
\mythetag{4.19} commute with the involution $\tau$, i\.
\,e\. if $\nabla(\tau(\bold X))=\tau(\nabla\bold X)$ for any 
spin-tensorial field $\bold X$.
\enddefinition
     Spinor connections concordant with the complex conjugation $\tau$
in the sense of the above definition~\mythedefinition{4.2} are also 
called {\it real connections}. 
\mytheorem{4.1} A spinor connection $(\Gamma,\Alpha,\bar{\Alpha})$
of the bundle of chiral spinors $SM$ is concordant with the complex 
conjugation $\tau$ if and only if 
$$
\xalignat 2
&\hskip -2em
\Gamma^k_{ij}=\overline{\Gamma^k_{ij}},
&&\bar{\Alpha}\vphantom{\Alpha}^k_{ij}=\overline{\Alpha^k_{ij}}.
\mytag{4.21}
\endxalignat
$$
\endproclaim
     The theorem~\mythetheorem{4.1} is proved by direct calculations 
on the base of the formula \mythetag{4.20}. The first relationship in
\mythetag{4.21} means that $\Gamma$-components of a real spinor
connection are real functions. They obey the transformation rules
\mythetag{4.16} coinciding with the transformation rules for the
components of an affine connection. 
\mycorollary{4.1} Any real spinor connection $(\Gamma,\Alpha,\bar{\Alpha})$
of the bundle of chiral spinors $SM$ comprises some affine connection
$\Gamma$ as its constituent part.
\endproclaim
\mydefinition{4.3} A spinor connection $(\Gamma,\Alpha,\bar{\Alpha})$
of the bundle of chiral spinors $SM$ is called {\it concordant with the
Infeld-van der Waerden field\/} if $\nabla\bold G=0$.
\enddefinition
\mydefinition{4.4} A spinor connection $(\Gamma,\Alpha,\bar{\Alpha})$
of the bundle of chiral spinors $SM$ is called {\it concordant with the
spin-metric tensor\/} if $\nabla\bold d=0$.
\enddefinition
\mydefinition{4.5} A spinor connection $(\Gamma,\Alpha,\bar{\Alpha})$
of the bundle of chiral spinors $SM$ is called {\it concordant with the
metric tensor\/} if $\nabla\bold g=0$.
\enddefinition
\mytheorem{4.2} Any real spinor connection $(\Gamma,\Alpha,\bar{\Alpha})$
of the bundle of chiral spinors $SM$ concordant with the Infeld-van der
Waerden field $\bold G$ and with the spin-metric tensor $\bold d$ is 
concordant with the metric tensor $\bold g$ too, i\.\,e\. for a real 
spinor connection $\nabla\bold G=0$ and $\nabla\bold d=0$ imply
$\nabla\bold g=0$.
\endproclaim
\demo{Proof} Since 
$(\Gamma,\Alpha,\bar{\Alpha})$ is real, from $\nabla\bold d=0$ we easily
derive that $\nabla\bbd=0$:
$$
\nabla\bbd=\nabla(\tau(\bold d))=\tau(\nabla\bold d)=0.
$$
Then we apply $\nabla_{\!r}$ to the identity \mythetag{3.9}. As a result 
we get
$$
\gathered
2\,\nabla_{\!r}g_{p\,q}
=\sum^2_{i=1}\sum^2_{j=1}\sum^2_{\bar i=1}\sum^2_{\bar j=1}
\left(\nabla_{\!r}d_{ij}\ \bd_{\,\bar i\kern 0.5pt\bar j}
\ G^{i\kern 0.5pt\bar i}_p\ G^{j\kern 0.5pt\bar j}_q
+\,d_{ij}\ \nabla_{\!r}\bd_{\,\bar i\kern 0.5pt\bar j}\,\times\right.\\
\left.\times\,G^{i\kern 0.5pt\bar i}_p\ G^{j\kern 0.5pt\bar j}_q
+d_{ij}\ \bd_{\,\bar i\kern 0.5pt\bar j}
\ \nabla_{\!r}G^{i\kern 0.5pt\bar i}_p\ G^{j\kern 0.5pt\bar j}_q
+d_{ij}\ \bd_{\,\bar i\kern 0.5pt\bar j}
\ G^{i\kern 0.5pt\bar i}_p\ \nabla_{\!r}G^{j\kern 0.5pt\bar j}_q
\right)=0.
\endgathered
\mytag{4.22}
$$
The identity \mythetag{4.22} means that $\nabla\bold g=0$. Thus, the 
theorem~\mythetheorem{4.2} is proved.\qed\enddemo
    The concordance condition $\nabla\bold g=0$ is well-known. Since
$\bold g$ is a spin-tensorial field of the type $(0,0|0,0|0,2)$, this
condition is written in terms of the $\Gamma$-components of a spinor 
connection only. Applying the formula \mythetag{4.20} to $\nabla_{\!r}
g_{ij}$, we get 
$$
\hskip -2em
L_{\boldsymbol\Upsilon_{\!r}}\!(g_{ij})-\sum^3_{s=0}\Gamma^s_{\!ri}
\ g_{sj}-\sum^3_{s=0}\Gamma^s_{\!rj}\ g_{is}=0.
\mytag{4.23}
$$
In general non-holonomic frame $(U,\,\boldsymbol\Upsilon_0,\,\boldsymbol
\Upsilon_1,\,\boldsymbol\Upsilon_2,\,\boldsymbol\Upsilon_3)$ the
$\Gamma$-components of a spinor connection are not symmetric. Therefore
we subdivide $\Gamma^k_{\!rs}$ into symmetric a skew-symmetric parts 
$\hat\Gamma^k_{\!rs}$ and $\check\Gamma^k_{\!rs}$ respectively:
$$
\hskip -2em
\Gamma^k_{\!rs}=\hat\Gamma^k_{\!rs}+\check\Gamma^k_{\!rs}.
\mytag{4.24}
$$
Lowering the upper index $k$ of $\Gamma^k_{\!rs}$, $\hat\Gamma^k_{\!rs}$, 
and $\check\Gamma^k_{\!rs}$, we define the following quantities:
$$
\xalignat 3
&\Gamma_{\!rsq}=\sum^3_{k=0}\Gamma^k_{\!rs}\ g_{kq},
&&\hat\Gamma_{\!rsq}=\sum^3_{k=0}\hat\Gamma^k_{\!rs}\ g_{kq},
&&\check\Gamma_{\!rsq}=\sum^3_{k=0}\check\Gamma^k_{\!rs}\ g_{kq}.
\qquad\quad
\mytag{4.25}
\endxalignat
$$
From \mythetag{4.24} we derive the analogous expansion of $\Gamma_{\!rsq}$
into two parts 
$$
\hskip -2em
\Gamma_{\!rsq}=\hat\Gamma_{\!rsq}+\check\Gamma_{\!rsq},
\mytag{4.26}
$$
where $\hat\Gamma_{\!rsq}=\hat\Gamma_{\!srq}$ and $\check\Gamma_{\!rsq}
=-\check\Gamma_{\!srq}$. Applying \mythetag{4.25} and \mythetag{4.26}
to \mythetag{4.23}, we get
$$
\hskip -2em
\hat\Gamma_{\!rij}+\hat\Gamma_{\!rj\kern 0.5pt i}
=L_{\boldsymbol\Upsilon_{\!r}}\!(g_{ij})
-\check\Gamma_{\!rij}-\check\Gamma_{\!rj\kern 0.5pt i}.
\mytag{4.27}
$$
By means of the cyclic transposition of indices from \mythetag{4.27} we
derive
$$
\align
&\hskip -2em
\hat\Gamma_{\!ijr}+\hat\Gamma_{\!i\kern 0.5pt rj}
=L_{\boldsymbol\Upsilon_{\!i}}\!(g_{jr})
-\check\Gamma_{\!ijr}-\check\Gamma_{\!i\kern 0.5pt rj},
\mytag{4.28}\\
\vspace{1ex}
&\hskip -2em
\hat\Gamma_{\!jri}+\hat\Gamma_{\!ji\kern 0.5pt r}
=L_{\boldsymbol\Upsilon_{\!j}}\!(g_{ri})
-\check\Gamma_{\!jri}-\check\Gamma_{\!ji\kern 0.5pt r}.
\mytag{4.29}
\endalign
$$
Now let's add \mythetag{4.28} and \mythetag{4.29}, then subtract
\mythetag{4.27} from the sum. As a result, taking into account the 
symmetry of $\hat\Gamma$ and the skew-symmetry of $\check\Gamma$, 
we get
$$
\hskip -2em
\hat\Gamma_{\!ijr}=\frac{L_{\boldsymbol\Upsilon_{\!i}}\!(g_{jr})
+L_{\boldsymbol\Upsilon_{\!j}}\!(g_{ri})-L_{\boldsymbol
\Upsilon_{\!r}}\!(g_{ij})}{2}-\check\Gamma_{\!i\kern 0.5pt rj}
-\check\Gamma_{\!jri}.
\mytag{4.30}
$$
Now, raising the lower index $r$ in \mythetag{4.30}, we obtain the
explicit formula for the symmetric part of the $\Gamma$-symbols
expressing them through the Lie derivatives $L_{\boldsymbol\Upsilon_{\!i}}
\!(g_{jr})$, $L_{\boldsymbol\Upsilon_{\!j}}\!(g_{ri})$, $L_{\boldsymbol
\Upsilon_{\!r}}\!(g_{ij})$ and through the skew-symmetric part of these 
$\Gamma$-symbols:
$$
\hskip -2em
\gathered
\hat\Gamma^k_{\!ij}=\sum^3_{r=0}\frac{g^{\kern 0.5pt kr}}{2}
\left(L_{\boldsymbol\Upsilon_{\!i}}\!(g_{jr})+L_{\boldsymbol
\Upsilon_{\!j}}\!(g_{ri})-L_{\boldsymbol\Upsilon_{\!r}}\!(g_{ij})
\right)\,-\\
-\sum^3_{r=0}\sum^3_{s=0}\check\Gamma^s_{\!i\kern 0.5pt r}\,g^{kr}\,
g_{sj}-\sum^3_{r=0}\sum^3_{s=0}\check\Gamma^s_{\!jr}\,g^{kr}\,g_{s
\kern 0.5pt i}.
\endgathered
\mytag{4.31}
$$
From \mythetag{4.31} and \mythetag{4.24} for the $\Gamma$-symbols 
themselves we derive
$$
\hskip -2em
\gathered
\Gamma^k_{\!ij}=\sum^3_{r=0}\frac{g^{\kern 0.5pt kr}}{2}
\left(L_{\boldsymbol\Upsilon_{\!i}}\!(g_{jr})+L_{\boldsymbol
\Upsilon_{\!j}}\!(g_{ri})-L_{\boldsymbol\Upsilon_{\!r}}\!(g_{ij})
\right)+\\
+\,\check\Gamma^k_{\!ij}
-\sum^3_{r=0}\sum^3_{s=0}\check\Gamma^s_{\!i\kern 0.5pt r}\,g^{kr}\,
g_{sj}-\sum^3_{r=0}\sum^3_{s=0}\check\Gamma^s_{\!jr}\,g^{kr}\,g_{s
\kern 0.5pt i}.
\endgathered
\mytag{4.32}
$$
Note that the skew-symmetric part of the $\Gamma$-symbols is determined
by the {\it torsion tensor} $\bold T$ (see \mycite{5}) and by the
structural constants $c^{\,k}_{ij}$ (see their definition \mythetag{4.15}):
$$
\hskip -2em
\check\Gamma^k_{\!ij}=\frac{T^k_{\!ij}-c^{\,k}_{ij}}{2}.
\mytag{4.33}
$$
Substituting \mythetag{4.33} into \mythetag{4.32}, we derive the ultimate
formula for $\Gamma^k_{\!ij}$:
$$
\hskip -2em
\gathered
\Gamma^k_{\!ij}=\sum^3_{r=0}\frac{g^{\kern 0.5pt kr}}{2}
\left(L_{\boldsymbol\Upsilon_{\!i}}\!(g_{jr})+L_{\boldsymbol
\Upsilon_{\!j}}\!(g_{ri})-L_{\boldsymbol\Upsilon_{\!r}}\!(g_{ij})
\right)-\\
-\,\frac{c^{\,k}_{ij}}{2}
+\sum^3_{r=0}\sum^3_{s=0}g^{kr}\,\frac{c^{\,s}_{i\kern 0.5pt r}}{2}\,
g_{sj}+\sum^3_{r=0}\sum^3_{s=0}g^{kr}\,\frac{c^{\,s}_{jr}}{2}
\,g_{s\kern 0.5pt i}\,+\\
+\,\frac{T^k_{ij}}{2}
-\sum^3_{r=0}\sum^3_{s=0}g^{kr}\,\frac{T^s_{i\kern 0.5pt r}}{2}
\,g_{sj}-\sum^3_{r=0}\sum^3_{s=0}g^{kr}\,\frac{T^s_{jr}}{2}
\,g_{s\kern 0.5pt i}.
\endgathered
\mytag{4.34}
$$
\mydefinition{4.6} A real spinor connection $(\Gamma,\Alpha,\bar{\Alpha})$
of the bundle of chiral spinors $SM$ is called a {\it metric connection},
if it is concordant 1) with the spin-metric tensor $\bold d$, 2) with the
Infeld-van der Waerden field $\bold G$, and 3) with the metric tensor 
$\bold d$.
\enddefinition
     The theorem~\mythetheorem{4.2} says that the conditions 1) and 2)
are sufficient for a real spinor connection $(\Gamma,\Alpha,\bar{\Alpha})$
to be a metric connection. The Einstein's theory of gravity, which is also
called the General Relativity, is a theory without torsion, i\.\,e\. the
torsion tensor $\bold T$ is taken to be zero in it: $\bold T=0$. Then
\mythetag{4.34} reduces to
$$
\hskip -2em
\gathered
\Gamma^k_{\!ij}=\sum^3_{r=0}\frac{g^{\kern 0.5pt kr}}{2}
\left(L_{\boldsymbol\Upsilon_{\!i}}\!(g_{jr})+L_{\boldsymbol
\Upsilon_{\!j}}\!(g_{ri})-L_{\boldsymbol\Upsilon_{\!r}}\!(g_{ij})
\right)-\\
-\,\frac{c^{\,k}_{ij}}{2}
+\sum^3_{r=0}\sum^3_{s=0}g^{kr}\,\frac{c^{\,s}_{i\kern 0.5pt r}}{2}\,
g_{sj}+\sum^3_{r=0}\sum^3_{s=0}g^{kr}\,\frac{c^{\,s}_{jr}}{2}
\,g_{s\kern 0.5pt i}.
\endgathered
\mytag{4.35}
$$
According to the corollary~\mythecorollary{4.1}, the $\Gamma$-components 
of a real metric spinor connection in {\it General Relativity} are the
components of the Levi-Civita connection for the metric $\bold g$ and
\mythetag{4.35} is a frame version of the well-known formula 
$$
\Gamma^k_{\!ij}=\sum^3_{r=0}\frac{g^{\kern 0.5pt kr}}{2}
\left(\frac{\partial g_{jr}}{\partial x^i}+\frac{\partial g_{ri}}
{\partial x^j}-\frac{\partial g_{ij}}{\partial x^r}\right)\!.
$$
The above calculations leading to the formula \mythetag{4.35} are standard.
They are similar to those performed in section~3 of \mycite{7}.\par
     Now let's study the concordance condition $\nabla\bold d=0$. In a
frame relative coordinate form this condition is written as follows:
$$
\hskip -2em
\nabla_{\!r}d_{ij}=
L_{\boldsymbol\Upsilon_{\!r}}\!(d_{ij})-\sum^2_{s=1}
\Alpha^{\!s}_{\kern 0.5pt ri}\,d_{\kern 0.5pt sj}
-\sum^2_{s=1}\Alpha^{\!s}_{\kern 0.5pt rj}\,d_{\kern 0.5pt is}=0.
\mytag{4.36}
$$ 
By lowering the upper index $s$ of $\Alpha^{\!s}_{\kern 0.5pt ri}$
we introduce the following quantities:
$$
\hskip -2em
\Alpha_{\kern 0.5pt rij}=\sum^2_{s=1}\Alpha^{\!s}_{\kern 0.5pt ri}\,
d_{\kern 0.5pt sj}.
\mytag{4.37}
$$
Then, due to \mythetag{4.37}, the equality \mythetag{4.36} reduces to
the following one:
$$
\hskip -2em
\Alpha_{\kern 0.5pt rij}-\Alpha_{\kern 0.5pt r\kern -0.5pt j\kern 0.5pt i}
=L_{\boldsymbol\Upsilon_{\!r}}\!(d_{ij}).
\mytag{4.38}
$$
The formula \mythetag{4.38} means that the skew-symmetric part of 
$\Alpha_{\kern 0.5pt rij}$ is determined by the Lie derivative
$L_{\boldsymbol\Upsilon_{\!r}}\!(d_{ij})$ so that we can write
$$
\hskip -2em
\Alpha_{\kern 0.5pt rij}=\hat{\Alpha}_{\kern 0.5pt rij}
+\frac{L_{\boldsymbol\Upsilon_{\!r}}\!(d_{ij})}{2},
\mytag{4.39}
$$
where $\hat{\Alpha}_{\kern 0.5pt rij}=\hat{\Alpha}_{\kern 0.5pt r
\kern -0.5pt j\kern 0.5pt i}$. Returning from \mythetag{4.39} back 
to the quantities $\Alpha^{\!s}_{\kern 0.5pt ri}$, we get
$$
\hskip -2em
\Alpha^{\!s}_{\kern 0.5pt ri}=\sum^2_{j=1}\hat{\Alpha}_{\kern 0.5pt rij}
\,d^{\kern 0.5pt js}+\sum^2_{j=1}\frac{L_{\boldsymbol\Upsilon_{\!r}}
\!(d_{ij})\ d^{\kern 0.5pt js}}{2}.
\mytag{4.40}
$$
Acting in a similar way, from $\nabla\bbd=0$ we easily derive the formulas
$$
\align
\hskip -2em
\bar{\Alpha}_{\kern 0.5pt rij}&=\Hat{\bar{\Alpha}}_{\kern 0.5pt rij}
+\frac{L_{\boldsymbol\Upsilon_{\!r}}\!(\bd_{ij})}{2},
\mytag{4.41}\\
\hskip -2em
\bar{\Alpha}\vphantom{\Alpha}^{\!s}_{\kern 0.5pt ri}
&=\sum^2_{j=1}\Hat{\bar{\Alpha}}_{\kern 0.5pt rij}
\,\bd^{\kern 0.5pt js}+\sum^2_{j=1}\frac{\,L_{\boldsymbol
\Upsilon_{\!r}}\!(\bd_{ij})\ \bd^{\kern 0.5pt js}}{2},
\mytag{4.42}
\endalign
$$
where $\Hat{\bar{\Alpha}}_{\kern 0.5pt rij}=\Hat{\bar{\Alpha}}_{\kern 0.5pt
r\kern -0.5pt j\kern 0.5pt i}$. Thus, we have managed to reduce the
concordance conditions $\nabla\bold d=0$ and $\nabla\bbd=0$ to the symmetry
conditions
$$
\xalignat 2
&\hskip -2em
\hat{\Alpha}_{\kern 0.5pt rij}=\hat{\Alpha}_{\kern 0.5pt r\kern -0.5pt j
\kern 0.5pt i},
&&\Hat{\bar{\Alpha}}_{\kern 0.5pt rij}=\Hat{\bar{\Alpha}}_{\kern 0.5pt
r\kern -0.5pt j\kern 0.5pt i}
\mytag{4.43}
\endxalignat
$$
and to the formulas \mythetag{4.40} and \mythetag{4.42} for the $\Alpha$
and $\bar{\Alpha}$-components of a spinor connection $(\Gamma,\Alpha,
\bar{\Alpha})$.\par
    The next step is to study the concordance conditions $\nabla\bold G=0$. 
Applying the formula \mythetag{4.20} to $\nabla\bold G=0$, we derive the
following equality:
$$
\hskip -2em
L_{\boldsymbol\Upsilon_{\!r}}\!(G^{i\kern 0.5pt
\bar i}_q)+\sum^2_{s=1}G^{s\kern 0.5pt\bar i}_q\,
\Alpha^i_{rs}+\sum^2_{\bar s=1}G^{i\kern 0.5pt\bar s}_q\,
\bar{\Alpha}\vphantom{\Alpha}^{\bar i}_{r\bar s}-\sum^3_{p=0}
G^{i\kern 0.5pt\bar i}_p\,\Gamma^p_{rq}=0.
\mytag{4.44}
$$
In order to transform \mythetag{4.44} we multiply it by $G^{\,q}_{j\bar j}$
and sum it over the index $q$, meanwhile taking into account the second
identity \mythetag{3.12}:
$$
\hskip -2em
2\,\Alpha^i_{\kern 0.5pt rj}\,\delta^{\,\bar i}_{\bar j}
+2\,\delta^{\,i}_{j}\,
\bar{\Alpha}\vphantom{\Alpha}^{\bar i}_{\kern 0.5pt r\bar j}
=\sum^3_{p=0}\sum^3_{q=0}
G^{i\kern 0.5pt\bar i}_p\,\Gamma^p_{rq}\,G^{\,q}_{j\bar j}
-\sum^3_{q=0}L_{\boldsymbol\Upsilon_{\!r}}\!(G^{i\kern 0.5pt
\bar i}_q)\,G^{\,q}_{j\bar j}.
\mytag{4.45}
$$
Then we apply the index lowering procedure to the indices $i$ and $\bar i$
in \mythetag{4.45}. As a result, taking into account \mythetag{4.37}, we
derive
$$
\hskip -2em
\gathered
2\,\Alpha_{\kern 0.5pt rj\kern 0.5pt i}\,\bd_{\bar j\kern 0.5pt \bar i}
+2\,d_{j\kern 0.5pt i}\,
\bar{\Alpha}_{\kern 0.5pt r\bar j\kern 0.5pt\bar i}
=\sum^2_{s=1}\sum^2_{\bar s=1}\sum^3_{p=0}\sum^3_{q=0}
d_{si}\,\bd_{\bar s\bar i}\,G^{s\bar s}_p\,
\Gamma^p_{rq}\,G^{\,q}_{j\bar j}\,-\\
-\sum^2_{s=1}\sum^2_{\bar s=1}\sum^3_{q=0}
d_{si}\,\bd_{\bar s\bar i}\,
L_{\boldsymbol\Upsilon_{\!r}}\!(G^{s\bar s}_q)\,G^{\,q}_{j\bar j}.
\endgathered
\mytag{4.46}
$$
Substituting \mythetag{4.39} and  \mythetag{4.41} into  \mythetag{4.46},
we derive
$$
\hskip -2em
\gathered
2\,\hat{\Alpha}_{\kern 0.5pt rj\kern 0.5pt i}
\,\bd_{\bar j\kern 0.5pt \bar i}+2\,d_{j\kern 0.5pt i}\,
\Hat{\bar{\Alpha}}_{\kern 0.5pt r\bar j\kern 0.5pt\bar i}
=\sum^2_{s=1}\sum^2_{\bar s=1}\sum^3_{p=0}\sum^3_{q=0}
d_{si}\,\bd_{\bar s\bar i}\,G^{s\bar s}_p\,\Gamma^p_{rq}
\,G^{\,q}_{j\bar j}\,-
\\
-\sum^2_{s=1}\sum^2_{\bar s=1}\sum^3_{q=0}
d_{si}\,\bd_{\bar s\bar i}\,
L_{\boldsymbol\Upsilon_{\!r}}\!(G^{s\bar s}_q)\,G^{\,q}_{j\bar j}
-L_{\boldsymbol\Upsilon_{\!r}}\!(d_{j\kern 0.5pt i})\,
\bd_{\bar j\kern 0.5pt\bar i}
-L_{\boldsymbol\Upsilon_{\!r}}\!(\bd_{\bar j\kern 0.5pt\bar i})\,
d_{j\kern 0.5pt i}.
\endgathered
\mytag{4.47}
$$
The left hand side of the equality \mythetag{4.47} is a sum of two
terms. The first term is symmetric in $i$ and $j$, see \mythetag{4.43},
while the other is skew-symmetric in these indices. Therefore, if we
subdivide the right hand side of \mythetag{4.47} into symmetric and
skew-symmetric parts, we can write \mythetag{4.47} as two separate
equalities:
$$
\gather
\hskip -2em
\gathered
\hat{\Alpha}_{\kern 0.5pt rj\kern 0.5pt i}
\,\bd_{\bar j\kern 0.5pt \bar i}=\sum^2_{s=1}\sum^2_{\bar s=1}
\sum^3_{p=0}\sum^3_{q=0}\frac{
d_{si}\,\bd_{\bar s\bar i}\,G^{s\bar s}_p\,\Gamma^p_{rq}
\,G^{\,q}_{j\bar j}+
d_{sj}\,\bd_{\bar s\bar i}\,G^{s\bar s}_p\,\Gamma^p_{rq}
\,G^{\,q}_{i\bar j}}{4}\,-\\
-\sum^2_{s=1}\sum^2_{\bar s=1}\sum^3_{q=0}
\frac{d_{si}\,\bd_{\bar s\bar i}\,
L_{\boldsymbol\Upsilon_{\!r}}\!(G^{s\bar s}_q)\,G^{\,q}_{j\bar j}
+d_{sj}\,\bd_{\bar s\bar i}\,
L_{\boldsymbol\Upsilon_{\!r}}\!(G^{s\bar s}_q)\,G^{\,q}_{i\bar j}}{4},
\endgathered\qquad
\mytag{4.48}\\
\vspace{3ex}
\gathered
\hskip -2em
d_{j\kern 0.5pt i}\,
\Hat{\bar{\Alpha}}_{\kern 0.5pt r\bar j\kern 0.5pt\bar i}
=\sum^2_{s=1}\sum^2_{\bar s=1}
\sum^3_{p=0}\sum^3_{q=0}\frac{
d_{si}\,\bd_{\bar s\bar i}\,G^{s\bar s}_p\,\Gamma^p_{rq}
\,G^{\,q}_{j\bar j}-
d_{sj}\,\bd_{\bar s\bar i}\,G^{s\bar s}_p\,\Gamma^p_{rq}
\,G^{\,q}_{i\bar j}}{4}\,-\\
-\sum^2_{s=1}\sum^2_{\bar s=1}\sum^3_{q=0}
\frac{d_{si}\,\bd_{\bar s\bar i}\,
L_{\boldsymbol\Upsilon_{\!r}}\!(G^{s\bar s}_q)\,G^{\,q}_{j\bar j}
-d_{sj}\,\bd_{\bar s\bar i}\,
L_{\boldsymbol\Upsilon_{\!r}}\!(G^{s\bar s}_q)\,G^{\,q}_{i\bar j}}{4}\,-
\\
\vphantom{\sum^2_{s=1}}
-\frac{L_{\boldsymbol\Upsilon_{\!r}}\!(d_{j\kern 0.5pt i})\,
\bd_{\bar j\kern 0.5pt\bar i}
+L_{\boldsymbol\Upsilon_{\!r}}\!(\bd_{\bar j\kern 0.5pt\bar i})\,
d_{j\kern 0.5pt i}}{2}.
\endgathered\qquad
\mytag{4.49}
\endgather
$$
In two-dimensional case any equality skew-symmetric in two indices is
equivalent to a scalar equality independent of these two indices. In 
the case of the equality \mythetag{4.49} we can multiply it by 
$d^{\kern 0.5pt ij}$ and sum over the indices $i$ and $j$. As a result 
we obtain the following equality equivalent to \mythetag{4.49}:
$$
\hskip -2em
\gathered
\Hat{\bar{\Alpha}}_{\kern 0.5pt r\bar j\kern 0.5pt\bar i}
=-\frac{L_{\boldsymbol\Upsilon_{\!r}}\!(\bd_{\bar j\kern 0.5pt
\bar i})}{2}+\sum^2_{s=1}\sum^2_{\bar s=1}
\sum^3_{p=0}\sum^3_{q=0}\frac{\bd_{\bar s\bar i}\,G^{s\bar s}_p
\,\Gamma^p_{rq}\,G^{\,q}_{s\bar j}}{4}\,-\\
-\sum^2_{s=1}\sum^2_{\bar s=1}\sum^3_{q=0}
\frac{\bd_{\bar s\bar i}\,L_{\boldsymbol\Upsilon_{\!r}}\!(G^{s\bar s}_q)
\,G^{\,q}_{s\bar j}}{4}
-\sum^2_{i=1}\sum^2_{j=1}
\frac{L_{\boldsymbol\Upsilon_{\!r}}\!(d_{j\kern 0.5pt i})
\,d^{\kern 0.5pt ij}\,\bd_{\bar j\kern 0.5pt\bar i}}{4}.
\endgathered
\mytag{4.50}
$$
By definition the left hand side of \mythetag{4.50} is symmetric in
$\bar i$ and $\bar j$. Therefore it should be equal to the symmetric 
part of the right hand side
$$
\hskip -2em
\gathered
\Hat{\bar{\Alpha}}_{\kern 0.5pt r\bar j\kern 0.5pt\bar i}
=\sum^2_{s=1}\sum^2_{\bar s=1}\sum^3_{p=0}\sum^3_{q=0}\frac{\bd_{\bar s
\bar i}\,G^{s\bar s}_p\,\Gamma^p_{rq}\,G^{\,q}_{s\bar j}+\bd_{\bar s
\bar j}\,G^{s\bar s}_p\,\Gamma^p_{rq}\,G^{\,q}_{s\bar i}}{8}\,-\\
-\sum^2_{s=1}\sum^2_{\bar s=1}\sum^3_{q=0}\frac{\bd_{\bar s\bar i}
\,L_{\boldsymbol\Upsilon_{\!r}}\!(G^{s\bar s}_q)\,G^{\,q}_{s\bar j}+
\bd_{\bar s\bar j}\,L_{\boldsymbol\Upsilon_{\!r}}\!(G^{s\bar s}_q)
\,G^{\,q}_{s\bar i}}{8},
\endgathered
\mytag{4.51}
$$
while the skew-symmetric part of the right hand side of \mythetag{4.50}
should be zero:
$$
\hskip -2em
\gathered
\sum^2_{s=1}\sum^2_{\bar s=1}\sum^3_{p=0}\sum^3_{q=0}\frac{\bd_{\bar s
\bar i}\,G^{s\bar s}_p\,\Gamma^p_{rq}\,G^{\,q}_{s\bar j}-\bd_{\bar s
\bar j}\,G^{s\bar s}_p\,\Gamma^p_{rq}\,G^{\,q}_{s\bar i}}{8}\,-\\
-\sum^2_{s=1}\sum^2_{\bar s=1}\sum^3_{q=0}\frac{\bd_{\bar s\bar i}
\,L_{\boldsymbol\Upsilon_{\!r}}\!(G^{s\bar s}_q)\,G^{\,q}_{s\bar j}-
\bd_{\bar s\bar j}\,L_{\boldsymbol\Upsilon_{\!r}}\!(G^{s\bar s}_q)
\,G^{\,q}_{s\bar i}}{8}\,-\\
-\frac{L_{\boldsymbol\Upsilon_{\!r}}\!(\bd_{\bar j\kern 0.5pt
\bar i})}{2}-\sum^2_{i=1}\sum^2_{j=1}
\frac{L_{\boldsymbol\Upsilon_{\!r}}\!(d_{j\kern 0.5pt i})
\,d^{\kern 0.5pt ij}\,\bd_{\bar j\kern 0.5pt\bar i}}{4}=0.
\endgathered
\mytag{4.52}
$$
Again, using the feature of the two-dimensional case, we can reduce
\mythetag{4.52} to an equality independent of $\bar i$ and $\bar j$.
For this purpose let's multiply it by $\bd^{\kern 0.5pt\bar i\bar j}$
and sum over the indices $\bar i$ and $\bar j$. As a result we get
$$
\hskip -2em
\gathered
\sum^2_{s=1}\sum^2_{\bar s=1}\sum^3_{p=0}\sum^3_{q=0}
\frac{G^{s\bar s}_p\,\Gamma^p_{rq}\,G^{\,q}_{s\bar s}}{4}
-\sum^2_{s=1}\sum^2_{\bar s=1}\sum^3_{q=0}\frac{L_{\boldsymbol
\Upsilon_{\!r}}\!(G^{s\bar s}_q)\,G^{\,q}_{s\bar s}}{4}\,-\\
-\sum^2_{\bar i=1}\sum^2_{\bar j=1}
\frac{L_{\boldsymbol\Upsilon_{\!r}}\!(\bd_{\bar j\kern 0.5pt
\bar i})\,\bd^{\kern 0.5pt\bar i\bar j}}{2}-\sum^2_{i=1}
\sum^2_{j=1}\frac{L_{\boldsymbol\Upsilon_{\!r}}\!(d_{j\kern 0.5pt i})
\,d^{\kern 0.5pt ij}}{2}=0.
\endgathered
\mytag{4.53}
$$
Taking into account the first identity \mythetag{3.12} and the formula
\mythetag{3.2}, we can transform the equality \mythetag{4.53} to a more
symmetric form:
$$
\hskip -2em
\gathered
\sum^3_{p=0}\frac{\Gamma^p_{rp}}{2}
-\sum^2_{s=1}\sum^2_{\bar s=1}\sum^3_{q=0}\frac{L_{\boldsymbol
\Upsilon_{\!r}}\!(G^{s\bar s}_q)\,d_{ns}\,\bd_{\bar n\bar s}
\,g^{qp}\,G^{n\bar n}_p}{4}\,-\\
-\sum^2_{\bar i=1}\sum^2_{\bar j=1}
\frac{L_{\boldsymbol\Upsilon_{\!r}}\!(\bd_{\bar j\kern 0.5pt
\bar i})\,\bd^{\kern 0.5pt\bar i\bar j}}{2}-\sum^2_{i=1}
\sum^2_{j=1}\frac{L_{\boldsymbol\Upsilon_{\!r}}\!(d_{j\kern 0.5pt i})
\,d^{\kern 0.5pt ij}}{2}=0.
\endgathered
\mytag{4.54}
$$
Note that the product $d_{ns}\,\bd_{\bar n\bar s}\,g^{qp}$ in 
\mythetag{4.54} is invariant under the simultaneous transposition
of $s\longleftrightarrow n$, $\bar s\longleftrightarrow\bar n$,
and $p\longleftrightarrow q$. Therefore, we can write
$$
\gather
\sum^3_{p=0}\frac{\Gamma^p_{rp}}{2}
-\sum^2_{\bar i=1}\sum^2_{\bar j=1}
\frac{L_{\boldsymbol\Upsilon_{\!r}}\!(\bd_{\bar j\kern 0.5pt
\bar i})\,\bd^{\kern 0.5pt\bar i\bar j}}{2}-\sum^2_{i=1}
\sum^2_{j=1}\frac{L_{\boldsymbol\Upsilon_{\!r}}\!(d_{j\kern 0.5pt i})
\,d^{\kern 0.5pt ij}}{2}\,-\\
-\sum^2_{s=1}\sum^2_{\bar s=1}\sum^3_{q=0}
\sum^2_{n=1}\sum^2_{\bar n=1}\sum^3_{p=0}\frac{L_{\boldsymbol
\Upsilon_{\!r}}\!(G^{s\bar s}_q)\,d_{ns}\,\bd_{\bar n\bar s}
\,g^{qp}\,G^{n\bar n}_p+G^{s\bar s}_q\,d_{ns}\,\bd_{\bar n\bar s}
\,g^{qp}\,L_{\boldsymbol\Upsilon_{\!r}}\!(G^{n\bar n}_p)}{8}=0.
\endgather
$$
The Lie derivative $L_{\boldsymbol\Upsilon_{\!r}}$ in the above equality
acts as a first order linear differential operator. For this reason we
can continue transforming the above equality:
$$
\allowdisplaybreaks
\gather
\sum^3_{p=0}\frac{\Gamma^p_{rp}}{2}
-\sum^2_{\bar i=1}\sum^2_{\bar j=1}
\frac{L_{\boldsymbol\Upsilon_{\!r}}\!(\bd_{\bar j\kern 0.5pt
\bar i})\,\bd^{\kern 0.5pt\bar i\bar j}}{2}-\sum^2_{i=1}
\sum^2_{j=1}\frac{L_{\boldsymbol\Upsilon_{\!r}}\!(d_{j\kern 0.5pt i})
\,d^{\kern 0.5pt ij}}{2}\,-\\
-\sum^2_{s=1}\sum^2_{\bar s=1}\sum^3_{q=0}
\sum^2_{n=1}\sum^2_{\bar n=1}\sum^3_{p=0}\frac{L_{\boldsymbol
\Upsilon_{\!r}}\!(G^{s\bar s}_q\,d_{ns}\,\bd_{\bar n\bar s}
\,G^{n\bar n}_p)\,g^{qp}}{8}\,+\\
+\sum^2_{s=1}\sum^2_{\bar s=1}\sum^3_{q=0}
\sum^2_{n=1}\sum^2_{\bar n=1}\sum^3_{p=0}\frac{G^{s\bar s}_q
\,L_{\boldsymbol\Upsilon_{\!r}}\!(d_{ns})\,\bd_{\bar n\bar s}
\,g^{qp}\,G^{n\bar n}_p+G^{s\bar s}_q\,d_{ns}
\,L_{\boldsymbol\Upsilon_{\!r}}\!(\bd_{\bar n\bar s})\,g^{qp}
\,G^{n\bar n}_p}{8}=0.
\endgather
$$
Taking into account the second identity \mythetag{3.12} and the formula
\mythetag{3.2}, we find that the last term of the above equality cancels
the second and the third terms in it. As a result, using \mythetag{3.9},
we can write it as follows:
$$
\hskip -2em
\sum^3_{p=0}\frac{\Gamma^p_{rp}}{2}
-\sum^3_{p=0}\sum^3_{q=0}\frac{L_{\boldsymbol\Upsilon_{\!r}}\!(g_{qp})
\,g^{qp}}{4}=0.
\mytag{4.55}
$$
Using the formula \mythetag{4.20}, we can write \mythetag{4.55} in a 
very simple form:
$$
\hskip -2em
\sum^3_{p=0}\sum^3_{q=0}g^{qp}\ \nabla_{\!r}g_{qp}=0.
\mytag{4.56}
$$
Thus, the formula \mythetag{4.52} is reduced to \mythetag{4.56}, 
while \mythetag{4.51} is equivalent to \mythetag{4.50} provided
\mythetag{4.56} is fulfilled. From \mythetag{4.50} and \mythetag{4.41}
we derive the following expression for $\bar{\Alpha}$-components of 
the spinor connection $(\Gamma,\Alpha,\bar{\Alpha})$:
$$
\hskip -2em
\gathered
\bar{\Alpha}\vphantom{\Alpha}^{\bar i}_{\kern 0.5pt r\bar j}
=\sum^2_{s=1}
\sum^3_{p=0}\sum^3_{q=0}\frac{G^{s\bar i}_p
\,\Gamma^p_{rq}\,G^{\,q}_{s\bar j}}{4}\,-\\
-\sum^2_{s=1}\sum^3_{q=0}
\frac{L_{\boldsymbol\Upsilon_{\!r}}\!(G^{s\bar i}_q)
\,G^{\,q}_{s\bar j}}{4}
-\sum^2_{i=1}\sum^2_{j=1}
\frac{L_{\boldsymbol\Upsilon_{\!r}}\!(d_{j\kern 0.5pt i})
\,d^{\kern 0.5pt ij}\,\delta^{\,\bar i}_{\bar j}}{4}.
\endgathered
\mytag{4.57}
$$
Both \mythetag{4.56} and \mythetag{4.57}, when taken together, are 
equivalent to \mythetag{4.49}.\par
     Having all done with \mythetag{4.49}, now we return back to 
the equality \mythetag{4.48}. The left hand side of this equality 
is skew-symmetric in $\bar i$ and $\bar j$. Like in the case of
\mythetag{4.47}, subdividing the right hand side of \mythetag{4.48}
into two parts symmetric and skew-symmetric in $\bar i$ and $\bar j$,
we write \mythetag{4.48} as two separate equalities. Here is the first 
of these two equalities. It is skew-symmetric in $\bar i$ and $\bar j$:
$$
\gathered
\hat{\Alpha}_{\kern 0.5pt rj\kern 0.5pt i}
\,\bd_{\bar j\kern 0.5pt \bar i}
=\sum^2_{s=1}\sum^2_{\bar s=1}
\sum^3_{p=0}\sum^3_{q=0}\frac{
d_{si}\,\bd_{\bar s\bar i}\,G^{s\bar s}_p\,\Gamma^p_{rq}
\,G^{\,q}_{j\bar j}+
d_{sj}\,\bd_{\bar s\bar i}\,G^{s\bar s}_p\,\Gamma^p_{rq}
\,G^{\,q}_{i\bar j}}{8}\,-\\
-\sum^2_{s=1}\sum^2_{\bar s=1}
\sum^3_{p=0}\sum^3_{q=0}\frac{
d_{si}\,\bd_{\bar s\bar j}\,G^{s\bar s}_p\,\Gamma^p_{rq}
\,G^{\,q}_{j\kern 0.5pt\bar i}+
d_{sj}\,\bd_{\bar s\bar j}\,G^{s\bar s}_p\,\Gamma^p_{rq}
\,G^{\,q}_{i\kern 0.5pt\bar i}}{8}\,-\\
-\sum^2_{s=1}\sum^2_{\bar s=1}\sum^3_{q=0}
\frac{d_{si}\,\bd_{\bar s\bar i}\,
L_{\boldsymbol\Upsilon_{\!r}}\!(G^{s\bar s}_q)\,G^{\,q}_{j\bar j}
+d_{sj}\,\bd_{\bar s\bar i}\,
L_{\boldsymbol\Upsilon_{\!r}}\!(G^{s\bar s}_q)\,G^{\,q}_{i\bar j}}{8}\,+\\
+\sum^2_{s=1}\sum^2_{\bar s=1}\sum^3_{q=0}
\frac{d_{si}\,\bd_{\bar s\bar j}\,
L_{\boldsymbol\Upsilon_{\!r}}\!(G^{s\bar s}_q)\,G^{\,q}_{j\kern 0.5pt\bar i}
+d_{sj}\,\bd_{\bar s\bar j}\,
L_{\boldsymbol\Upsilon_{\!r}}\!(G^{s\bar s}_q)\,G^{\,q}_{i\kern 0.5pt
\bar i}}{8}.
\endgathered\qquad
\mytag{4.58}
$$
Due to the skew-symmetry \mythetag{4.58} can be transformed to
an equality independent of $\bar i$ and $\bar j$ at all. Multiplying it 
by $\bd^{\kern 0.5pt\bar i\bar j}$ \pagebreak and summing over $\bar i$ 
and $\bar j$, we get the following equality analogous to the equality 
\mythetag{4.51}:
$$
\hskip -2em
\gathered
\hat{\Alpha}_{\kern 0.5pt rj\kern 0.5pt i}
=\sum^2_{s=1}\sum^2_{\bar s=1}\sum^3_{p=0}\sum^3_{q=0}
\frac{d_{si}\,G^{s\bar s}_p\,\Gamma^p_{rq}
\,G^{\,q}_{j\kern 0.5pt\bar s}+d_{sj}\,G^{s\bar s}_p\,
\Gamma^p_{rq}\,G^{\,q}_{i\kern 0.5pt\bar s}}{8}\,-\\
-\sum^2_{s=1}\sum^2_{\bar s=1}\sum^3_{q=0}
\frac{d_{si}\,L_{\boldsymbol\Upsilon_{\!r}}\!(G^{s\bar s}_q)
\,G^{\,q}_{j\kern 0.5pt\bar s}+d_{sj}\,L_{\boldsymbol
\Upsilon_{\!r}}\!(G^{s\bar s}_q)\,G^{\,q}_{i\kern 0.5pt\bar s}}{8}.
\endgathered
\mytag{4.59}
$$
By analogy to \mythetag{4.52} one can write the following equality:
$$
\hskip -2em
\gathered
\sum^2_{s=1}\sum^2_{\bar s=1}\sum^3_{p=0}\sum^3_{q=0}
\frac{d_{si}\,G^{s\bar s}_p\,\Gamma^p_{rq}\,
G^{\,q}_{j\kern 0.5pt\bar s}-d_{sj}\,G^{s\bar s}_p\,\Gamma^p_{rq}
\,G^{\,q}_{i\kern 0.5pt\bar s}}{8}\,-\\
-\sum^2_{s=1}\sum^2_{\bar s=1}\sum^3_{q=0}\frac{d_{si}
\,L_{\boldsymbol\Upsilon_{\!r}}\!(G^{s\bar s}_q)
\,G^{\,q}_{j\kern 0.5pt\bar s}-
d_{sj}\,L_{\boldsymbol\Upsilon_{\!r}}\!(G^{s\bar s}_q)
\,G^{\,q}_{i\kern 0.5pt\bar s}}{8}\,-\\
-\frac{L_{\boldsymbol\Upsilon_{\!r}}\!(d_{j\kern 0.5pt i})}{2}
-\sum^2_{\bar i=1}\sum^2_{\bar j=1}
\frac{L_{\boldsymbol\Upsilon_{\!r}}\!(\bd_{\bar j\kern 0.5pt\bar i})
\,\bd^{\kern 0.5pt\bar i\bar j}\,d_{j\kern 0.5pt i}}{4}=0.
\endgathered
\mytag{4.60}
$$
Acting in a similar way as in the case of \mythetag{4.52}, one can
show that the equality \mythetag{4.60} is equivalent to \mythetag{4.56}.
Adding \mythetag{4.60} to \mythetag{4.59}, we get
$$
\hskip -2em
\gathered
\hat{\Alpha}_{\kern 0.5pt rj\kern 0.5pt i}
=-\frac{L_{\boldsymbol\Upsilon_{\!r}}\!(d_{j\kern 0.5pt i})}{2}
+\sum^2_{s=1}\sum^2_{\bar s=1}\sum^3_{p=0}\sum^3_{q=0}
\frac{d_{si}\,G^{s\bar s}_p\,\Gamma^p_{rq}\,
G^{\,q}_{j\kern 0.5pt\bar s}}{4}\,-\\
-\sum^2_{s=1}\sum^2_{\bar s=1}\sum^3_{q=0}
\frac{d_{si}\,L_{\boldsymbol\Upsilon_{\!r}}\!(G^{s\bar s}_q)
\,G^{\,q}_{j\kern 0.5pt\bar s}}{4}
-\sum^2_{\bar i=1}\sum^2_{\bar j=1}
\frac{L_{\boldsymbol\Upsilon_{\!r}}\!(\bd_{\bar j\kern 0.5pt\bar i})
\,\bd^{\kern 0.5pt\bar i\bar j}\,d_{j\kern 0.5pt i}}{4}.
\endgathered
\mytag{4.61}
$$
The formula \mythetag{4.61} is similar to \mythetag{4.50}. From
\mythetag{4.61} and \mythetag{4.39} we derive
$$
\hskip -2em
\gathered
\Alpha^i_{\kern 0.5pt r j}
=\sum^2_{\bar s=1}\sum^3_{p=0}\sum^3_{q=0}
\frac{G^{i\bar s}_p\,\Gamma^p_{rq}
\,G^{\,q}_{j\kern 0.5pt\bar s}}{4}\,-\\
-\sum^2_{\bar s=1}\sum^3_{q=0}
\frac{L_{\boldsymbol\Upsilon_{\!r}}\!(G^{i\bar s}_q)
\,G^{\,q}_{j\kern 0.5pt\bar s}}{4}
-\sum^2_{\bar i=1}\sum^2_{\bar j=1}
\frac{L_{\boldsymbol\Upsilon_{\!r}}\!(\bd_{\bar j\kern 0.5pt\bar i})
\,\bd^{\kern 0.5pt\bar i\bar j}\,\delta^{\,i}_j}{4}.
\endgathered
\mytag{4.62}
$$\par
    The last step is to study the symmetric part of the equality 
\mythetag{4.48}. Symmetrizing \mythetag{4.48} with respect to $\bar i$ 
and $\bar j$, we find that its symmetric part is written as
$$
\gather
\hskip -2em
B_{ij\kern 0.5pt\bar i\bar j}+B_{j\kern 0.5pt i\kern 0.5pt\bar i\bar j}=0,
\mytag{4.63}\\
\vspace{-2ex}
\intertext{where}
\vspace{-2ex}
\hskip -2em
\gathered
B_{ij\kern 0.5pt\bar i\bar j}=\sum^2_{s=1}\sum^2_{\bar s=1}
\sum^3_{p=0}\sum^3_{q=0}\frac{
d_{si}\,\bd_{\bar s\bar i}\,G^{s\bar s}_p\,\Gamma^p_{rq}
\,G^{\,q}_{j\bar j}
+d_{sj}\,\bd_{\bar s\bar j}\,G^{s\bar s}_p\,\Gamma^p_{rq}
\,G^{\,q}_{i\kern 0.5pt\bar i}}{8}\,-\\
-\sum^2_{s=1}\sum^2_{\bar s=1}\sum^3_{q=0}
\frac{d_{si}\,\bd_{\bar s\bar i}\,
L_{\boldsymbol\Upsilon_{\!r}}\!(G^{s\bar s}_q)\,G^{\,q}_{j\bar j}
+d_{sj}\,\bd_{\bar s\bar j}\,
L_{\boldsymbol\Upsilon_{\!r}}\!(G^{s\bar s}_q)\,G^{\,q}_{i\kern 0.5pt
\bar i}}{8}.
\endgathered\qquad
\mytag{4.64}
\endgather
$$
Using \mythetag{4.25} and the formula \mythetag{3.2}, we transform
\mythetag{4.64} as follows:
$$
\gathered
B_{ij\kern 0.5pt\bar i\bar j}
=\msum{2}_{s,\ \bar s,\ n,\ \bar n}\msum{3}_{p,\ q,\ \alpha,\ \beta}
\frac{(d_{si}\,\bd_{\bar s\bar i}\,G^{s\bar s}_p\,
g^{p\alpha})\,\Gamma_{r\beta\alpha}\,(d_{nj}\,\bd_{\bar n\bar j}\,
G^{n\bar n}_q\,g^{\beta q})}{8}\,+\\
+\msum{2}_{s,\ \bar s,\ n,\ \bar n}\msum{3}_{p,\ q,\ \alpha,\ \beta}
\frac{(d_{sj}\,\bd_{\bar s\bar j}\,G^{s\bar s}_p\,
g^{p\alpha})\,\Gamma_{r\beta\alpha}\,(d_{ni}\,\bd_{\bar n\bar i}\,
G^{n\bar n}_q\,g^{\beta q})}{8}\,-\\
-\msum{2}_{s,\ \bar s,\ n,\ \bar n}\sum^3_{p,\ q}
\frac{(d_{si}\,\bd_{\bar s\bar i}\,
L_{\boldsymbol\Upsilon_{\!r}}\!(G^{s\bar s}_q))\,g^{qp}\,
(d_{nj}\,\bd_{\bar n\bar j}\,G^{n\bar n}_p)}{8}\,-\\
-\msum{2}_{s,\ \bar s,\ n,\ \bar n}\sum^3_{p,\ q}
\frac{(d_{si}\,\bd_{\bar s\bar i}\,
G^{s\bar s}_q\,g^{qp}\,
(d_{nj}\,\bd_{\bar n\bar j}
\,L_{\boldsymbol\Upsilon_{\!r}}\!(G^{n\bar n}_p))}{8}.
\endgathered\qquad
\mytag{4.65}
$$
Remember that the Lie derivative $L_{\boldsymbol\Upsilon_{\!r}}$ acts
as a first order linear differential operator. Therefore, the formula
\mythetag{4.65} can be written as
$$
\gathered
B_{ij\kern 0.5pt\bar i\bar j}=
-\msum{2}_{s,\ \bar s,\ n,\ \bar n}\sum^3_{p,\ q}
\frac{d_{si}\,\bd_{\bar s\bar i}\,d_{nj}\,\bd_{\bar n\bar j}\,
L_{\boldsymbol\Upsilon_{\!r}}\!(G^{s\bar s}_q\,g^{qp}\,
G^{n\bar n}_p)}{8}\,+\\
+\msum{2}_{s,\ \bar s,\ n,\ \bar n}\sum^3_{p,\ q}
\frac{(d_{si}\,\bd_{\bar s\bar i}\,G^{s\bar s}_q)\,
L_{\boldsymbol\Upsilon_{\!r}}\!(g^{qp})\,
(d_{nj}\,\bd_{\bar n\bar j}\,G^{n\bar n}_p)}{8}\,+\\
+\msum{2}_{s,\ \bar s,\ n,\ \bar n}\msum{3}_{p,\ q,\ \alpha,\ \beta}
\frac{(d_{si}\,\bd_{\bar s\bar i}\,G^{s\bar s}_p\,
g^{p\alpha})\,\Gamma_{r\beta\alpha}\,(d_{nj}\,\bd_{\bar n\bar j}\,
G^{n\bar n}_q\,g^{\beta q})}{8}\,+\\
+\msum{2}_{s,\ \bar s,\ n,\ \bar n}\msum{3}_{p,\ q,\ \alpha,\ \beta}
\frac{(d_{sj}\,\bd_{\bar s\bar j}\,G^{s\bar s}_p\,
g^{p\alpha})\,\Gamma_{r\beta\alpha}\,(d_{ni}\,\bd_{\bar n\bar i}\,
G^{n\bar n}_q\,g^{\beta q})}{8}.
\endgathered\qquad
\mytag{4.66}
$$
Remember that $g^{\kern 0.5pt qp}$ form the inverse matrix for 
$g_{qp}$. Therefore, we have
$$
\hskip -2em
L_{\boldsymbol\Upsilon_{\!r}}\!(g^{qp})=
-\sum^3_{\alpha=0}\sum^3_{\beta=0}g^{p\kern 0.5pt\alpha}\,
L_{\boldsymbol\Upsilon_{\!r}}\!(g_{\alpha\beta})\,g^{\beta q}
\mytag{4.67}
$$
Substituting \mythetag{4.67} into \mythetag{4.66} and applying 
\mythetag{3.10} and \mythetag{4.20} to it, we derive
$$
\gathered
B_{ij\kern 0.5pt\bar i\bar j}=
-\msum{2}_{s,\ \bar s,\ n,\ \bar n}\sum^3_{p,\ q}
\frac{d_{si}\,\bd_{\bar s\bar i}\,d_{nj}\,\bd_{\bar n\bar j}\,
L_{\boldsymbol\Upsilon_{\!r}}\!(d^{sn}\,\bd^{\bar s\bar n})}{4}\,-\\
-\msum{2}_{s,\ \bar s,\ n,\ \bar n}\msum{3}_{p,\ q,\ \alpha,\ \beta}
\frac{(d_{si}\,\bd_{\bar s\bar i}\,G^{s\bar s}_p\,
g^{p\alpha})\,\nabla_{\!r}g_{\alpha\beta}\,(d_{nj}\,\bd_{\bar n\bar j}\,
G^{n\bar n}_q\,g^{\beta q})}{8}.
\endgathered\qquad
\mytag{4.68}
$$
The first term in \mythetag{4.68} is skew-symmetric in $i$ and $j$.
It vanishes when we substitute \mythetag{4.68} into \mythetag{4.63}.
As a result the equality \mythetag{4.63} takes the form
$$
\pagebreak
\hskip -2em
\sum^3_{\alpha=0}\sum^3_{\beta=0}
G^\alpha_{i\kern 0.5pt\bar i}\ \nabla_{\!r}g_{\alpha\beta}
\ G^\beta_{j\bar j}+\sum^3_{\alpha=0}\sum^3_{\beta=0}
G^\alpha_{j\kern 0.5pt\bar i}\ \nabla_{\!r}g_{\alpha\beta}
\ G^\beta_{i\bar j}=0.
\mytag{4.69}
$$
Thus the concordance condition $\nabla\bold G=0$ is equivalent
to the formulas \mythetag{4.57}, and \mythetag{4.62} provided
$\nabla\bold d=0$, $\nabla\bbd=0$, and the equalities \mythetag{4.56} 
and \mythetag{4.69} are fulfilled. This result of the above 
calculations can be stated as the following theorem.
\mytheorem{4.3} The concordance conditions $\nabla\bold d=0$,
$\nabla\bbd=0$, and $\nabla\bold G=0$ for a spinor connection
$(\Gamma,\Alpha,\bar{\Alpha})$ are equivalent to the formulas
\mythetag{4.34}, \mythetag{4.57}, and \mythetag{4.62} for its
components in an arbitrary frame pair $(U,\,\boldsymbol\Upsilon_0,
\,\boldsymbol\Upsilon_1,\,\boldsymbol\Upsilon_2,\,\boldsymbol
\Upsilon_3)$ and $(U,\,\boldsymbol\Psi_1,\,\boldsymbol\Psi_2)$
of the bundles $TM$ and $SM$ respectively.
\endproclaim
\demo{Proof} Note that $\nabla\bold d=0$, $\nabla\bbd=0$, and 
$\nabla\bold G=0$ imply $\nabla\bold g=0$. The latter concordance
condition is equivalent to the formula \mythetag{4.34}. Moreover,
due to this condition the equalities \mythetag{4.56} and 
\mythetag{4.69} are fulfilled. Then from $\nabla\bold d=0$, 
$\nabla\bbd=0$, and $\nabla\bold G=0$ the formulas \mythetag{4.57} 
and \mythetag{4.62} are derived.\par
    Conversely, the formula \mythetag{4.34} leads to $\nabla\bold g=0$
and to the equalities \mythetag{4.56} and \mythetag{4.69}. The equality
\mythetag{4.69} is equivalent to \mythetag{4.63}. The equality 
\mythetag{4.56} is equivalent to \mythetag{4.60}. Being combined with
\mythetag{4.62}, the  equality \mythetag{4.60} leads to \mythetag{4.58}
and \mythetag{4.42}. Both \mythetag{4.63} and \mythetag{4.58} yield 
\mythetag{4.48}.\par
    The equality \mythetag{4.56} is equivalent to \mythetag{4.52}.
Being combined with \mythetag{4.57}, the  equality \mythetag{4.52} leads 
to \mythetag{4.49} and \mythetag{4.40}. The equalities \mythetag{4.40}
and \mythetag{4.42} are equivalent to $\nabla\bold d=0$ and $\nabla\bbd=0$.
An finally, from \mythetag{4.48} and \mythetag{4.49} by applying
\mythetag{4.40} and \mythetag{4.42} we derive $\nabla\bold G=0$. 
The theorem is proved.\qed\enddemo
\mycorollary{4.2} The components of a real metric connection $(\Gamma,
\Alpha,\bar{\Alpha})$ with zero torsion $\bold T=0$ for the bundle of
chiral spinors $SM$ are given by the explicit formulas \mythetag{4.35},
\mythetag{4.57}, and \mythetag{4.62}.
\endproclaim
\head
5. Dirac spinors.
\endhead
    Let $M$ be a space-time manifold and let $SM$ be a spinor bundle
over $M$ introduced by the definition~\mythedefinition{1.1}. By 
$S^{\sssize\dagger}\!M$ (see \mythetag{2.2}) we denote the Hermitian
conjugate bundle for $SM$. Taking both $SM$ and $S^{\sssize\dagger}\!M$, 
we construct their direct sum
$$
\hskip -2em
DM=SM\oplus S^{\sssize\dagger}\!M.
\mytag{5.1}
$$
The direct sum \mythetag{5.1} is called the {\it Dirac bundle} associated
with the spinor bundle $SM$. This is a four-dimensional complex bundle
over $M$. The bundles $SM$ and $S^{\sssize\dagger}\!M$, when treated as
the constituents of $DM$, are called {\it chiral bundles}.\par
     Due to the expansion the Dirac bundle $DM$ acquires from $SM$ its 
three basic spin-tensorial fields: the {\it spin-metric tensor\/} $\bold
d$, the {\it chirality operator\/} $\bold H$, and the {\it Hermitian
spin-metric tensor\/} $\bold D$, which is also called the {\it Dirac form}.
The definitions of these three fields can be found in section~3 of
\mycite{6}.
\mydefinition{5.1}A frame $(U,\,\boldsymbol\Psi_1,\,\boldsymbol\Psi_2,
\,\boldsymbol\Psi_3,\,\boldsymbol\Psi_4)$ of the Dirac bundle $DM$ is 
called an {\it orthonormal frame\/} if the spin-metric tensor $\bold d$ 
is represented by the following skew-symmetric matrix in this frame:
$$
\hskip -2em
d_{ij}=
d(\boldsymbol\Psi_i,\boldsymbol\Psi_j)
=\Vmatrix 0 & 1 & 0 & 0\\-1 & 0 & 0 & 0\\
0 & 0 & 0 & -1\\0 & 0 & 1 & 0\endVmatrix.
\mytag{5.2}
$$
\enddefinition
\mydefinition{5.2}A frame $(U,\,\boldsymbol\Psi_1,\,\boldsymbol\Psi_2,
\,\boldsymbol\Psi_3,\,\boldsymbol\Psi_4)$ of the Dirac bundle $DM$ is 
called an {\it anti-orthonormal frame\/} if the spin-metric tensor 
$\bold d$ is represented by the matrix opposite to the matrix
\mythetag{5.2} in this frame:
$$
\hskip -2em
d_{ij}=
d(\boldsymbol\Psi_i,\boldsymbol\Psi_j)
=\Vmatrix 0 & -1 & 0 & 0\\1 & 0 & 0 & 0\\
0 & 0 & 0 & 1\\0 & 0 & -1 & 0\endVmatrix.
\mytag{5.3}
$$
\enddefinition
\mydefinition{5.3} A frame $(U,\,\boldsymbol\Psi_1,\,\boldsymbol\Psi_2,
\,\boldsymbol\Psi_3,\,\boldsymbol\Psi_4)$ of the Dirac bundle $DM$ is
called a {\it chiral frame\/} if the chirality operator $\bold H$ 
given by the following matrix in this frame:
$$
\hskip -2em
H^{\kern 0.5pti}_{\kern -0.5pt j}=\Vmatrix 1 & 0 & 0 & 0\\0 & 1 & 0 & 0\\
0 & 0 & -1 & 0\\0 & 0 & 0 & -1\endVmatrix.
\mytag{5.4}
$$
\enddefinition
\mydefinition{5.4} A frame $(U,\,\boldsymbol\Psi_1,\,\boldsymbol\Psi_2,
\,\boldsymbol\Psi_3,\,\boldsymbol\Psi_4)$ of the Dirac bundle $DM$ is
called an {\it anti-chiral frame\/} if the chirality operator $\bold H$ 
given by the diagonal matrix opposite to the matrix \mythetag{5.4} in 
this frame:
$$
\hskip -2em
H^{\kern 0.5pti}_{\kern -0.5pt j}=\Vmatrix -1 & 0 & 0 & 0\\0 & -1 & 0 & 0\\
0 & 0 & 1 & 0\\0 & 0 & 0 & 1\endVmatrix.
\mytag{5.5}
$$
\enddefinition
\mydefinition{5.5} A frame $(U,\,\boldsymbol\Psi_1,\,\boldsymbol\Psi_2,
\,\boldsymbol\Psi_3,\,\boldsymbol\Psi_4)$ of the Dirac bundle $DM$ is
called a {\it self-adjoint frame\/} if the Hermitian spin-metric tensor 
$\bold D$ (the Dirac form) is represented by the following matrix in this
frame:
$$
\hskip -2em
D_{i\bar j}=D(\boldsymbol\Psi_{\bar j},\boldsymbol\Psi_i)
=\Vmatrix 0 & 0 & 1 & 0\\0 & 0 & 0 & 1\\
1 & 0 & 0 & 0\\0 & 1 & 0 & 0\endVmatrix.
\mytag{5.6}
$$
\enddefinition
\mydefinition{5.6} A frame $(U,\,\boldsymbol\Psi_1,\,\boldsymbol\Psi_2,
\,\boldsymbol\Psi_3,\,\boldsymbol\Psi_4)$ of the Dirac bundle $DM$ is
called an {\it anti-self-adjoint frame\/} if the Hermitian spin-metric
tensor $\bold D$ (the Dirac form) is represented by the matrix opposite
to the matrix \mythetag{5.6} in this frame:
$$
\hskip -2em
D_{i\bar j}=D(\boldsymbol\Psi_{\bar j},\boldsymbol\Psi_i)
=\Vmatrix 0 & 0 & -1 & 0\\0 & 0 & 0 & -1\\
-1 & 0 & 0 & 0\\0 & -1 & 0 & 0\endVmatrix.
\mytag{5.7}
$$
\enddefinition
In \mycite{6} the $P$ and $T$ reflection operations were studied and the 
following four types of frames in the Dirac bundle $DM$ were considered:
\roster
\item canonically orthonormal chiral frames;
\item $P$-reverse anti-chiral frames;
\item $T$-reverse anti-chiral frames;
\item $PT$-reverse chiral frames.
\endroster\par 
{\bf Canonically orthonormal chiral frames} are simultaneously orthonormal,
chiral, and self-adjoint frames. This is the basic type of frames most
closely related to the expansion \mythetag{5.1}. Each canonically orthonormal 
chiral frame of the Dirac bundle $(U,\,\boldsymbol\Psi_1,\,
\boldsymbol\Psi_2,\,\boldsymbol\Psi_3,\,\boldsymbol\Psi_4)$ is produced
from some orthonormal frame $(U,\,\boldsymbol\Psi_1,\,\boldsymbol\Psi_2)$ 
of the chiral bundle $SM$ as follows:
$$
\xalignat 4
&\boldsymbol\Psi_1=\boldsymbol\Psi^{\kern 0.3pt\sssize \text{chiral}}_1,
&&\boldsymbol\Psi_2=\boldsymbol\Psi^{\kern 0.3pt\sssize \text{chiral}}_2,
&&\boldsymbol\Psi_3=\overline{\boldsymbol\vartheta}
\vphantom{\boldsymbol\vartheta}^{\,1}_{\kern -0.3pt\sssize \text{chiral}},
&&\boldsymbol\Psi_4=\overline{\boldsymbol\vartheta}
\vphantom{\boldsymbol\vartheta}^{\,2}_{\kern -0.3pt\sssize \text{chiral}}.
\qquad
\mytag{5.8}
\endxalignat
$$\par
{\bf $P$-reverse anti-chiral frames} are self-adjoint, but anti-orthonormal
and anti-chiral. Any $P$-reverse anti-chiral frame $(U,\,\tilde{\boldsymbol
\Psi}_1,\,\tilde{\boldsymbol\Psi}_2,\,\tilde{\boldsymbol\Psi}_3,
\,\tilde{\boldsymbol\Psi}_4)$, is produced from some canonically
orthonormal chiral frame $(U,\,\boldsymbol\Psi_1,\,\boldsymbol\Psi_2,\,
\boldsymbol\Psi_3,\,\boldsymbol\Psi_4)$ by $P$-inversion: 
$$
\xalignat 4
&\tilde{\boldsymbol\Psi}_1=\boldsymbol\Psi_3,
&&\tilde{\boldsymbol\Psi}_2=\boldsymbol\Psi_4,
&&\tilde{\boldsymbol\Psi}_3=\boldsymbol\Psi_1,
&&\tilde{\boldsymbol\Psi}_4=\boldsymbol\Psi_2.
\quad
\mytag{5.9}
\endxalignat
$$\par
{\bf $T$-reverse anti-chiral frames} are orthonormal, anti-chiral, 
and anti-self-ad\-joint. Any $T$-reverse anti-chiral frame $(U,\,
\tilde{\boldsymbol\Psi}_1,\,\tilde{\boldsymbol\Psi}_2,\,
\tilde{\boldsymbol\Psi}_3,\,\tilde{\boldsymbol\Psi}_4)$, is produced 
from some canonically orthonormal chiral frame $(U,\,\boldsymbol\Psi_1,
\,\boldsymbol\Psi_2,\,\boldsymbol\Psi_3,\,\boldsymbol\Psi_4)$ by 
$T$-inversion: 
$$
\xalignat 4
&\tilde{\boldsymbol\Psi}_1=i\,\boldsymbol\Psi_3,
&&\tilde{\boldsymbol\Psi}_2=i\,\boldsymbol\Psi_4,
&&\tilde{\boldsymbol\Psi}_3=-i\,\boldsymbol\Psi_1,
&&\tilde{\boldsymbol\Psi}_4=-i\,\boldsymbol\Psi_2.
\qquad\quad
\mytag{5.10}
\endxalignat
$$\par
{\bf $PT$-reverse chiral frames} are anti-orthonormal, chiral, 
and anti-self-ad\-joint. Any $PT$-reverse chiral frame $(U,\,
\tilde{\boldsymbol\Psi}_1,\,\tilde{\boldsymbol\Psi}_2,\,
\tilde{\boldsymbol\Psi}_3,\,\tilde{\boldsymbol\Psi}_4)$, is produced 
from some canonically orthonormal chiral frame $(U,\,\boldsymbol\Psi_1,
\,\boldsymbol\Psi_2,\,\boldsymbol\Psi_3,\,\boldsymbol\Psi_4)$ by 
$PT$-inversion: 
$$
\xalignat 4
&\tilde{\boldsymbol\Psi}_1=i\,\boldsymbol\Psi_1,
&&\tilde{\boldsymbol\Psi}_2=i\,\boldsymbol\Psi_2,
&&\tilde{\boldsymbol\Psi}_3=-i\,\boldsymbol\Psi_3,
&&\tilde{\boldsymbol\Psi}_4=-i\,\boldsymbol\Psi_4.
\qquad\quad
\mytag{5.11}
\endxalignat
$$
All of the above facts are easily derived from \mythetag{5.2},
\mythetag{5.3}, \mythetag{5.4}, \mythetag{5.5}, \mythetag{5.6},
and \mythetag{5.7}. The $P$, $T$, and $PT$-inversions introduced 
in \mythetag{5.9}, \mythetag{5.10}, and \mythetag{5.11} are not 
actual operations over spinors, they are frame transformations
only. The formula \mythetag{5.8} defines a frame construction
operation.\par
    The frames of all of the above four types are canonically
associated with some frames in $TM$. The frame association is
given by the diagram
$$
\hskip -2em
\aligned
&\boxit{Canonically orthonormal}{chiral frames}\to
\boxit{Positively polarized}{right orthonormal frames}\\
&\boxit{$P$-reverse}{anti-chiral frames}\to
\boxit{Positively polarized}{left orthonormal frames}\\
&\boxit{$T$-reverse}{anti-chiral frames}\to
\boxit{Negatively polarized}{right orthonormal frames}\\
&\boxit{$PT$-reverse}{chiral frames}\to
\boxit{Negatively polarized}{left orthonormal frames}
\endaligned
\mytag{5.12}
$$\par
     Like $SM$, the Dirac bundle $DM$ is a complex vector bundle over
the smooth real space-time manifold $M$. For this reason there is a
semilinear involution of complex conjugation $\tau$ acting upon 
spin-tensorial fields associated with $DM$. This involution $\tau$ is
canonically associated with $DM$, it is introduces in a way similar
to $\tau$ for $SM$ (see details in \mycite{6}). Applying $\tau$ to
the spin-metric tensor $\bold d$, we get
$$
\hskip -2em
\bbd=\tau(\bold d).
\mytag{5.13}
$$
This is a spin-tensorial field of the type $(0,0|0,2|0,0)$, while
$\bold d$ itself is a field of the type $(0,2|0,0|0,0)$. The following
formulas are analogous to \mythetag{2.18} and \mythetag{2.19}:
$$
\gather
\hskip -1em
\tau(\bold X)=\msum{4}\Sb i_1,\,\ldots,\,i_\nu\\ j_1,\,\ldots,\,j_\gamma\\
\bar i_1,\,\ldots,\,\bar i_\alpha\\ \bar j_1,\,\ldots,\,
\bar j_\beta\endSb\msum{3}\Sb h_1,\,\ldots,\,h_m\\ k_1,\,\ldots,\,k_n
\endSb
\tau X^{i_1\ldots\,i_\nu\,\bar i_1\ldots\,\bar i_\alpha\,h_1\ldots\,
h_m}_{j_1\ldots\,j_\gamma\,\bar j_1\ldots\,\bar j_\beta\,k_1\ldots\, 
k_n}\ \boldsymbol\Psi^{j_1\ldots\,j_\gamma\,\bar j_1\ldots\,\bar j_\beta
\,k_1\ldots\, k_n}_{i_1\ldots\,i_\nu\,\bar i_1\ldots\,\bar i_\alpha\,
h_1\ldots\,h_m},\qquad
\mytag{5.14}\\
\vspace{-1ex}
\intertext{where}
\vspace{-2ex}
\tau X^{i_1\ldots\,i_\nu\,\bar i_1\ldots\,\bar i_\alpha\,h_1\ldots\,
h_m}_{j_1\ldots\,j_\gamma\,\bar j_1\ldots\,\bar j_\beta\,k_1\ldots\, 
k_n}
=\overline{X^{\bar i_1\ldots\,\bar i_\alpha\,i_1\ldots\,i_\nu\,h_1
\ldots\,h_m}_{\bar j_1\ldots\,\bar j_\beta\,j_1\ldots\,j_\gamma\,k_1
\ldots\,k_n}}.\qquad
\mytag{5.15}
\endgather
$$
Applying \mythetag{5.14} and \mythetag{5.15} to \mythetag{5.13}, we 
derive the components of $\bbd$:
$$
\hskip -2em
\bd_{\bar i\bar j}=\overline{d_{\bar i\bar j}}.
\mytag{5.16}
$$
The local fields $\boldsymbol\Psi^{j_1\ldots\,j_\gamma\,\bar j_1\ldots
\,\bar j_\beta\,k_1\ldots\, k_n}_{i_1\ldots\,i_\nu\,\bar i_1\ldots\,
\bar i_\alpha\,h_1\ldots\,h_m}$ in the expansion \mythetag{5.16} are
introduced in a way similar to that of \mythetag{2.13}. As for the
formulas \mythetag{5.14}, \mythetag{5.15}, and \mythetag{5.16}, they
hold for arbitrary frame pairs $(U,\,\boldsymbol\Psi_1,\,\boldsymbol
\Psi_2,\,\boldsymbol\Psi_3,\,\boldsymbol\Psi_4)$ and $(U,\,\boldsymbol
\Upsilon_0,\,\boldsymbol\Upsilon_1,\,\boldsymbol\Upsilon_2,\,\boldsymbol
\Upsilon_3)$, not only for those listed in the diagram \mythetag{5.12}.
\par
     Let's apply the formula \mythetag{5.15} to the formula \mythetag{5.6}
or to the formula \mythetag{5.7}. As a result we find that the Dirac field 
$\bold D$ is a real spin-tensorial field:
$$
\hskip -2em
\tau(\bold D)=\bold D.
\mytag{5.17}
$$
From \mythetag{5.17} for the Dirac form $D(\bold X,\bold Y)=C(\bold D
\otimes\tau(\bold X)\otimes\bold Y)$ we derive
$$
\hskip -2em
D(\bold X,\bold Y)=\sum^4_{i=1}\sum^4_{j=1}D_{i\bar j}\,
\overline{X^{\bar j}}\,Y^i=\overline{D(\bold Y,\bold X)}.
\mytag{5.18}
$$
The identity \mythetag{5.18} shows that the Dirac form $D(\bold X,
\bold Y)$ is a Hermitian form.\par
\head
6. Dirac's $\gamma$-field, $\gamma$-symbols, and $\gamma$-matrices.
\endhead
    The Infeld-van der Waerden symbols are not expanded to the 
Dirac bundle. Instead of them, here we have the Dirac's $\gamma$-field
$\boldsymbol\gamma$. This is a spin-tensorial field of the type
$(1,1|0,0|0,1)$. Its components are called the {\it Dirac's 
$\gamma$-symbols}. In a frame pair of any type listed on the diagram 
\mythetag{5.12}, $\gamma$-symbols are given explicitly:
$$
\xalignat 4
&\hskip -2em
\gamma^{\,1}_{1\,0}=0,
&&\gamma^{\,1}_{2\,0}=0,
&&\gamma^{\,1}_{3\,0}=1,
&&\gamma^{\,1}_{4\,0}=0,
\quad\\
&\hskip -2em
\gamma^{\,2}_{1\,0}=0,
&&\gamma^{\,2}_{2\,0}=0,
&&\gamma^{\,2}_{3\,0}=0,
&&\gamma^{\,2}_{4\,0}=1,
\quad\\
\vspace{-1.5ex}
&&&&&&&\mytag{6.1}\\
\vspace{-1.5ex}
&\hskip -2em
\gamma^{\,3}_{1\,0}=1,
&&\gamma^{\,3}_{2\,0}=0,
&&\gamma^{\,3}_{3\,0}=0,
&&\gamma^{\,3}_{4\,0}=0,
\quad\\
&\hskip -2em
\gamma^{\,4}_{1\,0}=0,
&&\gamma^{\,4}_{2\,0}=1,
&&\gamma^{\,4}_{3\,0}=0,
&&\gamma^{\,4}_{4\,0}=0,
\quad\\
\displaybreak
\vspace{2ex}
&\hskip -2em
\gamma^{\,1}_{1\,1}=0,
&&\gamma^{\,1}_{2\,1}=0,
&&\gamma^{\,1}_{3\,1}=0,
&&\gamma^{\,1}_{4\,1}=1,
\quad\\
&\hskip -2em
\gamma^{\,2}_{1\,1}=0,
&&\gamma^{\,2}_{2\,1}=0,
&&\gamma^{\,2}_{3\,1}=1,
&&\gamma^{\,2}_{4\,1}=0,
\quad\\
\vspace{-1.5ex}
&&&&&&&\mytag{6.2}\\
\vspace{-1.5ex}
&\hskip -2em
\gamma^{\,3}_{1\,1}=0,
&&\gamma^{\,3}_{2\,1}=-1,
&&\gamma^{\,3}_{3\,1}=0,
&&\gamma^{\,3}_{4\,1}=0,
\quad\\
&\hskip -2em
\gamma^{\,4}_{1\,1}=-1,
&&\gamma^{\,4}_{2\,1}=0,
&&\gamma^{\,4}_{3\,1}=0,
&&\gamma^{\,4}_{4\,1}=0,
\quad\\
&\hskip -2em
\gamma^{\,1}_{1\,2}=0,
&&\gamma^{\,1}_{2\,2}=0,
&&\gamma^{\,1}_{3\,2}=0,
&&\gamma^{\,1}_{4\,2}=-i,
\quad\\
&\hskip -2em
\gamma^{\,2}_{1\,2}=0,
&&\gamma^{\,2}_{2\,2}=0,
&&\gamma^{\,2}_{3\,2}=i,
&&\gamma^{\,2}_{4\,2}=0,
\quad\\
\vspace{-1.5ex}
&&&&&&&\mytag{6.3}\\
\vspace{-1.5ex}
&\hskip -2em
\gamma^{\,3}_{1\,2}=0,
&&\gamma^{\,3}_{2\,2}=i,
&&\gamma^{\,3}_{3\,2}=0,
&&\gamma^{\,3}_{4\,2}=0,
\quad\\
&\hskip -2em
\gamma^{\,4}_{1\,2}=-i,
&&\gamma^{\,4}_{2\,2}=0,
&&\gamma^{\,4}_{3\,2}=0,
&&\gamma^{\,4}_{4\,2}=0,
\quad\\
\vspace{2ex}
&\hskip -2em
\gamma^{\,1}_{1\,3}=0,
&&\gamma^{\,1}_{2\,3}=0,
&&\gamma^{\,1}_{3\,3}=1,
&&\gamma^{\,1}_{4\,3}=0,
\quad\\
&\hskip -2em
\gamma^{\,2}_{1\,3}=0,
&&\gamma^{\,2}_{2\,3}=0,
&&\gamma^{\,2}_{3\,3}=0,
&&\gamma^{\,2}_{4\,3}=-1,
\quad\\
\vspace{-1.5ex}
&&&&&&&\mytag{6.4}\\
\vspace{-1.5ex}
&\hskip -2em
\gamma^{\,3}_{1\,3}=-1,
&&\gamma^{\,3}_{2\,3}=0,
&&\gamma^{\,3}_{3\,3}=0,
&&\gamma^{\,3}_{4\,3}=0,
\quad\\
&\hskip -2em
\gamma^{\,4}_{1\,3}=0,
&&\gamma^{\,4}_{2\,3}=1,
&&\gamma^{\,4}_{3\,3}=0,
&&\gamma^{\,4}_{4\,3}=0.
\quad
\endxalignat
$$
The second lower index of the $\gamma$-symbols \mythetag{6.1},
\mythetag{6.2}, \mythetag{6.3}, and \mythetag{6.4} is a spacial
index. By fixing this index, we can arrange $\gamma$-symbols 
into four square matrices
$$
\hskip -2em
\gamma_k=\Vmatrix \gamma^{\,1}_{1\kern 0.2pt k} 
&\gamma^{\,1}_{2\kern 0.2pt k} &\gamma^{\,1}_{3\kern 0.2pt k}
&\gamma^{\,1}_{4\kern 0.2pt k}\\ 
\vspace{1.5ex}
\gamma^{\,2}_{1\kern 0.2pt k} &\gamma^{\,2}_{2\kern 0.2pt k} 
&\gamma^{\,2}_{3\kern 0.2pt k} &\gamma^{\,2}_{4\kern 0.2pt k}\\
\vspace{1.5ex}
\gamma^{\,3}_{1\kern 0.2pt k} &\gamma^{\,3}_{2\kern 0.2pt k} 
&\gamma^{\,3}_{3\kern 0.2pt k} &\gamma^{\,3}_{4\kern 0.2pt k}\\ 
\vspace{1.5ex}
\gamma^{\,4}_{1\kern 0.2pt k} &\gamma^{\,4}_{2\kern 0.2pt k} 
&\gamma^{\,4}_{3\kern 0.2pt k} &\gamma^{\,4}_{4\kern 0.2pt k}
\endVmatrix,
\qquad k=0,\,1,\,2,\,3.
\mytag{6.5}
$$
The matrices \mythetag{6.5} are called {\it Dirac matrices}. One can
write them explicitly:
$$
\xalignat 2
&\hskip -2em
\gamma_0=\Vmatrix 0&0&1&0\\0&0&0&1\\1&0&0&0\\0&1&0&0\endVmatrix,
&&\gamma_1=\Vmatrix 0&0&0&1\\0&0&1&0\\0&-1&0&0\\-1&0&0&0\endVmatrix,
\quad\\
\vspace{-1.5ex}
&&&\mytag{6.6}\\
\vspace{-1.5ex}
&\hskip -2em
\gamma_2=\Vmatrix 0&0&0&-i\\0&0&i&0\\0&i&0&0\\-i&0&0&0\endVmatrix,
&&\gamma_3=\Vmatrix 0&0&1&0\\0&0&0&-1\\-1&0&0&0\\0&1&0&0\endVmatrix.
\quad
\endxalignat
$$
The Dirac matrices \mythetag{6.6} are very popular in physics. However,
dealing with them, one should remember that each separate matrix has
no spin-tensorial interpretation.\par
    The most popular property of the Dirac matrices \mythetag{6.6} is
written in terms of their anticommutators $\{\gamma_i,\,\gamma_j\}=
\gamma_i\cdot\gamma_j-\gamma_j\cdot\gamma_i$:
$$
\hskip -2em
\{\gamma_i,\,\gamma_j\}=2\,g_{ij}\ \bold 1.
\mytag{6.7}
$$
Here $\bold 1$ is the unit matrix. \pagebreak Due to this property they
define the $4$-dimensional representation of the Clifford algebra
$\Cal Cl(1,3,\Bbb R)$ (see \mycite{8}). In terms of the $\gamma$-symbols
the formula \mythetag{6.7} is written as follows:
$$
\hskip -2em
\sum^4_{b=1}\gamma^{\,a}_{b\kern 0.5pt i}\,\gamma^{\,b}_{cj}
+\sum^4_{b=1}\gamma^{\,a}_{bj}\,\gamma^{\,b}_{c\kern 0.5pt i}
=2\,g_{ij}\,\delta^a_c.
\mytag{6.8}
$$
Apart from \mythetag{6.8}, there are also some analogs of the properties
of Infeld-van der Waerden symbols \mythetag{3.8}, \mythetag{3.9},
\mythetag{3.10}, and \mythetag{3.11}. Here is the most simple of them:
$$
\hskip -2em
\sum^4_{a=1}\sum^4_{b=1}\sum^4_{e=1}\sum^4_{h=1}
\gamma^{\,a}_{b\kern 0.5pt i}\,d_{ae}\,d^{\kern 0.5pt bh}
\,\gamma^{\,e}_{h\kern 0.2pt j}
=4\,g_{ij}.
\mytag{6.9}
$$
The identity \mythetag{6.9} is an analog of \mythetag{3.9}. It is 
derived from the following more simple identity with the use of 
the anticommutator relationship \mythetag{6.8}:
$$
\hskip -2em
\sum^4_{a=1}\sum^4_{b=1}
\gamma^{\,a}_{b\kern 0.5pt i}\,d_{ae}\,d^{\kern 0.5pt bh}
=\gamma^{\,h}_{e\kern 0.5pt i}.
\mytag{6.10}
$$
By $d^{\kern 0.5pt bh}$ in \mythetag{6.9} and \mythetag{6.10} we denote
the components of the dual spin-metric tensor. By tradition we denote it
by the same symbol $\bold d$. Its components $d^{\kern 0.5pt bh}$ form 
the matrix inverse to $d_{ae}$.\par
    The inverse Dirac's $\gamma$-field $\boldsymbol\gamma$ is a 
spin-tensorial field of the type $(1,1|0,0|1,0)$. Its components are
obtained from $\gamma^{\,i}_{jm}$ by raising the lower index $m$:
$$
\hskip -2em
\gamma^{\,i\kern 0.5pt m}_{j}=\sum^3_{k=0}\gamma^{\,i}_{jk}\ g^{km}.
\mytag{6.11}
$$
The formula \mythetag{6.11} is an analog of the formula \mythetag{3.2}.
The quantities $\gamma^{\,i\kern 0.5pt m}_{j}$ obtained through this 
formula are called the {\it inverse $\gamma$-symbols}. From \mythetag{6.9} 
and \mythetag{6.10}, taking into account \mythetag{6.11}, now we derive  
$$
\hskip -2em
\sum^4_{e=1}\sum^4_{h=1}\gamma^{\,h}_{e\kern 0.2pt j}
\ \gamma^{\,e\kern 0.2pt i}_{h}=4\,\delta^i_j.
\mytag{6.12}
$$
The formula \mythetag{6.12} is an analog of the first formula
\mythetag{3.12}. By raising the index $i$ in \mythetag{6.10} we 
obtain the following equality for the inverse $\gamma$-symbols:
$$
\hskip -2em
\sum^4_{a=1}\sum^4_{b=1}
\gamma^{\,a\kern 0.5pt i}_{b}\,d_{ae}\,d^{\kern 0.5pt bh}
=\gamma^{\,h\kern 0.5pt i}_{e}.
\mytag{6.13}
$$
Then rising both indices $i$ and $j$ in \mythetag{6.9} we derive
the identity
$$
\hskip -2em
\sum^4_{a=1}\sum^4_{b=1}\sum^4_{e=1}\sum^4_{h=1}
\gamma^{\,a\kern 0.5pt i}_{b}\,d_{ae}\,d^{\kern 0.5pt bh}
\,\gamma^{\,e\kern 0.2pt j}_{h}
=4\,g^{ij}.
\mytag{6.14}
$$
The identity \mythetag{6.14} is analogous to \mythetag{3.11}. As for
the identities \mythetag{6.10} and \mythetag{6.13}, they have no analogs 
in chiral spinors.\par
    The relation of Dirac's $\gamma$-field and the chirality operator
$\bold H$ is determined by the structure of $\gamma$-matrices 
\mythetag{6.6}. By means of direct calculations we prove that
$$
\hskip -2em
\{\gamma_m,\,\bold H\}=0.
\mytag{6.15}
$$
Like in \mythetag{6.7}, in \mythetag{6.15} we have the anticommutator
$\{\gamma_m,\,\bold H\}=\gamma_m\cdot\bold H+\bold H\cdot\gamma_m$.
In a coordinate form the identity \mythetag{6.15} is written as
$$
\hskip -2em
\sum^4_{b=1}\gamma^{\,a}_{b\kern 0.5pt m}\,H^b_c
+\sum^4_{b=1}H^a_b\,\gamma^{\,b}_{c\kern 0.5pt m}=0.
\mytag{6.16}
$$
As for $\bold d$ and $\bold H$, their relation is described by the
identity similar to \mythetag{6.16}:
$$
\hskip -2em
\sum^4_{b=1}d_{ab}\,H^b_c=\sum^4_{b=1}H^b_a\,d_{bc}.
\mytag{6.17}
$$
The identity \mythetag{6.17} means that the chirality operator 
$\bold H$ is a symmetric operator with respect to the bilinear 
form of the spin-metric tensor $\bold d$, i\.\,e\.
$$
\hskip -2em
d(\bold H(\bold X),\bold Y)=d(\bold X,\bold H(\bold Y))
\mytag{6.18}
$$
for any two spinors $\bold X$ and $\bold Y$. In the case if the Dirac
form (the form of the Hermitian spin-metric tensor $\bold D$) we have
the identity similar to \mythetag{6.18}:
$$
\hskip -2em
D(\bold H(\bold X),\bold Y)=-D(\bold X,\bold H(\bold Y)).
\mytag{6.19}
$$
The identity \mythetag{6.19} means that $\bold H$ is an anti-Hermitian
operator with respect to the Hermitian form $D$. In a coordinate form
\mythetag{6.19} is written as
$$
\hskip -2em
\sum^4_{\bar a=1}D_{i\bar a}\,\overline{H^{\bar a}_{\bar i}}=
-\sum^4_{a=1}H^a_i\,D_{a\bar i}.
\mytag{6.20}
$$
Returning back to the Dirac's $\gamma$-symbols, we can write
$$
\align
&\hskip -2em
\gathered
\sum^3_{m=0}\sum^3_{n=0}\gamma^{\,a}_{b\kern 0.5pt m}\,
\gamma^{\,e}_{h\kern 0.5pt n}\ g^{mn}=\delta^a_h\ \delta^e_b
-H^a_h\ H^e_b\,+\\
+\,d^{\kern 0.5pt ae}\,d_{\kern 0.5pt bh}-\sum^4_{r=1}\sum^4_{s=1}
H^a_r\,d^{\kern 0.5pt re}\,d_{\kern 0.5pt bs}\,H^s_h,
\endgathered
\mytag{6.21}\\
\vspace{2ex}
&\hskip -2em
\gathered
\sum^3_{m=0}\sum^3_{n=0}\gamma^{\,a\kern 0.5pt m}_b\,
\gamma^{\,e\kern 0.5pt n}_h\ g_{mn}=\delta^a_h\ \delta^e_b
-H^a_h\ H^e_b\,+\\
+\,d^{\kern 0.5pt ae}\,d_{\kern 0.5pt bh}-\sum^4_{r=1}\sum^4_{s=1}
H^a_r\,d^{\kern 0.5pt re}\,d_{\kern 0.5pt bs}\,H^s_h.
\endgathered
\mytag{6.22}
\endalign
$$
The identities \mythetag{6.21} and \mythetag{6.22} are analogs of 
\mythetag{3.8} and \mythetag{3.10}. \pagebreak Taking into account
\mythetag{6.11}, we can transform \mythetag{6.21} and \mythetag{6.22} 
to the following identity:
$$
\hskip -2em
\gathered
\sum^3_{m=0}\gamma^{\,a\kern 0.5pt m}_b\,
\gamma^{\,e}_{h\kern 0.5pt m}=\delta^a_h\ \delta^e_b
-H^a_h\ H^e_b\,+\\
+\,d^{\kern 0.5pt ae}\,d_{\kern 0.5pt bh}-\sum^4_{r=1}\sum^4_{s=1}
H^a_r\,d^{\kern 0.5pt re}\,d_{\kern 0.5pt bs}\,H^s_h.
\endgathered
\mytag{6.23}
$$
The identity \mythetag{6.23} is an analog of the second identity
\mythetag{3.12}.\par
     Let's take the components of the Dirac form (the Hermitian 
spin-metric tensor $\bold D$) and raise their indices. As a result 
we obtain
$$
\hskip -2em
D^{i\kern 0.5pt\bar i}=\sum^4_{a=1}\sum^4_{\bar a=1}
d^{\kern 0.5pt i\kern 0.5pt a}\,D_{a\bar a}
\ \bd^{\kern 0.5pt\bar a\bar i}.
\mytag{6.24}
$$
The matrix $D^{i\kern 0.5pt\bar i}$ is inverse to $D_{i\kern 0.5pt\bar i}$
in the sense of the following equalities:
$$
\xalignat 2
&\hskip -2em
\sum^4_{\bar a=1}D_{j\kern 0.5pt\bar a}\,D^{i\kern 0.5pt\bar a}
=\delta^{\kern 0.5pt i}_j,
&&\sum^4_{\bar a=1}D_{\kern -0.5pt a\bar j}\,D^{a\kern 0.5pt\bar i}
=\delta^{\kern 0.5pt\bar i}_{\bar j}.
\mytag{6.25}
\endxalignat
$$
Due to \mythetag{6.25} the spin-tensorial field $\bold D$ of the type
$(1,0|1,0|0,0)$ determined by the matrix \mythetag{6.24} is called
the {\it inverse Hermitian spin-metric tensor}. Using both
$D^{i\kern 0.5pt\bar i}$ and $D_{i\kern 0.5pt\bar i}$, we define the
following quantities:
$$
\xalignat 2
&\hskip -2em
\gamma^{\kern 0.5pt i\kern 0.5pt\bar i}_m=\sum^4_{a=1}
\gamma^{\kern 0.5pt i}_{a\kern 0.5pt m}\,D^{\kern 0.5pt a\bar i},
&&\gamma^m_{i\kern 0.5pt\bar i}=\sum^4_{a=1}
\gamma^{\kern 0.5pt a\kern 0.5pt m}_{i}\,D_{\kern -0.5pt a\bar i}.
\mytag{6.26}
\endxalignat
$$
The quantities \mythetag{6.26} are called {\it direct and inverse
Hermitian $\gamma$-symbols}. They define two spin-tensorial fields
of the types $(1,0|1,0|0,1)$ and $(0,1|0,1|1,0)$ respectively. We
denote these fields by the same symbol $\boldsymbol\gamma$, as well 
as the initial fields from which they are produced. The fields
\mythetag{6.26} are real fields:
$$
\hskip -2em
\tau(\boldsymbol\gamma)=\boldsymbol\gamma.
\mytag{6.27}
$$
Indeed, by means of direct calculations we can prove that
$$
\xalignat 2
&\hskip -2em
\gamma^{\kern 0.5pt i\kern 0.5pt\bar i}_m
=\overline{\gamma^{\kern 0.5pt\bar i\kern 0.5pt i}_m},
&&\gamma^m_{i\kern 0.5pt\bar i}
=\overline{\gamma^m_{\bar i\kern 0.5pt i}}.
\mytag{6.28}
\endxalignat
$$
The equalities \mythetag{6.28} are coordinate representations of the
equality \mythetag{6.27}. In terms of the initial $\gamma$-symbols
they can be written as
$$
\hskip -2em
\sum^4_{\bar a=1}D_{\kern -0.5pt i\kern 0.5pt\bar a}\,
\overline{\gamma^{\kern 0.5pt\bar a}_{\kern 0.5pt\bar i\kern 0.5pt m}}
=\sum^4_{a=1}\gamma^{\kern 0.5pt a}_{\kern 0.5pt i\kern 0.5pt m}
\,D_{\kern -0.5pt a\bar i}.
\mytag{6.29}
$$
The equality \mythetag{6.29} is similar to \mythetag{6.20}. Moreover, we
have
$$
\hskip -2em
\sum^4_{a=1}d_{i\kern 0.5pt a}\,
\gamma^{\kern 0.5pt a}_{j\kern 0.5pt m}
=-\sum^4_{a=1}\gamma^{\kern 0.5pt a}_{i\kern 0.5pt m}
\,d_{\kern 0.5pt aj}.
\mytag{6.30}
$$
The equality \mythetag{6.30} is derived from \mythetag{6.10}, it is
similar to \mythetag{6.20} and \mythetag{6.29}.\par
\head
7. Spinor connections for Dirac spinors.
\endhead
     Let $(U,\,\boldsymbol\Upsilon_0,\,\boldsymbol\Upsilon_1,\,
\boldsymbol\Upsilon_2,\,\boldsymbol\Upsilon_3)$ and $(U,\,\boldsymbol
\Psi_1,\,\boldsymbol\Psi_2,\,\boldsymbol\Psi_3,\,\boldsymbol\Psi_4)$
be two frames with a common domain $U$ of the bundles $TM$ and $DM$ 
respectively. Let $(\tilde U,\,\tilde{\boldsymbol\Upsilon}_0,\,
\tilde{\boldsymbol\Upsilon}_1,\,\tilde{\boldsymbol\Upsilon}_2,\,
\tilde{\boldsymbol\Upsilon}_3)$ and $(\tilde U,\,\tilde{\boldsymbol
\Psi}_1,\,\tilde{\boldsymbol\Psi}_2,\,\tilde{\boldsymbol\Psi}_3,
\,\tilde{\boldsymbol\Psi}_4)$ be other two such frames. Assume that 
$U\cap\tilde U\neq\varnothing$. Then at each point $p\in U\cap\tilde
U\neq\varnothing$ one can write the transition formulas \mythetag{1.4}.
Instead of \mythetag{1.6} here we write the following transition formulas:
$$
\xalignat 2
&\hskip -2em
\tilde{\boldsymbol\Psi}_i=\sum^4_{j=1}\goth S^j_i\,
\boldsymbol\Psi_j,
&&\boldsymbol\Psi_i=\sum^4_{j=1}\goth T^j_i\,
\tilde{\boldsymbol\Psi}_j.
\mytag{7.1}
\endxalignat
$$
Using the transition matrices from \mythetag{1.4} and \mythetag{7.1},
one can define the $\theta$-parameters. They are introduced by the
formulas which are almost the same as the formulas \mythetag{4.6}, 
\mythetag{4.7}, \mythetag{4.8}, and \mythetag{4.9} in section~4:
$$
\align
&\hskip -2em
\tilde\theta^k_{ij}=\sum^3_{a=0}T^k_a\,L_{\tilde{\boldsymbol\Upsilon}_i}
\!(S^a_j)=-\sum^3_{a=0}L_{\tilde{\boldsymbol\Upsilon}_i}
\!(T^k_a)\,S^a_j,
\mytag{7.2}\\
&\hskip -2em
\tilde\vartheta^k_{ij}=\sum^4_{a=1}\goth T^k_a\,
L_{\tilde{\boldsymbol\Upsilon}_i}(\goth S^a_j)
=-\sum^4_{a=1}L_{\tilde{\boldsymbol\Upsilon}_i}
\!(\goth T^k_a)\,\goth S^a_j,
\mytag{7.3}\\
&\hskip -2em
\theta^k_{ij}=\sum^3_{a=0}S^k_a\,L_{\boldsymbol\Upsilon_i}
\!(T^a_j)=-\sum^3_{a=0}L_{\boldsymbol\Upsilon_i}
\!(S^k_a)\,T^a_j,
\mytag{7.4}\\
&\hskip -2em
\vartheta^k_{ij}=\sum^4_{a=1}\goth S^k_a\,
L_{\boldsymbol\Upsilon_i}(\goth T^a_j)
=-\sum^4_{a=1}L_{\boldsymbol\Upsilon_i}
\!(\goth S^k_a)\,\goth T^a_j.
\mytag{7.5}
\endalign
$$
The only difference is that the indices $i$ and $k$ in \mythetag{7.3}
and \mythetag{7.5} run over the range from $1$ to $4$. The formulas
\mythetag{7.2} and \mythetag{7.4} coincide with \mythetag{4.6} and 
\mythetag{4.8} exactly. For this reason $\theta$-parameters
$\tilde\theta^k_{ij}$ and $\theta^k_{ij}$ here coincide with those in
section~4 and the relationships \mythetag{4.10}, \mythetag{4.11},
and \mythetag{4.14} for them are valid.
\mydefinition{7.1} A {\it spinor connection} in the Dirac bundle $DM$ 
is a geometric object such that in each frame pair $(U,\,\boldsymbol
\Upsilon_0,\,\boldsymbol\Upsilon_1,\,\boldsymbol\Upsilon_2,\,
\boldsymbol\Upsilon_3)$ and $(U,\,\boldsymbol\Psi_1,\,\boldsymbol
\Psi_2,\,\boldsymbol\Psi_3,\,\boldsymbol\Psi_4)$ of the bundles 
$TM$ and $DM$ it is given by three arrays of smooth complex-valued 
functions
$$
\align
\Gamma^k_{ij}&=\Gamma^k_{ij}(p),\quad i,j,k=0,\,\ldots,\,3,\\
\vspace{1ex}
\Alpha^k_{ij}&=\Alpha^k_{ij}(p),\quad i=0,\,\ldots,\,3,\quad 
j,k=1,\,\ldots,\,4,\\
\vspace{1ex}
\bar{\Alpha}\vphantom{\Alpha}^k_{ij}
&=\bar{\Alpha}\vphantom{\Alpha}^k_{ij}(p),\quad i=0,\,\ldots,\,3,
\quad j,k=1,\,\ldots,\,4,
\endalign
$$
where $p\in U$, such that when passing from the frame pair  $(U,\,
\boldsymbol\Upsilon_0,\,\boldsymbol\Upsilon_1,\,\boldsymbol
\Upsilon_2,\,\boldsymbol\Upsilon_3)$\linebreak and $(U,\,\boldsymbol
\Psi_1,\,\boldsymbol\Psi_2,\,\boldsymbol\Psi_3,\,\boldsymbol\Psi_4)$ 
to some other frame pair $(\tilde U,\,\tilde{\boldsymbol\Upsilon}_0,
\,\tilde{\boldsymbol\Upsilon}_1,\,\tilde{\boldsymbol\Upsilon}_2,\,
\tilde{\boldsymbol\Upsilon}_3)$ and\linebreak $(\tilde U,\,
\tilde{\boldsymbol\Psi}_1,\,\tilde{\boldsymbol\Psi}_2,\tilde{\boldsymbol
\Psi}_3,\,\tilde{\boldsymbol\Psi}_4)$ with $U\cap\tilde U\neq\varnothing$
these functions are transformed as follows:
$$
\allowdisplaybreaks
\align
&\hskip -2em
\Gamma^k_{ij}=\dsize\sum^3_{b=0}\sum^3_{a=0}\sum^3_{c=0}
S^k_a\,T^b_j\,T^c_i\ \tilde\Gamma^a_{c\,b}+\theta^k_{ij},
\mytag{7.6}\\
&\hskip -2em
\Alpha^k_{ij}=\dsize\sum^4_{b=1}\sum^4_{a=1}\sum^3_{c=0}
\goth S^k_a\,\goth T^b_j\,T^c_i\ \tilde{\Alpha}\vphantom{\Alpha}^a_{c\,b}
+\vartheta^k_{ij},
\mytag{7.7}\\
&\hskip -2em
\bar{\Alpha}\vphantom{\Alpha}^k_{ij}=\sum^4_{b=1}\sum^4_{a=1}
\sum^3_{c=0}\overline{\goth S^k_a}\ \overline{\goth T^b_j}\,T^c_i\  
\tilde{\bar{\Alpha}}\vphantom{\Alpha}^a_{c\,b}
+\overline{\vartheta^k_{ij}}.
\mytag{7.8}
\endalign
$$
\enddefinition
The components of the transition matrices $S$, $T$, $\goth S$, and $\goth T$
in \mythetag{7.6}, \mythetag{7.7}, and \mythetag{7.8} are taken from
\mythetag{1.4} and \mythetag{7.1}, while the quantities $\theta^k_{ij}$
and $\vartheta^k_{ij}$ are defined in \mythetag{7.4} and \mythetag{7.5}.
The covariant differential $\nabla$ associated with the spinor
connection $(\Gamma,\Alpha,\bar{\Alpha})$ introduced in the
definition~\mythedefinition{7.1} is a differential operator
$$
\hskip -2em
\nabla\!:\,D^\alpha_\beta\bar D^\nu_\gamma T^m_n M
\to D^\alpha_\beta\bar D^\nu_\gamma T^m_{n+1} M.
\mytag{7.9}
$$
In a coordinate form the operator \mythetag{7.9} is represented by a
covariant derivative:
$$
\hskip -5em
\gathered
\nabla_{\!k_{n+1}}X^{i_1\ldots\,i_\alpha\,\bar i_1\ldots\,\bar i_\nu
\,h_1\ldots\,h_m}_{j_1\ldots\,j_\beta\,\bar j_1\ldots\,\bar j_\gamma
\,k_1\ldots\, k_n}
=L_{\boldsymbol\Upsilon_{k_{n+1}}}\!\bigl(X^{i_1\ldots\,i_\alpha\,\bar i_1
\ldots\,\bar i_\nu\,h_1\ldots\,h_m}_{j_1\ldots\,j_\beta\,\bar j_1
\ldots\,\bar j_\gamma\,k_1\ldots\, k_n}\bigr)\,-\\
\vspace{2ex}
\gathered
\kern -9em
+\sum^\alpha_{\mu=1}\sum^4_{v_\mu=1}\Alpha^{i_\mu}_{k_{n+1}\,v_\mu}\ 
X^{i_1\ldots\,v_\mu\,\ldots\,i_\varepsilon\,\bar i_1\ldots\,\bar i_\nu
\,h_1\ldots\,h_m}_{j_1\ldots\,\ldots\,\ldots\,j_\beta\,\bar j_1\ldots\,
\bar j_\gamma\,k_1\ldots\,k_n}\,-\\
\kern 9em-\sum^\beta_{\mu=1}\sum^4_{w_\mu=1}\Alpha^{w_\mu}_{k_{n+1}\,j_\mu}
\ X^{i_1\ldots\,\ldots\,\ldots\,i_\alpha\,\bar i_1\ldots\,\bar i_\nu
h_1\ldots\,h_m}_{j_1\ldots\,w_\mu\,\ldots\,j_\beta\,\bar j_1\ldots\,
\bar j_\gamma k_1\ldots\,k_n}\,+\\
\kern -9em
+\sum^\nu_{\mu=1}\sum^4_{v_\mu=1}
\bar{\Alpha}\vphantom{\Alpha}^{\bar i_\mu}_{k_{n+1}\,v_\mu}\ 
X^{i_1\ldots\,i_\alpha\,\bar i_1\ldots\,v_\mu\,\ldots\,\bar i_\nu
\,h_1\ldots\,h_m}_{j_1\ldots\,j_\beta\,\bar j_1\ldots\,\ldots\,\ldots\,
\bar j_\gamma\,k_1\ldots\,k_n}\,-\\
\kern 9em-\sum^\gamma_{\mu=1}\sum^4_{w_\mu=1}
\bar{\Alpha}\vphantom{\Alpha}^{w_\mu}_{k_{n+1}\,\bar j_\mu}\
X^{i_1\ldots\,i_\alpha\,\bar i_1\ldots\,\ldots\,\ldots\,\bar i_\nu
\,h_1\ldots\,h_m}_{j_1\ldots\,j_\beta\,\bar j_1\ldots\,w_\mu\,\ldots\,
\bar j_\gamma\,k_1\ldots\,k_n}\,+\\
\kern -9em+\sum^m_{\mu=1}\sum^3_{v_\mu=0}\Gamma^{h_\mu}_{k_{n+1}\,v_\mu}\ 
X^{i_1\ldots\,i_\alpha\,\bar i_1\ldots\,\bar i_\nu\,
h_1\ldots\,v_\mu\,\ldots\,h_m}_{j_1\ldots\,j_\beta\,\bar j_1\ldots\,
\bar j_\gamma\,k_1\ldots\,\ldots\,\ldots\,k_n}\,-\\
\kern 9em-\sum^n_{\mu=1}\sum^3_{w_\mu=0}\Gamma^{w_\mu}_{k_{n+1}\,k_\mu}\
X^{i_1\ldots\,i_\alpha\,\bar i_1\ldots\,\bar i_\nu\,
h_1\ldots\,\ldots\,\ldots\,h_m}_{j_1\ldots\,j_\beta\,\bar j_1\ldots\,
\bar j_\gamma\,k_1\ldots\,w_\mu\,\ldots\,k_n}.
\endgathered\kern 4em
\endgathered\hskip -4em
\mytag{7.10}
$$
The formula \mythetag{7.9} is analogous to \mythetag{4.19}, while 
\mythetag{7.10} is an analog of \mythetag{4.20}.\par
\mydefinition{7.2} A spinor connection $(\Gamma,\Alpha,\bar{\Alpha})$
of the bundle of Dirac spinors $DM$ is called {\it concordant with 
the complex conjugation\/} or a {\it real connection} if the corresponding 
covariant differential \mythetag{7.9} commute with the involution $\tau$,
i\.\,e\. if $\nabla(\tau(\bold X))=\tau(\nabla\bold X)$ for any 
spin-tensorial field $\bold X$.
\enddefinition
\mydefinition{7.3} A spinor connection $(\Gamma,\Alpha,\bar{\Alpha})$
of the bundle of Dirac spinors $DM$ is called {\it concordant with the
Dirac's $\gamma$-field\/} if $\nabla\boldsymbol\gamma=0$.
\enddefinition
\mydefinition{7.4} A spinor connection $(\Gamma,\Alpha,\bar{\Alpha})$
of the bundle of Dirac spinors $DM$ is called {\it concordant with the
spin-metric tensor\/} if $\nabla\bold d=0$.
\enddefinition
\mydefinition{7.5} A spinor connection $(\Gamma,\Alpha,\bar{\Alpha})$
of the bundle of Dirac spinors $DM$ is called {\it concordant with the
metric tensor\/} if $\nabla\bold g=0$.
\enddefinition
\mydefinition{7.6} A spinor connection $(\Gamma,\Alpha,\bar{\Alpha})$
of the bundle of Dirac spinors $DM$ is called {\it concordant with the
chirality operator\/} if $\nabla\bold H=0$.
\enddefinition
\mydefinition{7.7} A spinor connection $(\Gamma,\Alpha,\bar{\Alpha})$
of the bundle of Dirac spinors $DM$ is called {\it concordant with the
Hermitian spin-metric tensor\/} if $\nabla\bold D=0$.
\enddefinition
\mytheorem{7.1} A spinor connection $(\Gamma,\Alpha,\bar{\Alpha})$
of the bundle of Dirac spinors $DM$ is concordant with the complex 
conjugation $\tau$ if and only if 
$$
\xalignat 2
&\hskip -2em
\Gamma^k_{ij}=\overline{\Gamma^k_{ij}},
&&\bar{\Alpha}\vphantom{\Alpha}^k_{ij}=\overline{\Alpha^k_{ij}}.
\mytag{7.11}
\endxalignat
$$
\endproclaim
     The theorem~\mythetheorem{7.1} is an analog of the 
theorem~\mythetheorem{4.1}. Its proof is obvious due to the formulas
\mythetag{5.14}, \mythetag{5.15}, and \mythetag{7.10}.\par
     The Dirac $\gamma$-symbols are more numerous than the Infeld-van 
der Waerden symbols. For this reason they contain more information and
they are more self-sufficient. Instead of the theorem~\mythetheorem{4.2}
here we have.
\mytheorem{7.2} Any spinor connection $(\Gamma,\Alpha,\bar{\Alpha})$
of the bundle of Dirac spinors $DM$ concordant with the Dirac 
$\gamma$-field $\boldsymbol\gamma$ is concordant with the metric tensor
$\bold g$ as well, i\.\,e\. $\nabla\boldsymbol\gamma=0$ implies $\nabla
\bold g=0$.
\endproclaim
\demo{Proof} In order to prove this theorem it is sufficient to apply
the identity \mythetag{6.8}. By setting $c=a$ and summing over the index
$a$ from \mythetag{6.8} we derive
$$
\hskip -2em
\sum^4_{a=1}\sum^4_{b=1}\gamma^{\,a}_{b\kern 0.5pt i}\,\gamma^{\,b}_{aj}
+\sum^4_{a=1}\sum^4_{b=1}\gamma^{\,a}_{bj}\,\gamma^{\,b}_{a\kern 0.5pt i}
=2\sum^4_{a=1}\sum^4_{b=1}\gamma^{\,a}_{b\kern 0.5pt i}\,\gamma^{\,b}_{aj}
=8\,g_{ij}.
\mytag{7.12}
$$
Applying the covariant derivative $\nabla_{\!k}$ to both sides of the
equality \mythetag{7.12} and taking into account that $\nabla_{\!k}
\gamma^{\,a}_{b\kern 0.5pt i}=0$ and $\nabla_{\!k}\gamma^{\,b}_{aj}$, 
we derive $\nabla_{\!k}g_{ij}=0$.
\qed\enddemo
\mytheorem{7.3} A spinor connection $(\Gamma,\Alpha,\bar{\Alpha})$ of 
the bundle of Dirac spinors $DM$ concordant with the Dirac $\gamma$-field
and with the spin-metric tensor $\bold d$ is concordant with the chirality
operator $\bold H$ too, i\.\,e\. $\nabla\boldsymbol\gamma=0$ and $\nabla
\bold d=0$ imply $\nabla\bold H=0$.
\endproclaim
\demo{Proof} In order to prove this theorem we apply the identity
\mythetag{6.23}. Due to the previous theorem~\mythetheorem{7.2} from
$\nabla_{\!k}\gamma^{\,i}_{jm}=0$ it follows that $\nabla_{\!k}g_{ij}=0$ 
and $\nabla_{\!k}g^{ij}=0$. Then from $\nabla_{\!k}g^{ij}=0$ and
$\nabla_{\!k}\gamma^{\,i}_{jm}=0$ due to the formula \mythetag{6.11} 
we get $\nabla_{\!k}\gamma^{\,i\kern 0.5pt m}_{j}=0$. Now applying 
$\nabla_{\!k}$ to both sides of \mythetag{6.23} and using 
$\nabla\bold d=0$, we derive
$$
\sum^4_{r=1}\sum^4_{s=1}\bigl(\nabla_{\!k}H^a_r\,d^{\kern 0.5pt re}
\,d_{\kern 0.5pt bs}\,H^s_h+H^a_r\,d^{\kern 0.5pt re}\,d_{\kern 0.5pt
bs}\,\nabla_{\!k}H^s_h\,\bigr)+\,\nabla_{\!k}H^a_h\ H^e_b+H^a_h\
\nabla_{\!k}H^e_b=0.
$$
Let's multiply this equality by $H^b_q$ and sum it over the index $b$. 
As a result we get
$$
\hskip -2em
\gathered
\sum^4_{r=1}\nabla_{\!k}H^a_r\,d^{\kern 0.5pt re}\,d_{qh}+
\sum^4_{b=1}\sum^4_{r=1}\sum^4_{s=1}H^a_r\,d^{\kern 0.5pt re}\,
d_{qb}\,H^b_s\,\nabla_{\!k}H^s_h\,+\\
+\,\nabla_{\!k}H^a_h\,\delta^e_q+\sum^4_{b=1}H^a_h\,H^b_q\,
\nabla_{\!k}H^e_b=0.
\endgathered
\mytag{7.13}
$$
In deriving \mythetag{7.13} we used the formula \mythetag{6.17} and the
identity $\bold H^2=\idop$. In a coordinate form the identity
$\bold H^2=\idop$ is written as follows:
$$
\hskip -2em
\sum^4_{b=1}H^e_b\,H^b_h=\delta^e_h.
\mytag{7.14}
$$
The formula \mythetag{7.14} is easily derived from the matrix 
representation of the chirality operator \mythetag{5.4} or 
\mythetag{5.5}. From \mythetag{7.14} we easily derive
$$
\hskip -2em
\sum^4_{e=1}\sum^4_{b=1}H^e_b\,H^b_e=4.
\mytag{7.15}
$$
Applying the covariant derivative $\nabla_{\!k}$ to \mythetag{7.15},
we find that
$$
\hskip -2em
\sum^4_{e=1}\sum^4_{b=1}\bigl(\nabla_{\!k}H^e_b\,H^b_e+
H^e_b\,\nabla_{\!k}H^b_e\,\bigr)=2\sum^4_{e=1}\sum^4_{b=1}
H^b_e\,\nabla_{\!k}H^e_b=0.
\mytag{7.16}
$$
Now let's set $q=e$ in the equality \mythetag{7.13} and recall 
that the matrix $d^{\kern 0.5pt re}$ is inverse to the matrix 
$d_{\kern 0.5pt eb}$ when summing over the index $e$. As a result 
from \mythetag{7.13} we get
$$
\hskip -2em
\gathered
6\,\nabla_{\!k}H^a_h
+\sum^4_{e=1}\sum^4_{b=1}H^a_h\,H^b_e\,
\nabla_{\!k}H^e_b=0.
\endgathered
\mytag{7.17}
$$
Applying \mythetag{7.16} to \mythetag{7.17}, we see that the second term
in \mythetag{7.17} is zero. Hence, we have $\nabla_{\!k}H^a_h=0$. The
theorem~\mythetheorem{7.3} is proved.
\qed\enddemo
\mytheorem{7.4} A real spinor connection $(\Gamma,\Alpha,\bar{\Alpha})$ of 
the bundle of Dirac spinors $DM$ concordant with the Dirac $\gamma$-field
and with the spin-metric tensor $\bold d$ is concordant with the Dirac
form $\bold D$ as well, i\.\,e\. $\nabla\boldsymbol\gamma=0$ and $\nabla
\bold d=0$ imply $\nabla\bold D=0$.
\endproclaim
\demo{Proof} The proof of this theorem is based on the following 
identity for $\gamma$-symbols:
$$
\gathered
\sum^3_{m=0}\sum^3_{n=0}\gamma^{\,a}_{b\kern 0.5pt m}
\ g_{mn}\ \overline{\gamma^{\,e}_{h\kern 0.5pt n}}=
-\sum^4_{s=1}\sum^4_{\bar r=1}D_{s\bar h}\,d^{\kern 0.5pt sa}
\,D_{b\bar r}\,\bd^{\kern 0.5pt\bar r\bar e}\,+\\
+\sum^4_{r=1}\sum^4_{s=1}\sum^4_{\bar r=1}\sum^4_{q=1}
H^q_s\,D_{q\bar h}\,d^{\kern 0.5pt sa}\,
H^r_b\,D_{r\bar r}\,\bd^{\kern 0.5pt\bar r\bar e}\,-\\
-\sum^4_{s=1}\sum^4_{\bar r=1}d^{\kern 0.5pt sa}
\,D_{s\bar r}\,\bd^{\kern 0.5pt\bar r\bar e}\,D_{b\bar h}
+\sum^4_{r=1}\sum^4_{s=1}\sum^4_{\bar r=1}\sum^4_{q=1}
d^{\kern 0.5pt sa}\,H^r_s\,D_{r\bar r}
\,\bd^{\kern 0.5pt\bar r\bar e}\,H^q_b\,D_{q\bar h}.
\endgathered\qquad
\mytag{7.18}
$$
The identity \mythetag{7.18} is proved by direct calculations.
Let's apply the covariant derivative $\nabla_{\!k}$ to both sides
of the identity \mythetag{7.18}. When doing it we should remember 
that $\nabla\boldsymbol\gamma=0$ and $\nabla \bold d=0$ imply 
$\nabla\bold g=0$ and $\nabla\bold H=0$ due to the 
theorems~\mythetheorem{7.2} and \mythetheorem{7.3}. Moreover,
$\nabla\boldsymbol\gamma=0$ implies $\nabla\tau(\boldsymbol\gamma)=0$
since $(\Gamma,\Alpha,\bar{\Alpha})$ is assumed to be a real 
spinor connection. As a result, taking into account all the above 
arguments, from \mythetag{7.18} we derive the following identity
for the components of the Dirac form $\bold D$:
$$
\gather
-\sum^4_{s=1}\sum^4_{\bar r=1}\nabla_{\!k}D_{s\bar h}
\,d^{\kern 0.5pt sa}\,D_{b\bar r}\,\bd^{\kern 0.5pt\bar r\bar e}
-\sum^4_{s=1}\sum^4_{\bar r=1}D_{s\bar h}\,d^{\kern 0.5pt sa}
\,\nabla_{\!k}D_{b\bar r}\,\bd^{\kern 0.5pt\bar r\bar e}\,+\\
+\sum^4_{r=1}\sum^4_{s=1}\sum^4_{\bar r=1}\sum^4_{q=1}
d^{\kern 0.5pt sa}\bigl(H^q_s\,\nabla_{\!k}D_{q\bar h}\,
H^r_b\,D_{r\bar r}+H^q_s\,D_{q\bar h}\,H^r_b
\,\nabla_{\!k}D_{r\bar r}\,\bigr)\,\bd^{\kern 0.5pt\bar r\bar e}
\,-\\
-\sum^4_{s=1}\sum^4_{\bar r=1}d^{\kern 0.5pt sa}\,
\nabla_{\!k}D_{s\bar r}\,\bd^{\kern 0.5pt\bar r\bar e}
\,D_{b\bar h}-\sum^4_{s=1}\sum^4_{\bar r=1}d^{\kern 0.5pt sa}
\,D_{s\bar r}\,\bd^{\kern 0.5pt\bar r\bar e}
\,\nabla_{\!k}D_{b\bar h}\,+\\
+\sum^4_{r=1}\sum^4_{s=1}\sum^4_{\bar r=1}\sum^4_{q=1}
d^{\kern 0.5pt sa}\bigl(H^r_s\,\nabla_{\!k}D_{r\bar r}
\,H^q_b\,D_{q\bar h}+H^r_s\,D_{r\bar r}\,H^q_b\,
\nabla_{\!k}D_{q\bar h}\,\bigr)\,\bd^{\kern 0.5pt\bar r\bar e}=0.
\endgather
$$
Note that each term in the above identity contains the sum over
the index $s$ and the factor $d^{\kern 0.5pt sa}$, where $a$ is
a free index. Similarly, each term contains the sum over the index 
$\bar r$ and the factor $\bd^{\kern 0.5pt\bar r\bar e}$, where 
$\bar e$ is a free index. The quantities $d^{\kern 0.5pt sa}$ and
$\bd^{\kern 0.5pt\bar r\bar e}$ form two non-degenerate matrices.
Therefore, we can cancel them in the above equality and simultaneously
omit the sums over $s$ and $\bar r$:  
$$
\gathered
-\nabla_{\!k}D_{s\bar h}\,D_{b\bar r}
-D_{s\bar h}\,\nabla_{\!k}D_{b\bar r}+\sum^4_{r=1}
\sum^4_{q=1}H^q_s\,\nabla_{\!k}D_{q\bar h}\,
H^r_b\,D_{r\bar r}\,+\\
+\sum^4_{r=1}\sum^4_{q=1}H^q_s\,D_{q\bar h}\,H^r_b
\,\nabla_{\!k}D_{r\bar r}-\nabla_{\!k}D_{s\bar r}
\,D_{b\bar h}-D_{s\bar r}\,\nabla_{\!k}D_{b\bar h}\,+\\
+\sum^4_{r=1}\sum^4_{q=1}H^r_s\,\nabla_{\!k}D_{r\bar r}
\,H^q_b\,D_{q\bar h}+\sum^4_{r=1}\sum^4_{q=1}H^r_s
\,D_{r\bar r}\,H^q_b\,\nabla_{\!k}D_{q\bar h}=0.
\endgathered\qquad
\mytag{7.19}
$$
Now let's multiply \mythetag{7.19} by $D^{\kern 0.5pt b\bar h}$ and sum
over the indices $b$ and $\bar h$. The quantities $D^{\kern 0.5pt b\bar h}$
are defined through the formula \mythetag{6.24}. Taking into account
\mythetag{6.25}, we derive
$$
\gathered
-5\,\nabla_{\!k}D_{s\bar r}+\sum^4_{r=1}\sum^4_{q=1}
\sum^4_{b=1}\sum^4_{\bar h=1}H^q_s\,\nabla_{\!k}D_{q\bar h}
\,D^{\kern 0.5pt b\bar h}\,H^r_b\,D_{r\bar r}\,-\\
-\sum^4_{b=1}\sum^4_{\bar h=1}
D_{s\bar r}\,\nabla_{\!k}D_{b\bar h}\,D^{\kern 0.5pt b\bar h}
+\sum^4_{r=1}\sum^4_{b=1}H^r_s\,\nabla_{\!k}D_{r\bar r}
\,H^b_b\,+\\
+\sum^4_{r=1}\sum^4_{q=1}\sum^4_{b=1}\sum^4_{\bar h=1}
H^r_s\,D_{r\bar r}\,H^q_b\,\nabla_{\!k}D_{q\bar h}
\,D^{\kern 0.5pt b\bar h}=0.
\endgathered\qquad
\mytag{7.20}
$$
In deriving \mythetag{7.20}, apart from \mythetag{6.25}, we used
the equality \mythetag{7.14}. Remember, that $\bold H$ is a traceless
operator (see its matrix \mythetag{5.4}). This means that
$$
\hskip -2em
\tr\bold H=\sum^4_{b=1}H^b_b=0.
\mytag{7.21}
$$
Due to \mythetag{7.21} the fourth term in \mythetag{7.20} is zero.
The third term in the left hand side of \mythetag{7.20} is zero 
too. This fact is proved with the use of the formula \mythetag{6.24}:
$$
\gather
0=\sum^4_{b=1}\sum^4_{\bar h=1}\nabla_{\!k}(D_{b\bar h}
\,D^{\kern 0.5pt b\bar h})=\sum^4_{b=1}\sum^4_{\bar h=1}
\sum^4_{a=1}\sum^4_{\bar a=1}\nabla_{\!k}(D_{b\bar h}
\,d^{\kern 0.5pt b\kern 0.5pt a}\,D_{a\bar a}
\,\bd^{\kern 0.5pt\bar a\bar h})=\\
=\sum^4_{b=1}\sum^4_{\bar h=1}\sum^4_{a=1}\sum^4_{\bar a=1}
\nabla_{\!k}D_{b\bar h}\,d^{\kern 0.5pt b\kern 0.5pt a}
\,D_{a\bar a}\,\bd^{\kern 0.5pt\bar a\bar h}
+\sum^4_{b=1}\sum^4_{\bar h=1}\sum^4_{a=1}\sum^4_{\bar a=1}
D_{b\bar h}\,d^{\kern 0.5pt b\kern 0.5pt a}
\,\nabla_{\!k}D_{a\bar a}\,\bd^{\kern 0.5pt\bar a\bar h}=\\
=\sum^4_{b=1}\sum^4_{\bar h=1}\sum^4_{a=1}\sum^4_{\bar a=1}
\nabla_{\!k}D_{b\bar h}\,d^{\kern 0.5pt b\kern 0.5pt a}
\,D_{a\bar a}\,\bd^{\kern 0.5pt\bar a\bar h}
+\sum^4_{b=1}\sum^4_{\bar h=1}\sum^4_{a=1}\sum^4_{\bar a=1}
\nabla_{\!k}D_{a\bar a}\,d^{\kern 0.5pt ab}\,D_{b\bar h}\,
\bd^{\kern 0.5pt\bar h\bar a}=\\
=2\,\sum^4_{b=1}\sum^4_{\bar h=1}\sum^4_{a=1}\sum^4_{\bar a=1}
\nabla_{\!k}D_{b\bar h}\,d^{\kern 0.5pt b\kern 0.5pt a}
\,D_{a\bar a}\,\bd^{\kern 0.5pt\bar a\bar h}
=2\,\sum^4_{b=1}\sum^4_{\bar h=1}\nabla_{\!k}D_{b\bar h}
\,D^{\kern 0.5pt b\bar h}.
\endgather
$$
In the above calculations we used \mythetag{6.25} and the 
skew-symmetry of $\bold d$ and $\bbd$. The last term in 
the left hand side of \mythetag{7.20} is also zero. This
fact is derived from \mythetag{7.21}:
$$
\gather
0=\sum^4_{b=1}\nabla_{\!k}H^b_b=\sum^4_{q=1}\sum^4_{b=1}
\sum^4_{\bar h=1}\nabla_{\!k}(H^q_b\,D_{q\bar h}
\,D^{\kern 0.5pt b\bar h})=\sum^4_{q=1}\sum^4_{b=1}
\sum^4_{\bar h=1}\sum^4_{a=1}\sum^4_{\bar a=1}\bigl(H^q_b\,
\times\\
\times\,\nabla_{\!k}D_{q\bar h}\,d^{\kern 0.5pt b\kern 0.5pt a}
\,D_{a\bar a}\,\bd^{\kern 0.5pt\bar a\bar h}+H^q_b\,D_{q\bar h}
\,d^{\kern 0.5pt b\kern 0.5pt a}\,\nabla_{\!k}D_{a\bar a}
\,\bd^{\kern 0.5pt\bar a\bar h}\bigr)=\sum^4_{q=1}\sum^4_{b=1}
\sum^4_{\bar h=1}\sum^4_{a=1}\sum^4_{\bar a=1}H^q_b\,
\times\\
\times\,\nabla_{\!k}D_{q\bar h}\,d^{\kern 0.5pt b\kern 0.5pt a}
\,D_{a\bar a}\,\bd^{\kern 0.5pt\bar a\bar h}+\sum^4_{q=1}
\sum^4_{b=1}\sum^4_{\bar h=1}\sum^4_{a=1}\sum^4_{\bar a=1}
H^a_b\,\nabla_{\!k}D_{a\bar a}\,d^{\kern 0.5pt b\kern 0.5pt q}
\,D_{q\bar h}\,\bd^{\kern 0.5pt\bar h\bar a}=\\
=2\,\sum^4_{q=1}\sum^4_{b=1}\sum^4_{\bar h=1}\sum^4_{a=1}
\sum^4_{\bar a=1}H^q_b\,\nabla_{\!k}D_{q\bar h}
\,d^{\kern 0.5pt b\kern 0.5pt a}\,D_{a\bar a}
\,\bd^{\kern 0.5pt\bar a\bar h}=2\,\sum^4_{q=1}\sum^4_{b=1}
\sum^4_{\bar h=1}H^q_b\,\nabla_{\!k}D_{q\bar h}
\,D^{\kern 0.5pt b\bar h}.
\endgather
$$
In these calculations we used \mythetag{6.25}, the skew-symmetry 
$\bold d$ and $\bbd$, and the equality
$$
\hskip -2em
\sum^4_{b=1}H^q_b\,d^{\kern 0.5pt b\kern 0.5pt a}=
\sum^4_{b=1}d^{\kern 0.5pt q\kern 0.5pt b}\,H^a_b.
\mytag{7.22}
$$
The equality \mythetag{7.22} is easily derived from \mythetag{6.17}.
And finally, we need to transform the second term in the left hand
side of \mythetag{7.20}:
$$
\pagebreak
\hskip 0.1em
\gathered
\sum^4_{r=1}\sum^4_{q=1}\sum^4_{b=1}\sum^4_{\bar h=1}H^q_s
\,\nabla_{\!k}D_{q\bar h}\,D^{\kern 0.5pt b\bar h}\,H^r_b
\,D_{r\bar r}=-\sum^4_{\bar s=1}\sum^4_{q=1}\sum^4_{b=1}
\sum^4_{\bar h=1}H^q_s\,\nabla_{\!k}D_{q\bar h}
\,\times\\
\times\,D^{\kern 0.5pt b\bar h}\,D_{b\bar s}\,
\overline{H^{\bar s}_{\bar r}}=-\sum^4_{q=1}\sum^4_{\bar h=1}
H^q_s\,\nabla_{\!k}D_{q\bar h}\,\overline{H^{\bar h}_{\bar r}}=
-\sum^4_{q=1}\sum^4_{\bar h=1}
H^q_s\,\nabla_{\!k}(D_{q\bar h}\,\overline{H^{\bar h}_{\bar r}})=\\
=\sum^4_{q=1}\sum^4_{r=1}H^q_s\,\nabla_{\!k}(H^r_q\,D_{r\bar r})
=\sum^4_{q=1}\sum^4_{r=1}H^q_s\,H^r_q\,\nabla_{\!k}D_{r\bar r}=
\nabla_{\!k}D_{s\bar r}.
\endgathered
\hskip -1.7em
\mytag{7.23}
$$
In deriving \mythetag{7.23} we used the equalities \mythetag{6.20},
\mythetag{6.25}, and \mythetag{7.14}. Now, substituting \mythetag{7.23}
back into \mythetag{7.20} and recalling that the third, the fourth,
and the fifth terms in the left hand side of this equality do vanish,
we find that \mythetag{7.20} is reduced to $\nabla_{\!k}D_{s\bar r}=0$.
Thus, $\nabla\bold D=0$. The theorem~\mythetheorem{7.4} is proved.
\qed\enddemo
\mydefinition{7.8} A real spinor connection $(\Gamma,\Alpha,\bar{\Alpha})$
of the bundle of Dirac spinors $DM$ concordant with the spin-metric tensor
$\bold d$ and with Dirac $\gamma$-field $\boldsymbol\gamma$ is called a
{\it metric connection}.
\enddefinition
     According to the above theorems~\mythetheorem{7.2}, \mythetheorem{7.3},
\mythetheorem{7.4}, a metric connection $(\Gamma,\Alpha,\bar{\Alpha})$ is
concordant with the metric tensor $\bold g$, with the chirality operator
$\bold H$, with the Dirac form $\bold D$, and with many other fields 
produced from these basic fields.
\head
8. Explicit formulas.     
\endhead     
     Let $(\Gamma,\Alpha,\bar{\Alpha})$ be a metric connection of the
bundle of Dirac spinors $DM$. Since $\nabla\bold g=0$ for this connection,
its $\Gamma$-components are uniquely determined by the torsion tensor
$\bold T$. They are given by the explicit formula \mythetag{4.34}. If
$\bold T=0$, this formula reduces to \mythetag{4.35}. Our goal in this
section is to study $\Alpha$-components of a metric connection. Its
$\bar{\Alpha}$-components then are expressed through $\Alpha$-components
according to the formula \mythetag{7.11}.\par
     In order to derive an explicit formula for the $\Alpha$-components 
of a metric connection we introduce some auxiliary objects. Using the
chirality operator $\bold H$, we define two projector operators
$\bulletBH$ and $\circBH$ by means of the following formulas:
$$
\xalignat 2
&\hskip -2em
\bulletBH=\frac{\idop+\bold H}{2},
&&\circBH=\frac{\idop-\bold H}{2}.
\mytag{8.1}
\endxalignat
$$
Their components are given by the formulas
$$
\xalignat 2
&\hskip -2em
\bulletH^i_j=\frac{\delta^i_j+H^i_j}{2},
&&\circH^i_j=\frac{\delta^i_j-H^i_j}{2},
\mytag{8.2}
\endxalignat
$$
Due to \mythetag{8.1} and \mythetag{8.2} we have the expansions
$$
\xalignat 2
&\hskip -2em
\idop=\bulletBH+\circBH,
&&\delta^i_j=\bulletH^i_j+\circH^i_j.
\mytag{8.3}
\endxalignat
$$
Now, relying on the expansion \mythetag{8.3}, we introduce the following
quantities:
$$
\xalignat 2
&\hskip -2em
\bulbulgamma^{\,i}_{\!jm}=\sum^4_{r=1}\sum^4_{s=1}\bulletH^i_r
\,\bulletH^s_j\,\gamma^{\,r}_{sm},
&&\circcircgamma^{\,i}_{\!jm}=\sum^4_{r=1}\sum^4_{s=1}\circH^i_r
\,\circH^s_j\,\gamma^{\,r}_{sm},\quad
\mytag{8.4}\\
&\hskip -2em
\circbulgamma^{\,i}_{\!jm}=\sum^4_{r=1}\sum^4_{s=1}\circH^i_r
\,\bulletH^s_j\,\gamma^{\,r}_{sm},
&&\bulcircgamma^{\,i}_{\!jm}=\sum^4_{r=1}\sum^4_{s=1}\bulletH^i_r
\,\circH^s_j\,\gamma^{\,r}_{sm}.\quad
\mytag{8.5}
\endxalignat
$$
From \mythetag{8.4}, \mythetag{8.5}, and from the expansion \mythetag{8.3}
we derive the reconstruction formula
$$
\hskip -2em
\gamma^{\,i}_{\!jm}=\bulbulgamma^{\,i}_{\!jm}+\bulcircgamma^{\,i}_{\!jm}
+\circbulgamma^{\,i}_{\!jm}+\circcircgamma^{\,i}_{\!jm}.
\mytag{8.6}
$$
By means of direct calculations we find that the quantities 
\mythetag{8.4} are identically zero: $\bulbulgamma^{\,i}_{\!jm}=0$ 
and $\circcircgamma^{\,i}_{\!jm}=0$. Therefore, the expansion 
\mythetag{8.6} reduces to the following one:
$$
\gamma^{\,i}_{\!jm}=\bulcircgamma^{\,i}_{\!jm}+\circbulgamma^{\,i}_{\!jm}.
$$\par
    Now let's proceed with the $\Alpha$-components of a metric connection
$(\Gamma,\Alpha,\bar{\Alpha})$. By analogy to the formulas \mythetag{8.4}
and \mythetag{8.5} we introduce the following quantities:
$$
\xalignat 2
&\hskip -2em
\bulbulA^i_{kj}=\sum^4_{r=1}\sum^4_{s=1}\bulletH^i_r
\,\bulletH^s_j\,\Alpha^r_{ks},
&&\circcircA^i_{kj}=\sum^4_{r=1}\sum^4_{s=1}\circH^i_r
\,\circH^s_j\,\Alpha^r_{ks},\quad
\mytag{8.7}\\
&\hskip -2em
\circbulA^i_{kj}=\sum^4_{r=1}\sum^4_{s=1}\circH^i_r
\,\bulletH^s_j\,\Alpha^r_{ks},
&&\bulcircA^i_{kj}=\sum^4_{r=1}\sum^4_{s=1}\bulletH^i_r
\,\circH^s_j\,\Alpha^r_{ks}.\quad
\mytag{8.8}
\endxalignat
$$
In this case the quantities \mythetag{8.7} are not necessarily 
zero. Therefore, we have all of the four terms \mythetag{8.7} and 
\mythetag{8.8} in the expansion
$$
\hskip -2em
\Alpha^i_{kj}=\bulbulA^i_{kj}+\bulcircA^i_{kj}
+\circbulA^i_{kj}+\circcircA^i_{kj}.
\mytag{8.9}
$$\par
     A metric connection $(\Gamma,\Alpha,\bar{\Alpha})$ is concordant
with $\boldsymbol\gamma$ and $\bold H$. Hence, from $\nabla\boldsymbol
\gamma=0$ and $\nabla\bold H=0$ we derive the following equalities:
$$
\xalignat 4
&\nabla_{\!k}\bulletH^i_j=0,
&&\nabla_{\!k}\circH^i_j=0, 
&&\nabla_{\!k}\bulcircgamma^{\,i}_{\!jm}=0,
&&\nabla_{\!k}\circbulgamma^{\,i}_{\!jm}=0.
\quad\qquad
\mytag{8.10}
\endxalignat
$$
Applying the formula \mythetag{7.10} to the first equality \mythetag{8.10},
we get
$$
\hskip -2em
L_{\boldsymbol\Upsilon_{\!k}}(\bulletH^i_s)+\sum^4_{r=1}
\Alpha^i_{k\kern 0.5pt r}
\,\bulletH^r_s-\sum^4_{r=1}\Alpha^r_{ks}\,\bulletH^i_r=0.
\mytag{8.11}
$$
Let's multiply \mythetag{8.11} by $\circH^s_j$ and sum it over the index 
$s$:
$$
\sum^4_{s=1}\circH^s_j\,L_{\boldsymbol\Upsilon_{\!k}}(\bulletH^i_s)
+\sum^4_{s=1}\sum^4_{r=1}\Alpha^i_{k\kern 0.5pt r}\,\bulletH^r_s\,
\circH^s_j-\sum^4_{s=1}\sum^4_{r=1}\bulletH^i_r\,\circH^s_j\,\Alpha^r_{ks}=0.
\qquad
\mytag{8.12}
$$
The projectors \mythetag{8.1} are complementary to each other. 
Therefore their product is zero: $\bulletBH\cdot\circBH=0$. Applying 
this equality to \mythetag{8.12} we derive
$$
\hskip -2em
\bulcircA^i_{kj}=\sum^4_{s=1}\circH^s_j
\,L_{\boldsymbol\Upsilon_{\!k}}(\bulletH^i_s).
\mytag{8.13}
$$
In a similar way from the second equality \mythetag{8.10}, we obtain
$$
\hskip -2em
\circbulA^i_{kj}=\sum^4_{s=1}\bulletH^s_j
\,L_{\boldsymbol\Upsilon_{\!k}}(\circH^i_s).
\mytag{8.14}
$$
Now let's apply the formula \mythetag{7.10} to the third and to the
fourth equalities \mythetag{8.10}. As a result we get the following
two formulas:
$$
\align
&\hskip -2em
L_{\boldsymbol\Upsilon_{\!k}}(\bulcircgamma^{\,a}_{\!bm})
+\sum^4_{r=1}\Alpha^a_{k\kern 0.5pt r}\,\bulcircgamma^{\,r}_{\!bm}
-\sum^4_{r=1}\Alpha^r_{kb}\,\bulcircgamma^{\,a}_{\!rm}
=\sum^3_{r=0}\Gamma^r_{km}\,\bulcircgamma^{\,a}_{\!br},
\mytag{8.15}\\
&\hskip -2em
L_{\boldsymbol\Upsilon_{\!k}}(\circbulgamma^{\,a}_{\!bm})
+\sum^4_{r=1}\Alpha^a_{k\kern 0.5pt r}\,\circbulgamma^{\,r}_{\!bm}
-\sum^4_{r=1}\Alpha^r_{kb}\,\circbulgamma^{\,a}_{\!rm}
=\sum^3_{r=0}\Gamma^r_{km}\,\circbulgamma^{\,a}_{\!br}.
\mytag{8.16}
\endalign
$$
In order to transform \mythetag{8.15} and \mythetag{8.16} we need some
identities for the $\gamma$-symbols \mythetag{8.5}. From the identity
\mythetag{6.23} we derive two formulas
$$
\align
&\hskip -2em
\sum^3_{m=0}\sum^3_{n=0}\bulcircgamma^{\,a}_{\!bm}\ g^{mn}
\ \bulcircgamma^{\,e}_{\!hn}=2\,\bulletd^{\kern 0.5pt ae}\,
\circd_{\kern 0.5pt bh},
\mytag{8.17}\\
&\hskip -2em
\sum^3_{m=0}\sum^3_{n=0}\circbulgamma^{\,a}_{\!bm}\ g^{mn}
\ \circbulgamma^{\,e}_{\!hn}=2\,\circd^{\kern 0.5pt ae}\,
\bulletd_{\kern 0.5pt bh}.
\mytag{8.18}
\endalign
$$
In the formulas \mythetag{8.17} and \mythetag{8.18} we used the 
following notations:
$$
\xalignat 2
&\hskip -2em
\bulletd^{\kern 0.5pt ae}=\sum^4_{r=1}\sum^4_{s=1}
\bulletH^a_r\ d^{\kern 0.5pt rs}\,\bulletH^e_s,
&&\circd^{\kern 0.5pt ae}=\sum^4_{r=1}\sum^4_{s=1}
\circH^a_r\ d^{\kern 0.5pt rs}\,\circH^e_s,\qquad\\
\vspace{-1.5ex}
&&&\mytag{8.19}\\
\vspace{-1.5ex}
&\hskip -2em
\bulletd_{\kern 0.5pt bh}=\sum^4_{r=1}\sum^4_{s=1}
\bulletH^r_b\ d_{\kern 0.5pt rs}\,\bulletH^s_h,
&&\circd_{\kern 0.5pt bh}=\sum^4_{r=1}\sum^4_{s=1}
\circH^r_b\ d_{\kern 0.5pt rs}\,\circH^s_h.\qquad
\endxalignat
$$
The matrices \mythetag{8.19} are skew-symmetric and degenerate.
However, in some cases they can be used for raising and lowering
spinor indices:
$$
\xalignat 2
&\hskip -2em
\bulbulA_{kij}=\sum^4_{s=1}\bulbulA^s_{ki}\,\bulletd_{\kern 0.5pt sj},
&&\circcircA_{kij}=\sum^4_{s=1}\circcircA^s_{ki}\,\circd_{\kern 0.5pt sj},\\
&\hskip -2em
\bulbulA^i_{kj}=\sum^4_{s=1}\bulbulA_{kjs}\,\bulletd^{\kern 0.5pt si},
&&\circcircA^i_{kj}=\sum^4_{s=1}\circcircA_{kjs}\,\circd^{\kern 0.5pt si}.
\mytag{8.20}
\endxalignat
$$
Applying the formula \mythetag{8.17} to \mythetag{8.15}, we derive the
following equality:
$$
\hskip -2em
\gathered
\frac{1}{2}\sum^3_{m=0}\sum^3_{n=0}
L_{\boldsymbol\Upsilon_{\!k}}(\bulcircgamma^{\,a}_{\!bm})\,
g^{mn}\,\bulcircgamma^{\,e}_{\!hn}+\sum^4_{r=1}
\Alpha^a_{k\kern 0.5pt r}\,\bulletd^{\kern 0.5pt re}\,
\circd_{\kern 0.5pt bh}\,-\\
-\sum^4_{r=1}\Alpha^r_{kb}\,\bulletd^{\kern 0.5pt ae}\,
\circd_{\kern 0.5pt rh}
=\frac{1}{2}\sum^3_{m=0}\sum^3_{n=0}\sum^3_{r=0}
\Gamma^r_{km}\,\bulcircgamma^{\,a}_{\!br}\,g^{mn}
\,\bulcircgamma^{\,e}_{\!hn},
\endgathered\qquad
\mytag{8.21}
$$
Let's multiply \mythetag{8.21} by $\bulletH^q_a\,\circH^b_j$ and sum
over the indices $a$ and $b$. Then we get
$$
\gathered
\frac{1}{2}\sum^3_{m=0}\sum^3_{n=0}\sum^4_{a=1}\sum^4_{b=1}
L_{\boldsymbol\Upsilon_{\!k}}(\bulcircgamma^{\,a}_{\!bm})\,
g^{mn}\,\bulcircgamma^{\,e}_{\!hn}\,\bulletH^q_a\,\circH^b_j
+\sum^4_{r=1}\bulbulA^q_{k\kern 0.5pt r}\,\bulletd^{\kern 0.5pt re}
\,\circd_{\kern 0.5pt jh}\,-\\
-\,\circcircA_{kjh}\,\bulletd^{\kern 0.5pt qe}
=\frac{1}{2}\sum^3_{m=0}\sum^3_{n=0}\sum^3_{r=0}\Gamma^r_{km}
\,\bulcircgamma^{\,q}_{\!jr}\,g^{mn}\,\bulcircgamma^{\,e}_{\!hn}.
\endgathered\qquad
\mytag{8.22}
$$
Now let's multiply \mythetag{8.22} by $\bulletd_{\kern 0.5pt eq}$ and
sum over the indices $e$ and $q$. This yields
$$
\gathered
2\,\circcircA_{kjh}
=\frac{1}{2}\sum^3_{m=0}\sum^3_{n=0}\sum^4_{a=1}\sum^4_{b=1}
\sum^4_{e=1}L_{\boldsymbol\Upsilon_{\!k}}(\bulcircgamma^{\,a}_{\!bm})
\,g^{mn}\,\bulcircgamma^{\,e}_{\!hn}\,\circH^b_j
\,\bulletd_{\kern 0.5pt ea}\,+\\
+\,\sum^4_{r=1}\bulbulA^r_{k\kern 0.5pt r}
\,\circd_{\kern 0.5pt jh}
-\,\frac{1}{2}\sum^3_{m=0}\sum^3_{n=0}\sum^3_{r=0}\sum^4_{e=1}
\sum^4_{q=1}\Gamma^r_{km}\,\bulcircgamma^{\,q}_{\!jr}\,g^{mn}
\,\bulcircgamma^{\,e}_{\!hn}\,\bulletd_{\kern 0.5pt eq}.
\endgathered\qquad
\mytag{8.23}
$$
Then let's raise the lower index $h$ in \mythetag{8.23} by means of
the matrix $\circd^{\kern 0.5pt hi}$, i\.\,e\. applying the second
formula \mythetag{8.20} to $\circcircA_{kjh}$ in the left hand side.
As a result we get
$$
\hskip -2em
\gathered
\circcircA^i_{kj}
=\frac{1}{4}\sum^3_{m=0}\sum^3_{n=0}\sum^4_{a=1}\sum^4_{b=1}
L_{\boldsymbol\Upsilon_{\!k}}(\bulcircgamma^{\,a}_{\!bm})
\,g^{mn}\,\circbulgamma^{\,i}_{\!an}\,\circH^b_j\,+\\
+\,\frac{1}{2}\sum^4_{r=1}\bulbulA^r_{k\kern 0.5pt r}
\,\circH^i_j-\,\frac{1}{4}\sum^3_{m=0}\sum^3_{n=0}\sum^3_{r=0}
\sum^4_{q=1}\Gamma^r_{km}\,\bulcircgamma^{\,q}_{\!jr}\,g^{mn}
\,\circbulgamma^{\,i}_{\!qn}.
\endgathered\qquad
\mytag{8.24}
$$\par
     Acting in a similar way, i\.\,e\. applying the formula 
\mythetag{8.18} to \mythetag{8.16} and then performing some
calculations analogous to the above ones, we derive the formula
$$
\hskip -2em
\gathered
\bulbulA^i_{kj}
=\frac{1}{4}\sum^3_{m=0}\sum^3_{n=0}\sum^4_{a=1}\sum^4_{b=1}
L_{\boldsymbol\Upsilon_{\!k}}(\circbulgamma^{\,a}_{\!bm})
\,g^{mn}\,\bulcircgamma^{\,i}_{\!an}\,\bulletH^b_j\,+\\
+\,\frac{1}{2}\sum^4_{r=1}\circcircA^r_{k\kern 0.5pt r}
\,\bulletH^i_j-\,\frac{1}{4}\sum^3_{m=0}\sum^3_{n=0}\sum^3_{r=0}
\sum^4_{q=1}\Gamma^r_{km}\,\circbulgamma^{\,q}_{\!jr}\,g^{mn}
\,\bulcircgamma^{\,i}_{\!qn}.
\endgathered\qquad
\mytag{8.25}
$$
Remember that for a metric connection $\nabla\bold d=0$ and 
$\nabla\bold H=0$. Then from \mythetag{8.2} and \mythetag{8.19}
we get the following equalities for $\bulletd_{\kern 0.5pt ab}$
and $\circd_{\kern 0.5pt ab}$:
$$
\xalignat 2
&\hskip -2em
\nabla_{\!k}\bulletd_{\kern 0.5pt ab}=0,
&&\nabla_{\!k}\circd_{\kern 0.5pt ab}=0.
\mytag{8.26}
\endxalignat
$$
In an expanded form these equalities \mythetag{8.26} are written as
$$
\align
&\hskip -2em
L_{\boldsymbol\Upsilon_{\!k}}(\bulletd_{\kern 0.5pt ab})
-\sum^4_{r=1}\Alpha^r_{ka}\,\bulletd_{\kern 0.5pt rb}
-\sum^4_{r=1}\Alpha^r_{kb}\,\bulletd_{\kern 0.5pt ar}=0,
\mytag{8.27}\\
&\hskip -2em
L_{\boldsymbol\Upsilon_{\!k}}(\circd_{\kern 0.5pt ab})
-\sum^4_{r=1}\Alpha^r_{ka}\,\circd_{\kern 0.5pt rb}
-\sum^4_{r=1}\Alpha^r_{kb}\,\circd_{\kern 0.5pt ar}=0.
\mytag{8.28}
\endalign
$$
Let's multiply $\mythetag{8.27}$ by $\bulletd^{\kern 0.5pt b\kern 0.5pt a}$
and multiply $\mythetag{8.28}$ by $\circd^{\kern 0.5pt b\kern 0.5pt a}$. 
Then let's sum both equalities over the indices $a$ and $b$. As a result
we get
$$
\xalignat 2
&\sum^4_{r=1}\bulbulA^r_{kr}=\frac{1}{2}\sum^4_{a=1}\sum^4_{b=1}
L_{\boldsymbol\Upsilon_{\!k}}(\bulletd_{\kern 0.5pt ab})\,
\bulletd^{\kern 0.5pt b\kern 0.5pt a},
&&\sum^4_{r=1}\circcircA^r_{kr}=\frac{1}{2}\sum^4_{a=1}\sum^4_{b=1}
L_{\boldsymbol\Upsilon_{\!k}}(\circd_{\kern 0.5pt ab})\,
\circd^{\kern 0.5pt b\kern 0.5pt a}.
\endxalignat 
$$
Substituting these formulas back into \mythetag{8.24} and \mythetag{8.25},
we derive
$$
\align
&\hskip -7em
\gathered
\circcircA^i_{kj}
=\frac{1}{4}\sum^3_{m=0}\sum^3_{n=0}\sum^4_{a=1}\sum^4_{b=1}
L_{\boldsymbol\Upsilon_{\!k}}(\bulcircgamma^{\,a}_{\!bm})
\,g^{mn}\,\circbulgamma^{\,i}_{\!an}\,\circH^b_j\,+\\
+\,\frac{1}{4}\sum^4_{a=1}\sum^4_{b=1}
L_{\boldsymbol\Upsilon_{\!k}}(\bulletd_{\kern 0.5pt ab})
\,\bulletd^{\kern 0.5pt b\kern 0.5pt a}\,\circH^i_j
-\,\frac{1}{4}\sum^3_{m=0}\sum^3_{n=0}\sum^3_{r=0}
\sum^4_{q=1}\Gamma^r_{km}\,\bulcircgamma^{\,q}_{\!jr}\,g^{mn}
\,\circbulgamma^{\,i}_{\!qn},
\endgathered\hskip -1em
\mytag{8.29}\\
&\hskip -7em
\gathered
\bulbulA^i_{kj}
=\frac{1}{4}\sum^3_{m=0}\sum^3_{n=0}\sum^4_{a=1}\sum^4_{b=1}
L_{\boldsymbol\Upsilon_{\!k}}(\circbulgamma^{\,a}_{\!bm})
\,g^{mn}\,\bulcircgamma^{\,i}_{\!an}\,\bulletH^b_j\,+\\
+\,\frac{1}{4}\sum^4_{a=1}\sum^4_{b=1}
L_{\boldsymbol\Upsilon_{\!k}}(\circd_{\kern 0.5pt ab})\,
\circd^{\kern 0.5pt b\kern 0.5pt a}\,\bulletH^i_j
-\,\frac{1}{4}\sum^3_{m=0}\sum^3_{n=0}\sum^3_{r=0}
\sum^4_{q=1}\Gamma^r_{km}\,\circbulgamma^{\,q}_{\!jr}\,g^{mn}
\,\bulcircgamma^{\,i}_{\!qn}.
\endgathered\hskip -1em
\mytag{8.30}
\endalign
$$
The formula \mythetag{8.29} is an analog of the formula \mythetag{4.57},
while \mythetag{8.30} is an analog of \mythetag{4.62}. Moreover, these
formulas are reduced to the corresponding formulas \mythetag{4.57} and
\mythetag{4.62} if we choose a spinor frame $(U,\,\boldsymbol\Psi_1,
\,\boldsymbol\Psi_2,\,\boldsymbol\Psi_3,\,\boldsymbol\Psi_4)$ concordant
with the expansion \mythetag{5.1}. This fact leads to the following 
theorem.
\mytheorem{8.1} Any metric connection $(\Gamma,\Alpha,\bar{\Alpha})$
of the bundle of Dirac spinors $DM$ is a unique extension of some metric
connection of the chiral bundle $SM$.
\endproclaim
    According to the formula \mythetag{8.9}, the final step in deriving 
an explicit formula for the $\Alpha$-components of a metric connection
$(\Gamma,\Alpha,\bar{\Alpha})$ is to add the above four formulas
\mythetag{8.13}, \mythetag{8.14}, \mythetag{8.29}, and \mythetag{8.30}.
This yields
$$
\gathered
\Alpha^i_{kj}=
\sum^4_{b=1}\circH^b_j\,L_{\boldsymbol\Upsilon_{\!k}}(\bulletH^i_b)+
\sum^4_{b=1}\bulletH^b_j\,L_{\boldsymbol\Upsilon_{\!k}}(\circH^i_b)\,+\\
+\sum^3_{m=0}\sum^3_{n=0}\sum^4_{a=1}\sum^4_{b=1}
\frac{L_{\boldsymbol\Upsilon_{\!k}}(\bulcircgamma^{\,a}_{\!bm})
\,g^{mn}\,\circbulgamma^{\,i}_{\!an}\,\circH^b_j+
L_{\boldsymbol\Upsilon_{\!k}}(\circbulgamma^{\,a}_{\!bm})
\,g^{mn}\,\bulcircgamma^{\,i}_{\!an}\,\bulletH^b_j}{4}\,-\\
-\sum^3_{m=0}\sum^3_{n=0}\sum^3_{r=0}\sum^4_{q=1}
\frac{\Gamma^r_{km}\,\bulcircgamma^{\,q}_{\!jr}\,g^{mn}
\,\circbulgamma^{\,i}_{\!qn}+\Gamma^r_{km}\,\circbulgamma^{\,q}_{\!jr}
\,g^{mn}\,\bulcircgamma^{\,i}_{\!qn}}{4}\,+\\
+\,\frac{1}{4}\sum^4_{a=1}\sum^4_{b=1}
L_{\boldsymbol\Upsilon_{\!k}}(\bulletd_{\kern 0.5pt ab})
\,\bulletd^{\kern 0.5pt b\kern 0.5pt a}\,\circH^i_j
+\frac{1}{4}\sum^4_{a=1}\sum^4_{b=1}
L_{\boldsymbol\Upsilon_{\!k}}(\circd_{\kern 0.5pt ab})\,
\circd^{\kern 0.5pt b\kern 0.5pt a}\,\bulletH^i_j.
\endgathered\qquad
\mytag{8.31}
$$
This formula \mythetag{8.31} can be simplified a little bit. For 
this purpose we need the following identities derived from the formula
\mythetag{6.23}:
$$
\align
&\hskip -2em
\sum^3_{m=0}\sum^3_{n=0}\bulcircgamma^{\,a}_{\!bm}\,g^{mn}
\,\circbulgamma^{\,e}_{\!hn}=2\,\bulletH^a_h\,\circH^e_b,
\mytag{8.32}\\
&\hskip -2em
\sum^3_{m=0}\sum^3_{n=0}\circbulgamma^{\,a}_{\!bm}\,g^{mn}
\,\bulcircgamma^{\,e}_{\!hn}=2\,\circH^a_h\,\bulletH^e_b.
\mytag{8.33}
\endalign
$$
The identities \mythetag{8.32} and \mythetag{8.33} are analogous 
to \mythetag{8.17} and \mythetag{8.18}. Before applying them to 
\mythetag{8.31} we need to set $e=i$ and $h=a$ in them and then
sum over the index $a$. As a result we transform them to the following
two identities:
$$
\xalignat 2
&\sum^4_{a=1}\sum^3_{m=0}\sum^3_{n=0}\bulcircgamma^{\,a}_{\!bm}
\,g^{mn}\,\circbulgamma^{\,i}_{\!an}=4\,\circH^i_b,
&&\sum^4_{a=1}\sum^3_{m=0}\sum^3_{n=0}\circbulgamma^{\,a}_{\!bm}
\,g^{mn}\,\bulcircgamma^{\,i}_{\!an}=4\,\bulletH^i_b.
\endxalignat
$$
Now, applying these identities to the second line in \mythetag{8.31},
we derive
$$
\gathered
\Alpha^i_{kj}=\sum^4_{a=1}\sum^4_{b=1}
\frac{L_{\boldsymbol\Upsilon_{\!k}}(\bulletd_{\kern 0.5pt ab})
\,\bulletd^{\kern 0.5pt b\kern 0.5pt a}\,\circH^i_j
+L_{\boldsymbol\Upsilon_{\!k}}(\circd_{\kern 0.5pt ab})\,
\circd^{\kern 0.5pt b\kern 0.5pt a}\,\bulletH^i_j}{4}\,+\\
+\sum^3_{m=0}\sum^3_{n=0}\sum^4_{a=1}
\frac{L_{\boldsymbol\Upsilon_{\!k}}(\bulcircgamma^{\,a}_{\!jm})
\,g^{mn}\,\circbulgamma^{\,i}_{\!an}+
L_{\boldsymbol\Upsilon_{\!k}}(\circbulgamma^{\,a}_{\!jm})
\,g^{mn}\,\bulcircgamma^{\,i}_{\!an}}{4}\,-\\
-\sum^3_{m=0}\sum^3_{n=0}\sum^3_{r=0}\sum^4_{a=1}
\frac{\Gamma^r_{km}\,\bulcircgamma^{\,a}_{\!jr}\,g^{mn}
\,\circbulgamma^{\,i}_{\!an}+\Gamma^r_{km}\,\circbulgamma^{\,a}_{\!jr}
\,g^{mn}\,\bulcircgamma^{\,i}_{\!an}}{4}.
\endgathered\qquad
\mytag{8.34}
$$
A metric connection is a real connection. According to
\mythetag{7.11}, the $\bar{\Alpha}$-components of a metric connection 
are obtained from $\Alpha^i_{kj}$ by complex conjugation:
$$
\hskip -2em
\bar{\Alpha}\vphantom{\Alpha}^i_{kj}=\overline{\Alpha^i_{kj}}.
\mytag{8.35}
$$
\mytheorem{8.2} For any skew-symmetric real spin-tensorial field $\bold T$
of the type $(0,0|0,0|1,2)$ there is a unique metric connection $(\Gamma,
\Alpha,\bar{\Alpha})$ in $DM$ with the torsion tensor $\bold T$. Its
components in an arbitrary frame pair $(U,\,\boldsymbol\Upsilon_0,
\,\boldsymbol\Upsilon_1,\,\boldsymbol\Upsilon_2,\,\boldsymbol\Upsilon_3)$
and $(U,\,\boldsymbol\Psi_1,\,\boldsymbol\Psi_2,\,\boldsymbol\Psi_3,
\,\boldsymbol\Psi_4)$ are given by the formulas \mythetag{4.34},
\mythetag{8.34}, and \mythetag{8.35}.
\endproclaim
\mycorollary{8.1} The components of a real metric connection $(\Gamma,
\Alpha,\bar{\Alpha})$ with zero torsion $\bold T=0$ for the bundle of
Dirac spinors $SM$ are given by the explicit formulas \mythetag{4.35},
\mythetag{8.34}, and \mythetag{8.35}.
\endproclaim
    The theorem~\mythetheorem{8.2} and the corollary~\mythecorollary{8.1}
are analogous to the theorem~\mythetheorem{4.3} and its 
corollary~\mythecorollary{4.2} in the case of chiral spinors. As for the
relation of metric connections of Dirac and chiral spinors, it is described
by the theorem~\mythetheorem{8.1}\par

\enddocument
\end